\documentclass[english,11pt]{article}

\usepackage[german,french,english]{babel}
\usepackage[cp850]{inputenc}
\usepackage{latexsym,graphicx, fancybox}
\usepackage{amssymb, amsmath, amsfonts}
\usepackage{color}
\usepackage{latexsym}

\font\tenmath=msbm10 scaled 1200
\font\sevenmath=msbm7 scaled 1200
\font\fivemath=msbm5 scaled 1200 

\newfam\mathfam \textfont\mathfam=\tenmath
\scriptfont\mathfam=\sevenmath \scriptscriptfont\mathfam=\fivemath

\def\R{{\mathbb R}}
\def\N{{\mathbb N}}
\def\E{{\mathbb E}}

\def\P{{\mathbb P}}

\def\Q{{\mathbb Q}}

\def\F{{\cal F}}

\newtheorem{Theorem}{Theorem}[section]
\newtheorem{Definition}[Theorem]{Definition}
\newtheorem{Proposition}[Theorem]{Proposition}

\newtheorem{Lemma}[Theorem]{Lemma}
\newtheorem{Corollary}[Theorem]{Corollary}
\newtheorem{Remark}[Theorem]{Remark}

\newfam\mathfam \textfont\mathfam=\tenmath
\scriptfont\mathfam=\sevenmath \scriptscriptfont\mathfam=\fivemath

\def \^#1{\if#1i{\accent"5E\i}\else{\accent"5E#1}\fi}

\def \cqfd{\quad_\Box}

\begin{document}
\selectlanguage{english}
\title{\bf Convex ordering for stochastic Volterra equations and  their  Euler schemes}
 
\author{ 
{\sc Benjamin Jourdain} \thanks{CERMICS, \'Ecole des Ponts, INRIA, Marne-la-Vall\'ee, France. E-mail: {\tt   benjamin.jourdain@enpc.fr}}~\footnotemark[3]
\and   
{\sc  Gilles Pag\`es} \thanks{Laboratoire de Probabilit\'es, Statistique et Mod\'elisation, UMR~8001, Sorbonne Universit\'e, case 158, 4, pl. Jussieu, F-75252 Paris Cedex 5, France. E-mail: {\tt  gilles.pages@sorbonne-universite.fr}}~\thanks{This research
benefited from the support of the ``Chaire Risques Financiers'', Fondation du Risque.}  }
\maketitle 
\renewcommand{\abstractname}{Abstract}
\begin{abstract}
  In this paper, we are interested in comparing solutions to stochastic Volterra equations for the convex order on the space of continuous $\R^d$-valued paths and for the monotonic convex order when $d=1$. Even if in general these solutions are neither semi-martingales nor Markov processes, we are able to exhibit conditions on their coefficients enabling the comparison. Our approach consists in first comparing their Euler schemes and then taking the limit as the time step vanishes. We consider two types of Euler schemes depending on the way the Volterra kernels are discretized. The conditions ensuring the comparison are slightly weaker for the first scheme than for the second one and this is the other way round for convergence. Moreover, we extend the integrability needed on the starting values in the existence and convergence results in the literature to be able to only assume finite first order moments, which is the natural framework for convex ordering.
\end{abstract}

\section*{Introduction}

The stochastic version of Volterra equations \begin{equation}
   X_t=X_0+\int_0^tK(t,s)b(s,X_s)ds+\int_0^tK(t,s)\sigma(s,X_s)dW_s\label{volsinglK}
\end{equation} goes back to the 1980's in settings corresponding to smooth kernels (see~\cite{Protter1985}). Because of the dependence of the kernel $K(t,s)$ on its first argument, the solution is in general neither a semi-martingale nor a Markov process. Rather than addressing singular kernels, more attention was paid to other extensions like anticipative integrands using Skorokhod integration in connection with Malliavin calculus (see~\cite{PardouxP1990}).
The more demanding case of singular kernels exploding on the diagonal (i.e. $K(t,s)\to +\infty$ as $s\to t$) appears in the  literature in the  early 2000's (see~\cite{CoutinD2001}) motivated, among other applications, by the revival of interest for the fractional Brownian motion as a way to model long memory  in financial mathematics and econometrics (see e.g.~\cite{ComteR1998}). More recently, these stochastic Volterra equations with singular kernels have been brought back to light in a long series of papers devoted to the  modeling of stochastic volatility in Finance which started with the empirical observation in~\cite{GatheralJR2018} that volatility paths have low $H$-H\"older regularity exponent ($H\simeq 0.1$).  Volterra equations with kernel $K(t,s)= (t-s)^{H-\frac 12}$ are a convenient model to reproduce this empirical feature since the H\"older exponent of their paths is $H$. Moreover, such Volterra equations with singular kernels  appear as the limiting dynamics when modeling the order book by ``nearly unstable'' Hawkes processes (see~\cite{JaissonR2016}).
 
 The resulting Volterra equation can be seen as the ``rough'' counterpart of the classical Cox-Ingersoll-Ross  (CIR)  process modeling the volatility in the classical/regular Heston stochastic volatility model. The rough Heston model obtained by substituting the Volterra process to the original CIR process can be used to price and hedge equity derivative products (see e.g.~\cite{BayerFG2016},~\cite{ElEuchR2018}).

 In the case when, more generally, $K(t,s)$ is a function of $t-s$, the stochastic Volterra equation is a stochastic convolution equation which can be studied using tools originally developped for the study of deterministic convolution equations, in particular resolvents (see \cite{AJLP2019,AJCLP2021}).

\smallskip
On the other hand, convex ordering of random vectors has been developed for long as a refined way to quantify  or compare risk in Econometrics and in Quantitative Finance.  It can be defined between two distributions $\mu$ and $\nu$ on $\R^d$ (having a finite first moment) by
\[
\mu\le_{cvx} \nu \quad\mbox{ if } \quad\forall \, \varphi:\R^d \to \R\mbox{ convex}, \quad \int\varphi(x)\mu(dx) \le \int \varphi(y)\nu(dy).
\]
This implies that both distributions
have the same expectation. Two $\R^d$-valued random vectors $X$ and $Y$ are ordered for the convex order if their distributions are. This approach had more recently applications to determine the sign and to  compare the sensitivities of derivatives products in an extensive (in fact functional) sense (see~\cite{Ruschendorf1},~\cite{Pag2016},~\cite{LiuPag2020},~\cite{LiuPag}). It is also an important tool to define and identify the {\em tracking error} in the  hedging process of derivative products when  the volatility dynamics is misspecified. The aim of this paper is to extend some functional convex ordering results established for Brownian, jumpy or McKean-Vlasov diffusion processes (see~\cite{LiuPag2020},~\cite{Pag2016}) in the above cited papers (and called {\em functional convex order or ordering})  to the solution of  Volterra equations. A typical result, consequence of Theorems~\ref{thm:main} and~\ref{thm:ci}  is the following: consider two solutions to the  following equations (defined on a probability space $(\Omega,{\cal A}, \P)$):
\[
X_t = X_0+\int_0^t (t-s)^{H-\frac 12}\sigma(X_s)dW_s \quad\mbox{ and }\quad Y_t =Y_0 +\int_0^t (t-s)^{H-\frac12}\tilde\sigma(Y_s)dW_s, \; t\!\in [0,T]
\]
($H\!\in (0,1]$) where  $W$ is standard Brownian motion and $\sigma, \tilde\sigma: \R\to \R$ are Lipschitz continuous.  If $X_0,\, Y_0\!\in L^1(\P)$ with $X_0\le_{cvx} Y_0$ and
\[
|\sigma| \; \mbox{ or }\;|\tilde\sigma| \mbox{ is convex $\qquad$ and } \qquad|\sigma| \le |\tilde\sigma|,
\] 
then for every lower semi-continuous (l.s.c.) convex functional $F: \big({\cal C}([0,T], \R), \|\cdot\|_{\sup}\big)\to \R$, 
$$
\E\, F\big((X_t)_{t\in [0,T]}\big)\le \E\, F((Y_t)_{t\in [0,T]}\big)
$$
{(where both expectations have a sense, see Section~\ref{subsec:1.2}).} Moreover if $|\sigma|$ is convex, then $x\mapsto \E\, F\big((X^x_t)_{t\in [0,T]}\big)$  is convex (where $X^x_0= x$). The extension of these results to the multimensional setting where $\sigma,\tilde\sigma:\R^d\to\R^{d\times q}$, $W$ is a $q$-dimensional Brownian motion and  $X_0,Y_0$ are $\R^d$-valued random vectors holds under appropriate generalizations of the condition $|\sigma|\le|\tilde\sigma|$ and the convexity of either $\sigma$ or $\tilde \sigma$.

\smallskip
The main strategy adopted in this paper to establish convex ordering properties on Volterra processes is two fold: first we establish these properties on a discrete time Euler scheme and then we transfer them to the continuous time Volterra process using a convergence theorem~---~here pathwise in $L^p$~--~of the time discretization scheme to the process  when the time step goes to $0$. This roadmap has already been used in a Markovian framework by various authors (see~\cite{Ruschendorf1, Ruschendorf2, Ruschendorf3,Pages2016} for diffusions and~\cite{LiuPag} for McKean-Vlasov diffusions). Such an approach is important for applications: an approximation scheme which can be easily simulated that preserves (and transfers)  convexity and convex ordering is a key asset since it avoids to create artificial arbitrage  opportunities in hedging portfolios of derivative products computed by discretization and more generally permits a better monitoring of market risks.  

\smallskip
Compared to the aforementioned  diffusions models, several  difficulties arise to implement this strategy and establish such a kind of result under natural $L^1(\P)$ integrability assumptions on the starting values when dealing with Volterra equations, especially with singular kernels like in the above example. First, of course, since the solution is not a semi-martingale  nor a Markov process,  powerful tools like stochastic calculus or Markov properties cannot be used/called upon to establish convex ordering. The other more hidden difficulty  follows from existing results about  integrability and pathwise regularity of the solutions to the Volterra equation. Thus following~\cite{ZhangXi2010}, the existence and uniqueness of a pathwise continuous and integrable solution is established under the stringent assumption that {\em $X_0$ and $Y_0$ lie in all $L^p(\P)$-spaces} i.e. $X_0, \,Y_0\!\in \cap_{p>0}L^p(\P)$ (under an additional appropriate integrability and regularity  condition on  the kernels which are satisfied by $K(t,s)= (t-s)^{H-\frac 12}$). We show in Theorem~\ref{prop:exunvolt} that this result can be extended to the natural case where $X_0$ and $Y_0$ lie in some $L^p(\P)$ for some $p\!\in (0, +\infty)$ by proving a representation formula of solutions as functionals of the Brownian path and the starting value (see Appendix~\ref{app:B'}). 

According to \cite{ZhangEuler} or~\cite{RiTaYa2020}, the same restrictions in terms of integrability come out for the convergence of the $K$-discrete Euler scheme~\eqref{eq:Eulergen1} of these equations where the second variable in the kernels is also chosen in the discretization grid. Moreover, the assumptions on the kernels $K(t,s)$ which ensure the $L^p(\P)$-convergence (in $\sup$-norm) of the continuous time extension of this Euler scheme toward the solution to the Volterra equation (see~Theorem~\ref{thm:Eulercvgce1}) are more stringent  than those ensuring the existence of a strong solution to the equation itself (see Theorem~\ref{prop:exunvolt}). To overcome this problem we will also consider  and analyze the convergence of the $K$-integrated Euler scheme~\eqref{eq:Euler2} which can also be simulated in many cases of interest (depending on the form of the kernel). We will show in Theorem~\ref{thm:Eulercvgce2} that, if $X_0\!\in L^p(\P)$ for some $p>0$, then it converges in $L^p(\P)$ for the  $\sup$-norm under  the same assumptions  (up to time regularity of the coefficients of the equation) as those ensuring  strong existence-uniqueness for the Volterra process. Propagating convexity and convex ordering through such a scheme appears as the adequate  tool to establish quite  general  results on convex ordering for the Volterra equation, especially as concerns the  assumptions on the kernels. The comparison betwen the two schemes remains balanced in the sense that, if relying on the $K$-discrete time scheme requires slightly more demanding assumptions on the kernels  to get  convergence,  this scheme is sometimes  easier and more stable  to simulate for very small steps. The integrability restriction on $X_0$ needed in~\cite{RiTaYa2020}  is also relaxed for this scheme   in Theorem~\ref{thm:Eulercvgce1}. Of course both schemes have a higher complexity in terms of simulation than their counterparts for regular diffusions due to the non-markovianity of the Volterra process (we refer to~\cite{AlfonsiK2021} for alternative approaches devised to improve the simulation of such  processes).

The first section is devoted to the statement of the main results: existence and uniqueness for the stochastic Volterra process, convergence to this process of the discretization schemes as the time step goes to $0$, convex ordering between two Volterra processes and monotonic convex ordering in dimension $d=1$. It also gives an application to VIX options in  the quadratic rough Heston model introduced in~\cite{GaJuRo2020}. In Section 2, we state and prove the convex ordering for the $K$-discrete and $K$-integrated Euler schemes. Section 3 addresses the monotonic convex ordering for the discretization schemes in dimension $d=1$. In Section 4, we deduce the ordering for the Volterra processes by implementing the second step of our strategy.

 \medskip
 \noindent {\sc Notation} $\bullet$ We denote $\R^{d\times q}$ the set of matrices with $d$ rows and $q$ columns. 

\smallskip
\noindent $\bullet$ We denote  $A^*$ the transpose matrix of $A$ (idem for vectors).

\smallskip
\noindent $\bullet$ We denote by ${\cal S}^+(d)= \big\{S\!\in \R^{d\times d}, \mbox{ symmetric, positive semi-definite}\big\}$ and ${\cal O}(d)= \big\{ P\!\in \R^{d\times d}:  PP^* =I_d\big\}$.

\smallskip
\noindent $\bullet$  For $S$, $T\!\in {\cal S}^+(d)$, $S\le T$ if $T-S\!\in {\cal S}^+(d)$.

\smallskip
\noindent $\bullet$  For $S\!\in {\cal S}^+(d)$, we denote by $\sqrt{S}$ the unique element of ${\cal S}^+(d)$ such that $\sqrt{S}\sqrt{S}=S$.

\smallskip 
\noindent $\bullet$  For $A\!\in \R^{d_1\times d_1}$ and $B\!\in \R^{d_2\times d_2}$, we denote  by $A\otimes B\in\R^{(d_1\times d_2)\times(d_1\times d_2)} $ the Kronecker product  of $A$ and $B$ defined by $$A\otimes B = \left(\begin{array}{cccc}A_{11} B &A_{12}B&\hdots &A_{1d_1}B\\A_{21} B &A_{22}B&\hdots &A_{2d_1}B\\\vdots &\vdots&\vdots&\vdots \\A_{d_11}B &A_{d_12}B &\hdots &A_{d_1d_1}B
                             \end{array}\right).
$$

\smallskip 
\noindent $\bullet$ For a function $f: E\to \R$, $\displaystyle \|f\|_{\sup}= \sup_{x\in E}|f(x)|$. 
 
 \smallskip 
\noindent $\bullet$ $\lambda_d$ denotes the Lebesgue measure on $(\R^d, {\cal B}or(\R^d))$, often denoted $\lambda$ when $d=1$.

 \smallskip 
\noindent $\bullet$ ${\cal C}([0,T], \R^d)$ (resp. ${\cal C}_0([0,T], \R^d)$) denotes the set of continuous functions (resp. null at $0$) from $[0,T]$ to $\R^d$. {${\cal B}or({\cal C}_{d})$ (resp. ${\cal B}or({\cal C}_{d,0})$) denotes its  Borel $\sigma$-field induced by the $\sup$-norm topology}.  

\smallskip 
\noindent $\bullet$ $\perp\!\!\!\perp$ stands for independence of random variables, vectors or processes.

\smallskip 
\noindent
$\bullet$ {$L_{\R^d}^0(\P)$ or simply $L^0(\P)$ denote the set of  $\R^d$-valued random vectors  defined on a probability space $(\Omega, {\cal A}, \P)$. For $p\in(0,+\infty)$, $L_{\R^d}^p(\P)$ or simply $L^p(\P)$ denote the set of  $\R^d$-valued random vectors $X$ defined on a probability space $(\Omega, {\cal A}, \P)$ such that $\|X\|_p:=\E^{1/p}|X|^p<+\infty$.}

\section{Convex ordering of Volterra {processes}: main result}
 \subsection{Volterra processes and their Euler schemes}

We are interested in the stochastic Volterra equation
 \begin{equation}\label{eq:Volterra}
   X_t=X_0+\int_0^t K_1(t,s) b(s,X_s)ds+\int_0^t K_2(t,s)\sigma(s,X_s)dW_s,\quad t\in[0,T],
 \end{equation}
 where $b:[0,T]\times\R^d\to \R^d$, $\sigma:[0,T]\times\R^d\to\R^{d\times q}$, $K_i: \big\{ (t,s): 0\le s<t\le T\big\}\to \R_+,\;i\in\{1,2\}$ are measurable and $(W_t)_{t\in[0,T]}$ is a standard $q$-dimensional Brownian motion independent from the $\R^d$-valued random vector $X_0$. Let $(\F_t)_{t\in [0,T]}$ be a filtration such that $X_0\!\in {\cal F}_0$ and $W$ is an ${\cal F}_t$-Brownian motion.
 The following existence and uniqueness result   for~\eqref{eq:Volterra} is an improved version of~\cite[Theorem 3.1 and Theorem 3.3]{ZhangXi2010}, at least for Volterra equations having the above form.

 \begin{Theorem}\label{prop:exunvolt}   
Assume that the kernels $K_i$, $i=1,2$, satisfy the integrability assumption
 \begin{equation}\label{eq:KisL^2}
\big({\cal K }^{int}_\beta\big)\hskip1,5cm 
\sup_{t\in [0,T]}\int_0^t\left(K_1(t,s)^
{{\frac{2\beta}{\beta+1}}}
  +K_2(t,s)^{2\beta}\right)ds<+\infty \hskip 1,5cm 
\end{equation}
for some $\beta>1$ and the continuity assumption
\begin{align}\nonumber
({\cal K}^{cont}_{\theta}) \;\;  
\exists\, \kappa< +\infty,&\;\forall \,\delta{\in (0, T)},\\
\label{eq:contK}&\eta(\delta):=\max_{i=1,2} \sup_{t\in [0,T]} \left[\int_0^t |K_i(\big(t+\delta)\wedge T,s\big)-K_i(t,s)|^ids \right]^{\frac 1i} \le  \kappa\,\delta^{\theta}
\end{align}
for  some $\theta\in (0,1]$.  Assume that the functions $\R^d\ni x\mapsto b(t,x)$ and $\R^d\ni x\mapsto \sigma(t,x)$ are Lipschitz with Lipschitz  coefficient uniform in $t\in[0,T]$, i.e.
\[
\exists \, C_{_T}= C_{b,\sigma,T}\; \mbox{ such that }\; \forall\, t\!\in [0,T], \; \forall\, x,\, y\!\in \R^d,\quad |b(t,x)-b(t,y)|+\|\sigma(t,x)-\sigma(t,y)\|\le C_{_T} |x-y|
\] 
and $\displaystyle\sup_{t\in[0,T]}(|b(t,0)|+\|\sigma(t,0)\|)<+\infty$. Finally assume that 
$X_0\!\in  L^p(\P)$ for some $p\!\in (0, +\infty)$.

Then the Volterra equation~\eqref{eq:Volterra} admits, up to a $\P$-indistinguishability, a unique  $({\cal F}_t)$-adapted solution $X=(X_t)_{t\in [0,T]}$, pathwise continuous,  in the sense that, $\P$-$a.s.$,
\[
\forall\, t\!\in [0,T],  \quad X_t=X_0+\int_0^t K_1(t,s) b(s,X_s)ds+\int_0^t K_2(t,s)\sigma(s,X_s)dW_s.
\]
 This solution satisfies
\begin{equation}\label{eq:Lpincrements}
\forall\, s,\, t\!\in [0,T],\quad \|X_t-X_s\|_p \le C_{p,T} {(1+\|X_0\|_p)}|t-s|
^{\theta\wedge \frac{\beta-1}{2\beta}}.
\end{equation}
Moreover,  

\begin{equation}\label{eq:Holderpaths}
\forall\, a\!\in \big(0,\theta\wedge \frac{\beta-1}{2\beta}\big), \qquad\left \| \sup_{s\neq t\in [0,T]}\frac{|X_t-X_s|}{|t-s|^{a}}\right\|_p  <C_{a,p,T}  (1+\|X_0\|_p)
\end{equation}
for some positive real constant $C_{a,p,T}=C_{a,b,\sigma, K_1,K_2, \theta,p,T} $. In particular
\begin{equation}\label{eq:supLpbound}
\Big\| \sup_{t\in [0,T]} |X_t|\Big\|_p \le C'_{a,p,T}(1+\|X_0\|_p).
\end{equation}
 Finally, if the condition
\begin{equation}\label{eq:contKtilde}
(\widehat {\cal K}^{cont}_{\widehat \theta})\;\;\exists\,\widehat\kappa< +\infty,\;\forall\delta\!\in (0,T],\; \widehat \eta(\delta) := \max_{i=1,2} \sup_{t\in [0,T]} \left[\int_{(t-\delta)^+}^t \hskip-0,25cm K_i\big(t,u\big)^i du\right]^{1/i}\le \widehat \kappa \,\delta^{\,\widehat \theta}
\end{equation}
is satisfied for some $\widehat\theta\in (0,1]$, then one can replace $\frac{\beta-1}{2\beta}$ by $\widehat \theta$ in~\eqref{eq:Lpincrements} and~\eqref{eq:Holderpaths}.
\end{Theorem}

\noindent {\bf Comments \& Remarks}. $\bullet$ In~\cite[Theorem 3.3]{ZhangXi2010}, the starting value $X_0$  of the equation is assumed to lie in $\cap_{r>0}L^r(\P)$. In fact a careful reading shows that the proof establishes a slightly more general result than announced: if $X_0\!\in L^p(\P)$ for $p>p_{eu}:=\frac{1}{\theta} \vee \frac{2\beta}{\beta-1}$, then existence and uniqueness (up to $\P$-indistinguishability) of an $(\F_t)$-adapted pathwise $a$-H\"older continuous solution for some small enough $a>0$ holds true.  No bound on $a$ is provided.  The upper-bounds~\eqref{eq:Lpincrements},~\eqref{eq:Holderpaths} and~\eqref{eq:supLpbound}  are established for $p\ge 2$ and their dependence  in $\|X_0\|_p$ is not (always) specified. 

\smallskip
\noindent $\bullet$  Hence the above  theorem appears~--~in our  setting~--~as an extension of~\cite[Theorem 3.3]{ZhangXi2010} in terms of existence and uniqueness and in terms of moment control since we simply assume $X_0\!\in L^p(\P)$  for a fixed $p\!\in (0, +\infty)$. The extended existence and uniqueness is proved in Appendix \ref{app:eu}. Moreover, the moment  control in~\eqref{eq:Lpincrements},~\eqref{eq:Holderpaths} and~\eqref{eq:supLpbound} with respect to the starting  value $X_0$ are mostly new. These controls combined with the representation of the solution as a functional $F(X_0,W)$  of the starting value $X_0$ and the Brownian motion $W$ (like for diffusions) is the key to our extension achieved through Lemma~\ref{lem:gap filled}. This extension is proved in Appendix~\ref{app:A''}. Its main interest is of course to include the case $p=1$ in view of   our convex ordering and convexity results.

\smallskip
\noindent  $\bullet$ Assumptions $({\cal K}^{int}_{\beta})$ and  $(\widehat {\cal K}^{cont}_{\widehat \theta})$ are very close  in the singular case, e.g. when $K_i(t,s)= \varphi_i(t-s)$, $0\le s<t \le T$ with $\varphi_i$ decreasing from $+\infty$ at  $0$.

\smallskip
If $({\cal K}^{int}_{\beta})$ is in force, then $(\widehat {\cal K}^{cont}_{\widehat \theta})$ holds with $\widehat \theta =  \frac{\beta-1}{2\beta}$ since H\"older's inequality implies, for $\delta \!\in [0,T]$,  
\[
\int_{(t-\delta)^+}^tK_1(t,u) du \le \Big(\int_0^t K_1(t,u)^{\frac{2\beta}{\beta+1}}du \Big)^{\frac{\beta+1}{2\beta}}  (\delta\wedge t)^{\frac{\beta-1}{2\beta}}\le C_{1,\beta,T}\,\delta^{\frac{\beta-1}{2\beta}} 
\]
and 
\[
\Big(\int_{(t-\delta)^+}^tK_2(t,u)^2 du\Big)^{1/2}  \le \Big(\int_0^t K_2(t,u)^{2\beta}du \Big)^{\frac{1}{2\beta}}  (\delta\wedge t)^{\frac{\beta-1}{2\beta}}\le C_{2,\beta,T}\,\delta^{\frac{\beta-1}{2\beta}} .
\]
Consequently, in that case, we may always assume that $\widehat \theta\ge \frac{\beta-1}{2\beta}$. 

There is no converse  in general. However if, for every $t\!\in[0,T]$, 
\[
[0,t) \ni s\mapsto K_i(t,s) \mbox{ is non-decreasing}, \; i=1,2,
\]
(which includes usual singular kernels) then $(\widehat {\cal K}^{cont}_{\widehat \theta})$ implies $({\cal K}^{int}_{\beta})$ for every $\beta$ such that $\frac{\beta-1}{2\beta}<\widehat \theta$. In that case $\widehat \theta$ appears as the supremum (not attained) of  $\frac{\beta-1}{2\beta}$ over all $\beta$ is such that $({\cal K}^{int}_{\beta})$ holds. This easily follows  from the fact that, taking advantage of the monotonicity of the kernels $K_i$, one has under $(\widehat {\cal K}^{cont}_{\widehat \theta})$, for $0\le s<t\le T$,
$$
(t-s)^{1/i}K_i(t,s)\le\left(\int_s^t K_i(t,u)^idu\right)^{1/i}\le \widehat\kappa(t-s)^{\widehat\theta}\mbox{ so that }K_i(t,s)\le \widehat\kappa(t-s)^{\widehat\theta-1/i}
$$and $({\cal K}^{int}_{\beta})$  holds as soon as $\frac{2\beta}{\beta+1}(\widehat \theta- 1) >- 1$ and $2\beta(\widehat \theta-\frac12)>-1$ which reads for both kernels $\frac{\beta-1}{2\beta}< \widehat \theta$.

\medskip
A significant difference between regular diffusion processes and Volterra processes from a technical viewpoint  comes from the presence of the kernels which introduces some memory in the dynamics of the process, depriving us of the Markov property and usual tools of stochastic calculus.

Moreover, on the more practical  level of simulation, it  naturally  suggests two ways to devise a discrete time Euler scheme associated to a regular mesh $(t^n_k=\frac{kT}{n})_{0\le k\le n}$ of $[0,T]$ with step $h= \frac Tn$.

\medskip 
\noindent $\rhd$ {\em $K$-discrete Euler scheme}: this Euler scheme is defined inductively by 
\begin{align}\label{eq:Eulergen1}
   \bar X_{t^n_k}=X_0+ \sum_{\ell=1}^{k} \Big(K_1(t^n_k,t^n_{\ell-1})b(t^n_{\ell-1},\bar X_{t^n_{\ell-1}})\tfrac{T}{n} +K_2(t^n_k,t^n_{\ell-1})\sigma(t^n_{\ell-1},\bar X_{t^n_{\ell-1}})&(W_{t^n_\ell}-W_{t^n_{\ell-1}})\Big),\\
\nonumber    &\qquad k=0,\ldots,n.
 \end{align}
 The notation is slightly ambiguous, 
 but we choose not to alleviate this ambiguity by denoting $\bar X^n$ in place of $\bar X$, since, on the one hand, we only consider the Euler scheme with $n$ steps throughout the paper and, on the other hand, for $x_0\in\R^d$, we will at some places use the notation $\bar X^{x_0}$ to emphasize the dependence of the Euler scheme on the initial condition $X_0=x_0$.
 
 Note that simulation of this scheme only requires to compute pointwise values of the kernels $K_i$ and simulate standard Brownian increments.
 
\smallskip  In view of producing a priori error bounds/convergence rate  for the approximation of the stochastic Volterra process by its Euler scheme~\eqref{eq:Euler1}, we naturally  extend it into its continuous time (or ``genuine'' version) by setting
 \begin{equation}\label{eq:Euler1}
\bar X_t =X_0+\int_0^t K_1(t,\underline s)b(\underline s,\bar X_{\underline s})ds +\int_0^t K_2(t,\underline s)\sigma(\underline s,\bar X_{\underline s})dW_s,\;t\in[0,T],
 \end{equation}
where  $\underline s = \underline s_n=  t^n_k := \frac{kT}{n}$ when $s\!\in [t^n_k, t^n_{k+1})$. It reads 
{``in extension'' when $t\!\in [t^n_k,t^n_{k+1})$,}
 \begin{align}
 \label{eq:Eulergen2bis} \bar X_{t} &=X_0 +  K_1(t,t^n_k) \,b(t^n_{k},\bar X_{t^n_{k}})(t-t^n_k)+\sigma(t^n_{k},\bar X_{t^n_{k}})K_2(t,t^n_k)(W_{t}-W_{t^n_k}) \\
   \nonumber  &\quad + \sum_{\ell=1}^{k} \Big(  K_1(t,t^n_{\ell-1})\,b(t^n_{\ell-1},\bar X_{t^n_{\ell-1}}) \tfrac Tn  +\sigma(t^n_{\ell-1},\bar X_{t^n_{\ell-1}})  K_2(t,t^n_{\ell-1})(W_{t^n_{\ell}}-W_{t^n_{\ell-1}})\, \Big),\;k=0,\ldots,n.
 \end{align}
 \noindent $\rhd$ {\em $K$-integrated  discrete time Euler scheme}:  This scheme reads as follows
 \begin{align}\label{eq:Euler2}
   \bar X_{t^n_k}=X_0+ \sum_{\ell=1}^{k} \Big(\int_{t^n_{\ell-1}}^{t^n_{\ell}} \hskip-0,25cm K_1(t^n_k,s)ds \,b(t^n_{\ell-1},\bar X_{t^n_{\ell-1}})  +\sigma(t^n_{\ell-1},\bar X_{t^n_{\ell-1}})\int_{t^n_{\ell-1}}^{t^n_{\ell}}\hskip-0,25cm &K_2(t^n_k,s)dW_s\, \Big),\\
\nonumber    &\qquad k=0,\ldots,n.
 \end{align}
When a closed or a semi-closed form is available for the vectors $\displaystyle\Big[\int_{t^n_{\ell-1}}^{t^n_{\ell}}K_1(t^n_k,s)ds\Big]_{1\le \ell\le k}$, $k=1,\ldots,n$ and the covariance matrices $\displaystyle \Big [\int_{t^n_{\ell-1}}^{t^n_{\ell}}K_2(t^n_k,s)K_2(t^n_{k'},s)ds \Big]_{\ell\le k, k'\le n}$, $\ell=1,\ldots,n$, then the above $K$-integrated Euler  below becomes  simulable (see practitioner's corner further on).

Especially when  he kernels $K_i$ are singular ( i.e.  $\lim_{s\to t-}K_i(t,s) =+\infty$), the  main asset of this second scheme is that the terms corresponding to $\ell=k$ are computed/simulated much more accurately than their above counterparts in the $K$-discrete scheme (see~the practitioner's corner at the end of the current section  in the case $K_i(t,s)=(t-s)^{\alpha_i}$). This may have a positive significant impact in practical simulations when the step $\frac Tn$ is  a bit coarse.

This variant also admits a continuous time version given with the same notations by 
 \begin{equation}\label{eq:Eulergen2}
\bar X_t =X_0+\int_0^t K_1(t,s)b(\underline s,\bar X_{\underline s})ds +\int_0^t K_2(t,s)\sigma(\underline s,\bar X_{\underline s})dW_s,\; t \!\in[0,T],
 \end{equation}
 {which also reads ``in extension'' when $t\!\in [t^n_k,t^n_{k+1})$,}
 \begin{align}
  \nonumber  \bar X_{t} &=X_0 + \int_{t^n_{k}}^{t}  K_1(t,s)ds \,b(t^n_{k},\bar X_{t^n_{k}})  +\sigma(t^n_{k},\bar X_{t^n_{k}})\int_{t^n_{k}}^{t}   K_2(t,s)dW_s \\
\label{eq:Eulergen2bis}   &\quad + \sum_{\ell=1}^{k} \Big(\int_{t^n_{\ell-1}}^{t^n_{\ell}} \hskip-0,25cm K_1(t,s)ds \,b(t^n_{\ell-1},\bar X_{t^n_{\ell-1}})  +\sigma(t^n_{\ell-1},\bar X_{t^n_{\ell-1}})\int_{t^n_{\ell-1}}^{t^n_{\ell}} K_2(t,s)dW_s\, \Big),\;k=0,\ldots,n.
 \end{align}
 
 
 \smallskip
 \noindent
 $\rhd$ {\em Lack of Markovianity}. Because of the lack of Markovianity, for both variants, $\bar X_{t^n_k}$ is in general not a function of $(\bar X_{t^n_{k-1}},W_{t^n_{k}}-W_{t^n_{k-1}})$. Nevertheless, $\bar X_{t^n_k}$ can be computed in a unique way from the values of $(\bar X_{0},\cdots,\bar X_{t^n_{k-1}})$ and either the  vector of Brownian increments $\left(W_{t^n_{1}}-W_{t^n_0},\cdots,W_{t^n_{k}}-W_{t^n_{k-1}}\right)$ for~\eqref{eq:Eulergen1} or the Gaussian vector $\left(\int_{t_{\ell-1}}^{t_{\ell}}K_2(t^n_k,s)dW_s\right)_{\ell= 1,\ldots,k}$ for~\eqref{eq:Euler2}  so that these Euler schemes are well defined by induction. Note that the complexity for the computation of $(\bar X_{t^n_k})_{k=0,\ldots,n}$ is ${\mathcal O}(n^2)$ and that for the $K$-integrated scheme, one should generate the independent Gaussian vectors $Y^{(\ell)}=\left(\int_{t_{\ell-1}}^{t_{\ell}}K_2(t^n_k,s)dW_s\right)_{k=\ell,\ldots,n}$, $\ell=1, \ldots, n$,  (see practitioner's corner further on).

\medskip
We start  by a theorem of convergence for the $K$-integrated Euler scheme since less assumptions are needed than for the convergence of the $K$-discrete Euler scheme.

 \begin{Theorem}[$K$-integrated Euler scheme]\label{thm:Eulercvgce2}  
 Let $T>0$ and let $p\!\in (0, +\infty)$ being fixed. Assume $b$ and $\sigma$ satisfy the following time-space H\"older-Lipschitz continuity assumption
 \begin{align}
\nonumber ({\cal LH}_{\gamma})\quad  &\exists\, C_{b,\sigma}< +\infty,\; \forall\, s,t\!\in [0,T], \; \forall\, x,\,y \!\in \R^d,\\
\label{eq:HolLipbsig} & |b(t,y)-b(s,x)| +\|\sigma(t,y)-\sigma(s,x)\| \le C_{b,\sigma}\big((1+|x|+|y|)|t-s|^{\gamma} + |x-y|\big)
 \end{align}
 for some $\gamma\in (0,1]$. Assume the kernels $K_i$, $i=1,2$,   satisfy the integrability condition~$({\cal K}^{int}_{\beta})$  (see~\eqref{eq:KisL^2}) for some $\beta>1$ and  the continuity condition $({\cal K}^{cont}_{\theta} )$ (see~\eqref{eq:contK}) for some $\theta\in (0,1]$.  {If $\bar X$ denotes  the  $K$-integrated  Euler scheme~\eqref{eq:Eulergen2} with time step $\frac Tn$, then $\bar X$ has a pathwise continuous modification.} 

 Assume furthermore
{that  $(\widehat {\cal K}^{cont}_{\widehat \theta})$} holds
for some $\widehat\theta\in (0,1]$. 
%
Then there exists a real constant $C= C_{K_1,K_2,\beta, p,b,\sigma,T} \!\in (0, +\infty)$  (not depending on $n$) such that, for every $n\ge 1$,  
 \begin{equation}\label{eq:incrementXbar}
\forall\, s,\, t\!\in [0, T], \qquad 
\big\|\bar X_t -\bar X_s\big\|_p \le C(1+\|X_0\|_p)|t-s|^{\theta\wedge \widehat \theta }
 \end{equation}
and
\begin{equation}\label{eq:Xt-barX_t}
\max_{k=0,\ldots,n}\big\|X_{t_k} -\bar X_{t_k}\big\|_p \le \sup_{t\in [0,T]}\big\|X_t -\bar X_t\big\|_p\le C(1+\|X_0\|_p) \big(\tfrac Tn\big)^{\gamma\wedge \theta \wedge  \widehat\theta}.
\end{equation}
%
%
Moreover, there exists for every $\varepsilon\in(0,1)$  a real constant $C_{\varepsilon}= C_{\varepsilon, K_1,K_2,\beta,\gamma, p,b,\sigma,T} \!\in (0, +\infty)$ such that 
\begin{equation}\label{eq:sup|Xt-barX_t|}
\Big \|\max_{k=0,\ldots,n}\big|X_{t_k} -\bar X_{t_k}\big| \Big\|_p \le \Big\|\sup_{t\in [0,T]}\big|X_t -\bar X_t\big| \Big\|_p\le C_{\varepsilon}(1+\|X_0\|_p)\big(\tfrac Tn\big)^{(\gamma\wedge\theta\wedge \widehat\theta) (1-\varepsilon)}.
 \end{equation}
\end{Theorem}
The proof of this theorem is postponed to Appendices~\ref{app:B} and~\ref{app:E}. 
\begin{Remark}  If $(\widehat {\cal K}^{cont}_{\widehat \theta})$ is not satisfied, the above theorem remains valid by replacing {\em mutatis mutandis} $\widehat \theta$ by $\frac{\beta-1}{2\beta}$ in all three inequalities
\end{Remark}

\begin{Theorem}[$K$-discrete Euler scheme]\label{thm:Eulercvgce1} Let $T>0$ and let $p\!\in (0, +\infty)$ being fixed. 
Assume $\|X_0\|_p<+\infty$, $b$ and $\sigma$ satisfy~$({\cal LH}_{\gamma})$.  Assume that the kernels $K_i$, $i=1,2$, satisfy~$({\cal K}^{int}_{\beta})$ for some $\beta >1$,

\begin{align}
(\underline{\cal K}^{cont}_{\underline\theta}) \; 
&\exists\,\kappa< +\infty,\;\forall \delta\!\in [0,T),\;                      \max_{i=1,2}\sup_{t\in [0,T], \,n\ge 1}
\displaystyle \left( \int_0^t \big|K_i((t+\delta)\wedge T,\underline u_n)-K_i(t,\underline u_n)\big|^idu\right)^{1/i}\le  \kappa\,\delta^{\underline \theta}\label{eq:contKbar}
\end{align} 
for some $\underline\theta\in(0,1]$ and 
\begin{align}
(\underline{\widehat{\cal K}}^{cont}_{\underline{\widehat\theta}}) \; 
&\exists\,\widehat\kappa< +\infty,\;\forall \delta\!\in[0,T/2],\; \max_{i=1,2}\sup_{t\in [0,T], \,n\ge 1}\left( \int_{(t-\delta)^+}^t\ K_i(t,\underline u_n )^idu \right)^{1/i}\le \widehat\kappa \,\delta^{\underline{\widehat\theta}}\label{eq:contKtildebar}
\end{align}
for some $\underline{\widehat\theta}\!\in (0,1]$.

Let $\bar X$ denote  the $K$-discrete time Euler scheme~\eqref{eq:Euler1} with step $\frac Tn$.

\medskip
\noindent $(a)$  Under the above assumptions, $\bar X$ has continuous paths and there exists a real constant $C= C_{K_1,K_2,\beta, p,b,\sigma,T} \!\in (0, +\infty)$  such that, for every $n\ge 1$,  
 \begin{equation}\label{eq:incrementXbar2}
\forall\, s,\, t\!\in [0, T], \qquad 
 \big\|\bar X_t -\bar X_s\big\|_p \le C(1+\|X_0\|_p)|t-s|^{\underline\theta\wedge \underline{\widehat \theta}}.
\end{equation}
\noindent $(b)$ Assume moreover that
\begin{align*}
(\underline{\check{\cal K}}^{cont}_{\underline{\check\theta}})\;\;    \exists\, \check\kappa< +\infty,\;\forall n\ge 1,\;\max_{i=1,2}\sup_{t\in[0,T]}\left( \int_0^t\big|K_i(t,u)-K_i(t,\underline u_n)\big|^idu\right)^{1/i}\le\check\kappa\left(\frac Tn \right)^{\underline{\check\theta}}
\end{align*}
for some $\underline{\check\theta}>0$.
Then there exists a finite constant $C$ such that \begin{equation}\label{eq:Xt-barX_t2}
\max_{k=0,\ldots,n}\big\|X_{t_k} -\bar X_{t_k}\big\|_p \le \sup_{t\in [0,T]}\big\|X_t -\bar X_t\big\|_p\le C(1+\|X_0\|_p) \big(\tfrac Tn\big)^{\gamma\wedge \underline\theta\wedge \underline{\widehat \theta}\wedge \underline{\check \theta}} .
\end{equation}Moreover,  for every $\varepsilon\in(0,1)$,  there exists  a real constant $C_{\varepsilon}= C_{\varepsilon, K_1,K_2,\beta,\gamma, \gamma',p,b,\sigma,T} \!\in (0, +\infty)$ such that 
 \begin{equation}\label{eq:sup|Xt-barX_t|2}
\Big \|\max_{k=0,\ldots,n}\big|X_{t_k} -\bar X_{t_k}\big| \Big\|_p \le \Big\|\sup_{t\in [0,T]}\big|X_t -\bar X_t\big| \Big\|_p\le C_{\varepsilon}(1+\|X_0\|_p)\big(\tfrac Tn\big)^{(\gamma\wedge\underline\theta \wedge \underline{\widehat\theta}\wedge \underline{\check\theta}) (1-\varepsilon)}.
\end{equation}
\end{Theorem}

The main difference with the original reference~\cite{RiTaYa2020} is that we claim that the rate holds for any $p>0$. This gap is filled in Appendix~\ref{app:E} using  a general representation formula of the solutions and a ``splitting''  lemma (Lemma~\ref{lem:gap filled}) both established in Appendix~\ref{app:B'}.

\begin{Remark}\label{rem:hypkdisckint}Note that under $(\underline{\check{\cal K}}^{cont}_{\underline{\check\theta}})$
for some $\underline{\check\theta}>0$, we have $(\underline{\cal K}^{cont}_{\underline\theta}) \Rightarrow({\cal K}^{cont}_{\underline\theta})$ and $(\underline{\widehat {\cal K}}^{cont}_{\underline{\widehat\theta}}) \Rightarrow(\widehat{\cal K}^{cont}_{\underline{\widehat\theta}})$. {Hence one may always assume that $\theta \ge \underline\theta$ and $\widehat \theta\ge \underline{\widehat \theta}$ even if in practice these two inequalities may hold as equalities (see Practitioner's corner below when $K_i(t,s)= (t-s)^{\alpha_i}$, $\alpha_i>-1/i$.}
\end{Remark}

\noindent {\bf Practitioner's corner}. $(a)$ If $K_i(t,s) = (t-s)^{\alpha_i}$ with $\alpha_1> -1$ and $\alpha_2>-\frac 12$, then one checks (see e.g.~Example 2.1 in~\cite{RiTaYa2020}) that the assumptions of the above theorems are satisfied with $\beta \!\in \big( 1, \frac 1{(2\alpha_1+1)^-}\wedge\frac{1}{2\alpha_2^-} \big)$, $\theta =\widehat \theta=\underline{\theta}=\underline{\widehat\theta}= \underline{\check\theta}=\min(\alpha_1+1,\alpha_2+\frac 12,1)$.
 
Hence if $\alpha_1=\alpha_2=\alpha \!\in (-\frac 12,0)$, one has $\frac{\beta-1}{2\beta}= \frac 12 -\frac{1}{2\beta}< \alpha+\frac 12=  {\widehat \theta= \theta=\underline{\theta}=\underline{\widehat\theta}= \underline{\check\theta}}$. As a consequence, the final rate of convergence is $O\big(\big(\frac Tn\big)^{\gamma\wedge (\alpha +\frac 12)}  \big)$ at a fixed time $t$ and $O\big(\big(\frac Tn\big)^{(\gamma\wedge (\alpha +\frac 12))(1-\varepsilon)}  \big)$ for every $\varepsilon\in(0,1)$ uniformly over the interval $[0,T]$. 
 
{In particular, when considering a Volterra equation driven by a pseudo-fractional Brownian motion with Hurst constant $H\!\in (0,1]$ with $b$ and $\sigma$ Lipschitz in time ($\gamma=1$) which is the standard framework of rough  volatility models (see~\cite{GatheralJR2018}) when $H\!\in (0,\frac 12)$ or of  long memory volatility models when $H\!\in (\frac 12, 1)$ (see~\cite{ComteR1998}), then $\alpha = H-\frac 12\!\in (-\frac 12,\frac 12)$. Note that in the long memory setting,  rates $O\big((\frac Tn)^{\alpha+ \frac 12}  \big)$ or $O\big((\frac Tn)^{(\alpha+ \frac 12)(1-\varepsilon)}  \big)$ (for small enough $\varepsilon>0$) are  faster than $O\big((\frac Tn)^{ \frac 12}  \big)$.}

\medskip 
\noindent $(b)$ {\em Simulation of the $K$-integrated scheme~\eqref{eq:Euler2}}. The exact simulation of the Euler scheme~\eqref{eq:Euler2} boils down on the  one hand to computing the weight matrix
\[
\Big(\int_{t^n_{\ell-1}}^{t^n_{\ell}}K_1(t^n_k, s) ds\Big)_{1\le \ell\le k\le n}
\]
and, one the other hand,  to simulating in an exact way the independent random vectors
\[
\displaystyle Y^{(\ell)} = \Big(\int_{t^n_{\ell-1}}^{t^n_{\ell}}K_2(t^n_k, s) dW_s\Big)_{\ell\le k\le n}, \quad \ell=1,\ldots,n.
\]

The $(n-\ell+1)\times (n-\ell+1)$  (symmetric) covariance  matrix of $Y^{(\ell)}$ is given by 
\[
\Sigma^{(\ell)}= \left [ \int_{0}^{T/n}K_2(t^n_i, t^n_{\ell-1}+u) K_2(t^n_j, t^n_{\ell-1}+u)du\right]_{\ell\le i,j\le n}.
\] 
If furthermore, for every $\ell=1,\ldots,n$, the functions $K_2(t^n_k,t^n_{\ell-1}+\cdot)$, $k=\ell,\ldots,n$, are linearly independent (for  Lebesgue $a.e.$ equality on $[0,T/n]$), then the covariance matrices $\Sigma^{(\ell)}$ are all positive definite and admit a Cholesky decomposition
\[
\Sigma^{(\ell)}= T^{(\ell)}T^{(\ell)*}\quad\mbox{ with }\quad   T^{(\ell)}\; \mbox{ lower triangular}
\]
so that  
\[
(Y^{(\ell)})_{\ell=1,\ldots,n}\stackrel{d}{=} (T^{(\ell)}Z^{(\ell)})_{\ell=1,\ldots,n},\mbox{ where }Z^{(\ell)},\; \ell=1,\ldots,n,\mbox{ are independent and }  Z^{(\ell)} \sim{\cal N}(0,I_{n-\ell+1}).
\]
In case the functions $K_2(t^n_k,t^n_{\ell-1}+\cdot)$, $k=\ell,\ldots,n$ are not linearly independent, one considers the (unique) positive symmetric square root of the covariance matrix $\Sigma^{(\ell)}$ as an alternative.

\smallskip
\noindent $(c)$ {\em The case of the kernels $K_i(t,s)= (t-s)^{\alpha_i}$, $\alpha_1> -1$ and $\alpha_2>-\frac 12$}. For the ``drift kernel'' $K_1$ with $\alpha_1>-1$, closed forms are straightforward since, for every $\ell= 1,\ldots,n$ and every $i\!\in\{0, \ldots,n-\ell\}$
\[
\int_{t^n_{\ell-1}}^{t^n_{\ell}}K_1(t^n_{\ell+i}, s) ds = \frac{1}{\alpha_1 +1}\Big(\frac Tn\Big)^{\alpha_1+1}\Big((i+1)^{\alpha_1+1}-i^{\alpha_1+1} \Big).
\]
As for the ``pseudo-diffusion kernel'' $K_2$ with $\alpha_2>-1/2$,  elementary computations show that
\begin{equation}\label{eq:Sigmaell}
\Sigma^{(\ell)}_{\ell+i,\ell+j} = \Big(\frac Tn\Big)^{2\alpha_2+1}\int_0^1 \big((i+u)(j+u)\big)^{\alpha_2}du, \quad  0 \le i \le  j\le n-\ell.
\end{equation}
When $i=j$, 
\[
\Sigma^{(\ell)}_{\ell+i, \ell+i} = \frac{1}{2\alpha_2 +1}\Big(\frac Tn\Big)^{2\alpha_2+1}\Big((i+1)^{2\alpha_2+1}-i^{2\alpha_2+1} \Big).
\]
When $i\ne j$, $\Sigma^{(\ell)}_{\ell+i, \ell+j}$ can be computed offline. In the particular case $i=0$, one may use
for this purpose that
\[
\Sigma^{(\ell)}_{\ell,\ell+j} =  \frac{1}{\alpha_2 +1}\Big(\frac Tn\Big)^{2\alpha_2+1}\int_0^1(j+v^{1/(\alpha_2+1)})^{\alpha_2}dv.
\]

\subsection{Convex ordering: from random vectors to continuous processes}\label{subsec:1.2}
\begin{Definition}Let $U, V: (\Omega, \mathcal{F}, \mathbb{P})\rightarrow \big(\R^d, Bor(\R^d)\big)$ be two integrable random vectors.

\noindent $(a)$ {\em Convex ordering.}  We say that $U$ is dominated by $V$ for the \textit{ convex ordering} -- denoted by $U\preceq_{cvx} V$ -- if, for every convex function $f: \R^d\rightarrow\R$,
\begin{equation}\label{eq:defcvorder}
\E \,f(U)\leq \E\, f(V).
\end{equation}

\noindent $(b)$  {\em Increasing convex ordering ($d=1$).} We say that $U$ is dominated by $V$ for the \textit{increasing convex ordering} - denoted by $U\preceq_{icv} V$ -- if~\eqref{eq:defcvorder} holds for every non-decreasing convex function $f: \R\rightarrow\R$.

\noindent $(c)$  {\em Decreasing convex ordering ($d=1$).} We say that $U$ is dominated by $V$ for the \textit{decreasing convex ordering} - denoted by $U\preceq_{dcv} V$ -- if~\eqref{eq:defcvorder} holds for every non-increasing convex function $f: \R\rightarrow\R$.\end{Definition}
Note that any convex function $f:\R^d\to\R$ is bounded from below by an affine function so that for $U$ integrable, $\E\max(-f(U),0)< +\infty$ and $\E \,f(U)$ always makes sense in $\R\cup\{+\infty\}$.

According to Lemma A.1~\cite{AlfJou2020} for the regular convex ordering (resp. Lemma 1.1~\cite{LiuPag} for the increasing and decreasing convex orderings), the restriction of~\eqref{eq:defcvorder} to convex (resp. convex non-decreasing, resp. convex non-increasing) functions $f$ with at most affine growth is enough to characterize the order between $U$ and $V$. Reasoning like in the proof of Lemma~\ref{onlylipschitz} below, one easily checks that the test functions $f$ may even be supposed to be Lipschitz continuous. \begin{Lemma}\label{onlylineargrowth} For all integrable $\R^d$-valued (resp. $\R$-valued) random vectors $U$ and $V$, we have $U\preceq_{cvx} V$ (resp. $U\preceq_{icv}V$, resp. $U\preceq_{dcv} V$) if and only if, for every Lipschitz convex (resp. convex and non-decreasing, resp. convex non-increasing) function $f:\R^d\to \R$ (resp. $f:\R\to\R$, resp. $f:\R\to\R$) 
  , $\E \,f(U)\leq \E\, f(V)$.
\end{Lemma}

Let us extend the definition to continuous processes.
\begin{Definition}Let $X, Y: (\Omega, \mathcal{F}, \mathbb{P})\rightarrow {\cal C}([0,T],\R^d)$ be two integrable continuous processes such that  $\E[\|X\|_{\sup}+\|Y\|_{\sup}]< +\infty$.

  \noindent $(a)$ {\em Convex ordering.}  We say that $X$ is dominated by $Y$ for the \textit{ convex ordering} -- denoted by $X\preceq_{cvx} Y$ -- if, for every l.s.c. (for the uniform convergence topology) convex functional $F: {\cal C}([0,T],\R^d)\rightarrow\R$,
\begin{equation}\label{eq:defcvorderproc}
\E \,F(X)\leq \E\, F(Y).
\end{equation}
\noindent $(b)$  {\em Increasing convex ordering ($d=1$).} We say that $X$ is dominated by $Y$ for the \textit{increasing convex ordering} - denoted by $X\preceq_{icv} Y$ -- if~\eqref{eq:defcvorderproc} holds for every l.s.c. convex functional $F: {\cal C}([0,T],\R)\rightarrow\R$ non-decreasing for the pointwise partial order on continuous functions.

\noindent $(c)$  {\em Decreasing convex ordering ($d=1$).} We say that $X$ is dominated by $Y$ for the \textit{decreasing convex ordering} - denoted by $X\preceq_{dcv} Y$ -- if~\eqref{eq:defcvorderproc} holds for every l.s.c. convex functional $F: {\cal C}([0,T],\R)\rightarrow\R$ non-increasing for the pointwise partial order on continuous functions.\end{Definition}
Note that any  l.s.c. convex functional $F:{\cal C}([0,T],\R^d)\to\R$ is bounded from below by a continuous affine functional according to Lemma 7.5~\cite{AliBord} (which applies since ${\cal C}([0,T],\R^d)$ endowed with the supremum norm is, like any normed vector space, a locally convex Hausdorff space) so that for $X$ integrable, $\E\max(-F(X),0)< +\infty$ and $\E \,F(X)$ always makes sense in $\R\cup\{+\infty\}$.
\begin{Lemma}\label{onlylipschitz} For all integrable ${\cal C}([0,T],\R^d)$-valued (resp. ${\cal C}([0,T],\R)$-valued) processes $X$ and $Y$, we have $X\preceq_{cvx} Y$ (resp. $X\preceq_{icv}Y$, resp. $X\preceq_{dcv} Y$) if and only if~\eqref{eq:defcvorderproc} holds for every Lipschitz convex (resp. Lipschitz convex and non-decreasing for the pointwise partial order on continuous functions, resp. Lipschitz convex and non-increasing for the pointwise partial order on continuous functions) functional $F:{\cal C}([0,T],\R^d)\to \R$ (resp. $F:{\cal C}([0,T],\R)\to\R$, resp. $F:{\cal C}([0,T],\R)\to\R$).
\end{Lemma}
\noindent {\bf Proof.}
Let $F:{\cal C}([0,T],\R^d)\to \R$ be  l.s.c. and convex (resp.  l.s.c., convex and non-decreasing for the pointwise partial order on continuous functions, resp.  l.s.c., convex and non-increasing for the pointwise partial order on continuous functions). We are going to { establish  the existence of a }non-decreasing and bounded from below by a continuous affine functional sequence of Lipschitz convex (resp.  Lipschitz convex and non-decreasing for the pointwise order on continuous functions, resp. Lipschitz convex and non-increasing for the pointwise order on continuous functions) functionals  having $F$ as a  pointwise limit. The monotone convergence theorem then ensures that $\E\, F(X)\le \E \, F(Y)$. Since $F$ is arbitrary, $X\preceq_{cvx} Y$ (resp. $X\preceq_{icv}Y$, resp. $X\preceq_{dcv} Y$).

Let $F:{\cal C}([0,T],\R^d)\to \R$ be  l.s.c. and convex. 
For $n\in\N$, we define $F_n$ as the inf-convolution of $F$ with $n$-times the norm: 
$$
F_n(x)=\inf_{z\in {\cal C}([0,T],\R^d)}\left(F(z)+n\|z-x\|_{\sup}\right).
$$ 
By Lemma 7.5~\cite{AliBord}, there exists a continuous affine functional $G:{\cal C}([0,T],\R^d)\to \R$ such that $F\ge G$. Hence, denoting by $\|G\|$ the norm of the linear part of $G$, we have for each $x\in{\cal C}([0,T],\R^d)$,
\begin{equation}
   \forall z\in {\cal C}([0,T],\R^d),\;F(z)+n\|z-x\|_{\sup}\ge G(x)+G(z)-G(x)+n\|z-x\|_{\sup}\ge G(x)+(n-\|G\|)\|z-x\|_{\sup},\label{minofn}
\end{equation}
so that $F_n(x)\ge G(x)>-\infty$ as soon as $n\ge \|G\|$, which we now suppose from now on.
Let $x,y\in{\cal C}([0,T],\R^d)$ be such that $F_n(x)\ge F_n(y)$ and $(z_k)_{k\in\N}$ be a ${\cal C}([0,T],\R^d)$-valued sequence such that $F_n(y)=\lim_{k\to\infty}\left(F(z_k)+n\|z_k-y\|_{\sup}\right)$. Then, by the definition of $F_n(x)$ and the triangle inequality,
$$F_n(x)-F_n(y)\le \liminf_{k\to\infty}\left(F(z_k)+n\|z_k-x\|_{\sup}-F(z_k)-n\|z_k-x\|_{\sup}\right)\le n\|x-y\|_{\sup},$$
so that $F_n$ is Lipschitz with constant $n$. Since, by convexity of $F$ and $\|\cdot\|_{\sup}$, for $x,y,z,w\in{\cal C}([0,T],\R^d)$ and $\alpha\in[0,1]$,
\begin{align*}
   F(\alpha z+(1-\alpha)w)+n\|\alpha z+(1-\alpha)w-\alpha x-(1-\alpha)y\|_{\sup}\le &\alpha\left(F(z)+n\|z-x\|_{\sup}\right)\\&+(1-\alpha)\left(F(w)+n\|w-y\|_{\sup}\right),
\end{align*}
taking the infimum over $z$ and $w$, we obtain that the function $F_n$ is convex. It is also bounded from above by $F$ (choice $z=x$ in the definition of $F_n(x)$). Moreover, still by Lemma 7.5~\cite{AliBord}, for every  $x\in {\cal C}([0,T],\R^d)$ and $\alpha<F(x)$ there exists a continuous affine functional $G_{x,\alpha}$ on ${\cal C}([0,T],\R^d)$ such that $G_{x,\alpha}(x)=\alpha$ and $\forall z\in {\cal C}([0,T],\R^d)$, $F(z)>G_{x,\alpha}(z)$. Replacing $G$ by $G_{x,\alpha}$ in~\eqref{minofn}, we get that for $n\ge \|G_{x,\alpha}\|$,
$F_n(x)\ge \alpha$. Hence $\forall x\in {\cal C}([0,T],\R^d)$, $\lim_{n\to\infty}F_n(x)=F(x)$. Moreover this convergence is monotone since, clearly, $F_n\le F_{n+1}$.

If $d=1$ and $F$ is non-decreasing (resp. non-increasing) for the pointwise partial order on continuous functions, then let $x,y\in{\cal C}([0,T],\R)$ with $x\le y$ (resp. $y\le x$) for this partial order and let $(z_k)_{k\in\N}$ be a ${\cal C}([0,T],\R)$-valued sequence such that $F_n(y)=\lim_{k\to\infty}\left(F(z_k)+n\|z_k-y\|_{\sup}\right)$. For each $k\in\N$, we have $z_k+x-y\le z_k$ (resp. $z_k\le z_k+x-y$) for the partial order so that $F(z_k+x-y)\le F(z_k)$ and, by definition of $F_n(x)$, $$F_n(x)\le \liminf_{k\to +\infty}\left(F(z_k+x-y)+n\|z_k+x-y-x\|_{\sup}\right)\le \lim_{k\to +\infty}\left(F(z_k)+n\|z_k-y\|_{\sup}\right)=F_n(y).$$
Hence the monotonicity for the pointwise partial order is transferred from $F$ to $F_n$.\hfill$\cqfd$

\subsection{Main result}
We introduce a second Volterra process $(Y_t)_{t\in[0,T]}$, similar to the first one,  solving the following stochastic Volterra equation  
\begin{equation}\label{eq:Volterra2}
   Y_t= Y_0+\int_0^t   \widetilde K_1(t,s)  \widetilde b(s, Y_s)ds+\int_0^t \widetilde  K_2(t,s) \widetilde \sigma(s, Y_s)dW_s,\;t\in[0,T],
 \end{equation}
 where $ \widetilde  K_1,\widetilde  K_2:\{(t,s):0\le s<t\le T\}\to\R_+$, $\widetilde b:[0,T]\times\R^d\to\R^d$ and $ \widetilde \sigma:[0,T]\times\R^d\to\R^{d\times q}$  and $Y_0$ satisfy the assumptions made in Theorem~\ref{prop:exunvolt} on $(K_1,K_2,b,\sigma,X_0)$ so that~\eqref{eq:Volterra2} admits a unique pathwise continuous solution.
 
 \smallskip
 \noindent {$\rhd$ {\em Assumptions for the $K$-discrete Euler scheme}.} We make the following comparison assumption mixing the kernels $K_2$ and $\widetilde K_2$  and the diffusion coefficients $\sigma$ and $\widetilde \sigma$.
\begin{align}\label{eq:sigmacKompar}
\big(\mathcal{CK}^{disc}_2\sigma\big)\;\;\forall\, j\in\N^*,&\;\forall \,x\!\in\R^d,\;\forall \, s_0,\,s_1, \cdots, \,s_j\!\in [0,T]\mbox{ with }0\le s_0< s_1\le \cdots\le s_j\le T,\notag\\
&{\mathbf K}_2{\mathbf K}_2^*(s_j,\cdots,s_1,s_0)\otimes \sigma\sigma^*(s_0,x)\le \widetilde {\mathbf K}_2\widetilde {\mathbf K}_2^*(s_j,\cdots,s_1,s_0)\otimes \widetilde\sigma\widetilde\sigma^*(s_0,x), 
\end{align} 
where 
$$ 
{\mathbf K}_2(s_j,\cdots,s_1,s_0)=\left(\begin{array}{c} K_2(s_{1},s_0)\\ K_2(s_{2},s_0)\\\vdots\\ K_2(s_j,s_0)
                            \end{array}\right)\;\mbox{ and }\;\widetilde  {\mathbf K}_2(s_j,\cdots,s_1,s)=
                            \left(\begin{array}{c}\widetilde K_2(s_{1},s_0)\\
                            \widetilde K_2(s_{2},s_0)\\\vdots\\ \widetilde K_2(s_j,s_0)
                             \end{array}\right).
$$
\begin{Remark}  $\bullet$ If $d=1$ the condition~\eqref{eq:sigmacKompar} reads
\begin{align}
\nonumber \big(\mathcal{CK}^{disc}_2\sigma_{1d}\big)\; \forall (s,x)\in[0,T)\times\R,&\;\exists \lambda(s,x)\in [0,1],\;\forall t\in(s,T],\\
 &\;K_2(t,s)|\sigma(s,x)|=\lambda(s,x) \widetilde K_2(t,s)|\widetilde \sigma(s,x)|\label{comparksig}
\end{align}
 \noindent   since one checks that, for $y,z\in\R^j$, $yy^*\le zz^*$ if only if there exists $\lambda\in [-1,1]$ such that $y=\lambda z$.

\smallskip
\noindent $\bullet$   If $\widetilde K_2 = K_2$ then the above condition~\eqref{eq:sigmacKompar}  reads ${\mathbf K}_2{\mathbf K}_2^*(s_j,\cdots,s_1,s)\otimes \big(\widetilde\sigma\widetilde\sigma^*(s,x) - \sigma\sigma^*(s,x) \big)\in{\mathcal S}^+(j\times d)$. It follows from Lemma~\ref{lem:Hadamard} in Appendix~\ref{sec:matrices} that this condition is satisfied  as soon as 
\begin{align}
\big(\mathcal{C}\sigma\big)\qquad    \forall (s,x)\!\in[0,T]\times \R^d,\; &\sigma\sigma^*(s,x)\le\widetilde\sigma\widetilde\sigma^*(s,x).\label{comsig}
\end{align}
%
\noindent $\bullet$    If $\widetilde\sigma=\sigma$ and for each $s\in[0,T]$ there exists $x\in\R^d$ such that $\sigma\sigma^*(s,x)$ is positive definite, then~\eqref{eq:sigmacKompar} holds if and only if
\begin{equation}
\big(\mathcal{CK}_2\big)\qquad  \exists\,\lambda:[0,T]\to[0,1],\;\forall \,0\le s<t\le T,\; K_2(t,s)=\lambda(s) \widetilde K_2(t,s).\label{comK}
\end{equation}
{Note that if both kernels are normalized in the sense that $\int_0^tK_2(t,s)ds = \int_0^t\widetilde K_2(t,s)ds$ for $t\!\in [0,T]$, then $\lambda(s)= 1$ $ds$-$a.e.$ i.e. $K_2(t, \cdot)= \widetilde K_2(t,\cdot)$ $a.e.$ for every $t\!\in[0,T]$}.
\end{Remark}
\noindent {$\rhd$ {\em Assumptions   for the $K$-integrated Euler scheme}. While Condition~\eqref{eq:sigmacKompar} is tailored to deal with the approximation of the Volterra equation~\eqref{eq:Volterra} by the $K$-discrete Euler scheme, its counterpart for the $K$-integrated Euler scheme writes
\begin{align}\label{eq:sigmacKompar2}
\hskip-0,5cm\big(\mathcal{CK}^{int}_2\sigma\big)&\;\;\forall\, j\!\in\N^*,\;\forall \,x\!\in\R^d,\;\forall \, s_0,\,s_1, \cdots, \,s_j\!\in [0,T]\mbox{ with }0\le s_0<s_1\le \cdots\le s_j\le T,\notag\\
&\int_{s_0}^{s_1}{\mathbf K}_2{\mathbf K}_2^*(s_j,\cdots,s_1,s)ds\otimes \sigma\sigma^*(s_0,x)\le \int_{s_0}^{s_1} \widetilde {\mathbf K}_2\widetilde {\mathbf K}_2^*(s_j,\cdots,s_1,s)ds\otimes \widetilde\sigma\widetilde\sigma^*(s_0,x).
\end{align}If $K_2$ and $\widetilde K_2$ are continuous in their second variable, for $j\ge 2$, dividing by $\frac{1}{s_1-s_0}$ the inequality with $(s_j,\cdots,s_1,s)$ replaced by $(s_j,\cdots,s_2,s)$ (which is a consequence of the full inequality \eqref{eq:sigmacKompar2}), and letting $s_1\downarrow s_0$, we recover $\big(\mathcal{CK}^{disc}_2\sigma\big)$ with $(s_j,\cdots,s_1,s_0)$ replaced by $(s_j,\cdots,s_2,s_0)$. Without the continuity assumption, no comparison of $\big(\mathcal{CK}^{int}_2\sigma\big)$ and  $\big(\mathcal{CK}^{disc}_2\sigma\big)$ is possible.

Both conditions are implied by 
\begin{align}\label{eq:sigmacKompargen}
\big(\mathcal{CK}_2\sigma\big)\;\;\forall\, j\in\N^*,&\;\forall \,x\!\in\R^d,\;\forall \, s_0,s,\,s_1, \cdots, \,s_j\!\in [0,T]\mbox{ with }0\le s_0\le s< s_1\le \cdots\le s_j\le T,\notag\\
&{\mathbf K}_2{\mathbf K}_2^*(s_j,\cdots,s_1,s)\otimes \sigma\sigma^*(s_0,x)\le \widetilde {\mathbf K}_2\widetilde {\mathbf K}_2^*(s_j,\cdots,s_1,s)\otimes \widetilde\sigma\widetilde\sigma^*(s_0,x).
\end{align}
Indeed $\big(\mathcal{CK}^{disc}_2\sigma\big)$ is deduced for the choice $s=s_0$ and $\big(\mathcal{CK}^{int}_2\sigma\big)$ by integration with respect to $s$ between $s_0$ and $s_1$.

\smallskip
\noindent {$\rhd$  {\em Sufficient condition for $\big(\mathcal{CK}_2\sigma\big)$}:\[
 \mbox{$\big(\mathcal{CK}_2\big)$ (see~\eqref{comK}) $\quad$ and  $\quad$ $\big(\mathcal{C}\sigma\big)$ (see~\eqref{comsig})}.
\]}
\noindent {{\bf Proof}. Combining the decomposition 
\begin{align*}
   \widetilde {\mathbf K}_2\widetilde {\mathbf K}_2^*(s_j,\cdots,s_1,s)\,\otimes\, &\widetilde\sigma\widetilde\sigma^*(s_0,x)-{\mathbf K}_2{\mathbf K}_2^*(s_j,\cdots,s_1,s)\otimes \sigma\sigma^*(s_0,x)\\=&\left(\widetilde {\mathbf K}_2\widetilde {\mathbf K}_2^*(s_j,\cdots,s_1,s)-{\mathbf K}_2{\mathbf K}_2^*(s_j,\cdots,s_1,s)\right)\otimes \widetilde\sigma\widetilde\sigma^*(s_0,x)\\&+{\mathbf K}_2{\mathbf K}_2^*(s_j,\cdots,s_1,s)\otimes\left(\widetilde\sigma\widetilde\sigma^*(s_0,x)-\sigma\sigma^*(s_0,x)\right)
\end{align*}
with Lemma~\ref{lem:Hadamard}, one indeed easily checks that~\eqref{comK} and~\eqref{comsig} imply~$\big(\mathcal{CK}_2\sigma\big)$ (see~\eqref{eq:sigmacKompargen}). \hfill$\Box$}

 \medskip
{In terms of regularity of the kernels we need for  the $K$-discrete Euler scheme the global assumption
\begin{equation}
(\underline{\cal K}^{cont})\; \exists\,\underline\theta,\underline{\widehat\theta},\underline{\check\theta}\in(0,1]\mbox{ such that $(K_1,K_2)$ and $(\widetilde K_1,\widetilde K_2)$ satisfy }(\underline{\cal K}^{cont}_{\underline\theta}), (\underline{\widehat {\cal K}}^{cont}_{\underline{\widehat\theta}})\mbox{ and }(\underline{\check{\cal K}}^{cont}_{\underline{\check\theta}}).
\end{equation}}
{When dealing with  the $K$-integrated Euler scheme, we only need  
\begin{equation}
({\cal K}^{cont})\quad\quad\quad\exists\,\theta,\widehat\theta\in(0,1]\mbox{ such that $(K_1,K_2)$ and $(\widetilde K_1,\widetilde K_2)$ satisfy }({\cal K}^{cont}_{\theta})\mbox{ and }(\widehat{\cal K}^{cont}_{{\widehat\theta}}).\end{equation} 
which is, according to Remark~\ref{rem:hypkdisckint}, less stringent than the former $(\underline{\cal K}^{cont})$.}

\medskip                
\noindent $\rhd$ {\em Convexity assumption on $\sigma$}. We make the following {\em convexity} assumption on the matrix field $\sigma$:

\smallskip
 For every $(t,x,y, \alpha)\in[0,T]\times\R^d\times\R^d\times[0,1],\;$ there exists  $V=V_{t,x,y, \alpha}\!\in{\cal O}(q)$ such that
\begin{equation}    \label{eq:cxsig}
\big({\mathcal Conv}\big)\;\;\sigma\sigma^*(t,\alpha x+(1-\alpha)y)\le \left(\alpha \sigma(t,x)+(1-\alpha)\sigma(t,y)V\right)\left(\alpha \sigma(t,x)+(1-\alpha)\sigma(t,y)V\right)^*.
 \end{equation} 
 This condition has already been introduced in a less  general form in~\cite{LiuPag2020} (with $V=I_q$) and in the present form in~\cite{JP}. It is further discussed at the end of the current section.
\begin{Theorem}[Convex ordering] \label{thm:main} Let $X=(X_t)_{t\in [0,T]}$ and  $Y=(Y_t)_{t\in [0,T]}$ be two  Volterra  processes defined by~\eqref{eq:Volterra} and~\eqref{eq:Volterra2} respectively with $(b,\sigma)$ and $(\widetilde b,\widetilde \sigma)$ satisfying $({\cal LH}_{\gamma})$ for some $\gamma\in(0,1]$ and $(K_1,K_2)$ and $(\widetilde K_1,\widetilde K_2)$ both satisfying $\big({\cal K }^{int}_\beta\big)$ for some $\beta>1$.

 \smallskip
\noindent  $(a)$ {\em Regular convex ordering}. Assume common drift kernel $K_1=\widetilde K_1$ and common affine drift function
\begin{equation}\label{eq:b(t,x)}
\forall (t,x)\in [0,T]\times\R^d,\;b(t,x)=\widetilde b(t,x)=\mu(t)+\nu(t)x\mbox{ for functions }\mu:[0,T]\to\R^d,\;\nu:[0,T]\to\R^{d\times d}.
 \end{equation}
Assume furthermore that $\sigma$ or $\widetilde \sigma$ satisfies $\big({\mathcal Conv}\big)$
and that 
\[
\mbox{either $({\cal K}^{cont})$, $ \big(\mathcal{CK}^{int}_2\sigma\big)$ hold
or $(\underline{\cal K}^{cont})$ and $ \big(\mathcal{CK}^{disc}_2\sigma\big)$ hold.}
\]
Finally assume that $X_0$, $Y_0\!\in L^1(\P)$ and that $X_0\preceq_{cvx} Y_0$. Then, 
\[
X\preceq_{cvx} Y.
\]
\noindent $(b)$ {\em Increasing convex ordering when $d=q=1$}. Assume  either $x\mapsto b(t,x)$ and $x\mapsto |\sigma(t,x)|$ are convex and non-decreasing for every $t\!\in [0,T]$ or $x\mapsto \widetilde b(t,x)$ and $x\mapsto |\widetilde \sigma(t,x)|$ are convex and non-decreasing for every $t\!\in [0,T]$. Also assume either 
$({\cal K}^{cont})$, $\big(\mathcal{CK}_2\sigma\big)$ and 
\begin{align}
   \forall\, 0\le s_0<s_1\le s_2\le T,\;\forall x\in\R,\;b(s_0,x)\int_{s_0}^{s_1}K_1(s_2,s)ds\le\widetilde b(s_0,x)\int_{s_0}^{s_1}\widetilde K_1(s_2,s)ds
  \label{compkbs}
  \end{align}
or  $(\underline{\cal K}^{cont})$, $\big(\mathcal{CK}^{disc}_2\sigma_{1d}\big)$ (see~\eqref{comparksig}) and
\begin{align}
&\forall\, 0\le s< t\le T,\;\forall x\in\R,\; K_1(t,s)b(s,x)\le \widetilde K_1(t,s)\widetilde b(s,x)\label{comparkb}.
  \end{align}
Finally assume that $X_0, \; Y_0\!\in L^1(\P)$ with $X_0\preceq_{icv} Y_0$. Then, 
\[
X\preceq_{icv} Y
\]

\smallskip
\noindent $(c)$ {\em Decreasing convex ordering when $d=q=1$}. Assume  either $x\mapsto -b(t,x)$ and $x\mapsto |\sigma(t,x)|$ are convex and non-increasing for every $t\!\in [0,T]$ or $x\mapsto -\widetilde b(t,x)$ and $x\mapsto |\widetilde \sigma(t,x)|$ are convex and non-increasing for every $t\!\in [0,T]$. Also assume either
 $({\cal K}^{cont})$, $\big(\mathcal{CK}_2\sigma\big)$ and 
\begin{align}
   \forall\, 0\le s_0<s_1\le s_2\le T,\;\forall x\in\R,\;b(s_0,x)\int_{s_0}^{s_1}K_1(s_2,s)ds\ge\widetilde b(s_0,x)\int_{s_0}^{s_1}\widetilde K_1(s_2,s)ds
\label{compkbbs}
\end{align}
or  $(\underline{\cal K}^{cont})$, $\big(\mathcal{CK}^{disc}_2\sigma_{1d}\big)$ and
\begin{align}
  &\forall\, 0\le s<t\le T,\;\forall x\in\R,\; K_1(t,s)b(s,x)\ge \widetilde K_1(t,s)\widetilde b(s,x)\label{comparkbb}.
  \end{align}
Finally assume that $X_0, \; Y_0\!\in L^1(\P)$ with $X_0\preceq_{dcv} Y_0$. Then
\[
X\preceq_{dcv} Y.
\]
\end{Theorem}
\begin{Remark}
  Under~\eqref{eq:b(t,x)}, for $b=\widetilde b$ to satisfy $({\cal LH}_{\gamma})$ (see~\eqref{eq:HolLipbsig}), it is enough that both functions $\mu$ and $\nu$ are H\"older continuous with exponent $\gamma$.
\end{Remark}
\begin{Theorem}[Convexity]\label{thm:ci}  Let $(b,\sigma)$ satisfy $({\cal LH}_{\gamma})$ for some $\gamma\in(0,1]$ and $(K_1,K_2)$ satisfy $\big({\cal K }^{int}_\beta\big)$ for some $\beta>1$ and $({\cal K}^{cont}_{\theta})$ and $(\widehat{\cal K}^{cont}_{{\widehat\theta}})$ for some $\theta,\widehat\theta>0$. Let $X^x=(X^x_t)_{t\in [0,T]}$ denote the solution starting from $x\in\R^d$ to the Volterra stochastic differential equation~\eqref{eq:Volterra}.

 \smallskip
\noindent  $(a)$ {\em Convexity in the initial condition}. Assume \begin{equation*}\forall (t,x)\in [0,T]\times\R^d,\;b(t,x)=\mu(t)+\nu(t)x\mbox{ for functions }\mu:[0,T]\to\R^d,\;\nu:[0,T]\to\R^{d\times d},
 \end{equation*}and $\sigma$ satisfies the convexity assumption~$\big({\mathcal Conv}\big)$
 . Then, for every convex functional $F: {\cal C}([0,T],\R^d)\to \R$ with at most polynomial growth with respect to the $\sup$-norm, 
$\R^d\ni x\mapsto \E\, F(X^x)$ is convex. 

\smallskip
\noindent $(b)$ {\em Non-decreasing convexity in the initial condition when $d=q=1$}. Assume that $x\mapsto b(t,x)$ and $x\mapsto |\sigma(t,x)|$ are convex and non-decreasing for every $t\!\in [0,T]$. Then, for every convex functional $F: {\cal C}([0,T], \R)\to \R$ with at most polynomial growth with respect to the $\sup$-norm and non-decreasing for the pointwise partial order on continuous functions, $\R\ni x\mapsto \E\, F(X^x)$  is convex and non-decreasing.

\smallskip
\noindent $(c)$ {\em Non-Increasing convexity in the initial condition when $d=q=1$}. Assume that $x\mapsto -b(t,x)$ and $x\mapsto |\sigma(t,x)|$ are convex and non-increasing for every $t\!\in [0,T]$. Then, for every convex functional $F: {\cal C}([0,T], \R)\to \R$ with at most polynomial growth with respect to the $\sup$-norm and non-increasing for the pointwise partial order on continuous  functions, $\R\ni x\mapsto \E\, F(X^x)$  is convex and non-increasing.
\end{Theorem}
Let us state  some equivalent formulations of the convexity condition  $\big({\mathcal Conv}\big)$ 
that may appear more ``symmetric'' or more ``intrinsic'' and explain why, in the previous theorem, we suppose  when $d=q=1$ that $x\mapsto|\sigma(t,x)|$ is convex for every $t\in[0,T]$. 
\begin{Proposition} $(a)$ The convexity condition~$\big({\mathcal Conv}\big)$ (see~\eqref{eq:cxsig})
  is equivalent to 
 
 \smallskip
 \noindent $(i)$ For every $(t,x,y, \alpha)\in[0,T]\times\R^d\times\R^d\times[0,1],\;$ there exists  $U=U_{t,x,y, \alpha}$, $V=V_{t,x,y, \alpha}\!\in{\cal O}(q)$ such that
\begin{equation*}
   \label{eq:cxsig2}\sigma\sigma^*(t,\alpha x+(1-\alpha)y)\le \left(\alpha \sigma(t,x)U+(1-\alpha)\sigma(t,y)V\right)\left(\alpha \sigma(t,x)U+(1-\alpha)\sigma(t,y)V\right)^*.
 \end{equation*} 
 
 \smallskip
 \noindent $(ii)$ For every $ (t,x,y, \alpha)\!\in[0,T]\times\R^d\times\R^d\times[0,1]$
 \begin{align}\exists\,\varsigma,\widetilde\varsigma\in\R^{d\times q},\begin{cases}\varsigma\varsigma^*=\sigma\sigma^*(t,x),\;\widetilde\varsigma\widetilde\varsigma^*=\sigma\sigma^*(t,y),\mbox{ and }\\
    \sigma\sigma^*(t,\alpha x+(1-\alpha)y)\le \left(\alpha \varsigma+(1-\alpha)\widetilde \varsigma \right)\left(\alpha \varsigma+(1-\alpha)\widetilde \varsigma \right)^{\!*}
     \end{cases}.\label{cxsigbis}
\end{align}

\smallskip
\noindent $(b)$ If $q\le d$,  the convexity condition~\eqref{eq:cxsig} implies that:

  \smallskip
 for every $(t,x,y, \alpha)\in[0,T]\times\R^d\times\R^d\times[0,1],\;$ there exists  $V=V_{t,x,y, \alpha}\!\in{\cal O}(d)$ such that
\begin{equation*}
  \sigma\sigma^*(t,\alpha x+(1-\alpha)y)\le \left(\alpha \sqrt{\sigma\sigma^*}(t,x)+(1-\alpha)\sqrt{\sigma\sigma^*}(t,y)V\right)\left(\alpha \sqrt{\sigma\sigma^*}(t,x)+(1-\alpha)\sqrt{\sigma\sigma^*}(t,y)V\right)^{\!*},
 \end{equation*} 
 with equivalence when $q=d$.

\smallskip
\noindent $(c)$ When $d=q=1$, the convexity condition~\eqref{eq:cxsig} reads $x\mapsto |\sigma(t,x)|$ is convex for each $t\in[0,T]$.
\end{Proposition}

\noindent {\bf Proof.} $(a)$ The equivalence of $(i)$ and $(ii)$ is an easy consequence of the fact that $A,B\in \R^{d\times q}$ are such that $AA^*=BB^*$ if and only if $A=BO$ for some  $O\!\in{\cal O}(q)$ (see Lemma~\ref{lem:matrices} in Appendix~\ref{sec:matrices}). 
Condition $(i)$ implies~\eqref{eq:cxsig} (with $U=I_q$). These are in fact equivalent since, when $U\in{\mathcal O}(q)$, by introducing $U^*U=I_q$ between the two terms of the product in the left-hand side,
\begin{align*}
  &\left(\alpha \sigma(t,x)U+(1-\alpha)\sigma(t,y)V\right)\left(\alpha \sigma(t,x)U+(1-\alpha)\sigma(t,y)V\right)^*\\&\phantom{\alpha \sigma(t,x)U+(1-\alpha)\sigma(t,y)V}=\left(\alpha \sigma(t,x)+(1-\alpha)\sigma(t,y)VU^*\right)\left(\alpha \sigma(t,x)+(1-\alpha)\sigma(t,y)VU^*\right)^*
\end{align*}
with $VU^*\in{\mathcal O}(q)$ if $V$ also belongs to ${\mathcal O}(q)$.

\smallskip\noindent $(b)$ {When $q\le d$, since the matrix $\pi\in\R^{q\times d}$ with only non vanishing entries $\pi_{ii}=1$ for $i\in\{1,\cdots,q\}$, satisfies $\pi\pi^*=I_q$, we have that for $V\in{\cal O}(q)$, 
\begin{align*}
   (\alpha \sigma(t,x)+(1-\alpha)&\sigma(t,y)V)\left(\alpha \sigma(t,x)+(1-\alpha)\sigma(t,y)V\right)^*\\&=\left(\alpha \sigma(t,x)\pi+(1-\alpha)\sigma(t,y)V\pi\right)\left(\alpha \sigma(t,x)\pi+(1-\alpha)\sigma(t,y)V\pi\right)^*.
\end{align*}
Since the matrices $\varsigma=\sigma(t,x)\pi\in\R^{d\times d}$ and $\widetilde\varsigma=\sigma(t,y)V\pi\in\R^{d\times d}$ are such that $\varsigma\varsigma^*=\sigma\sigma^*(t,x)$ and $\widetilde\varsigma\widetilde\varsigma^*=\sigma\sigma^*(t,y)$, we conclude by $(a)$$(ii)$ applied to $\sqrt{\sigma\sigma^*}$ in place of $\sigma$. The reverse implication when $q=d$ also follows from $(a)$$(ii)$.}

\smallskip\noindent  $(c)$ When $d=q=1$, then the right-hand side of~\eqref{eq:cxsig} is maximal for $U$ equal to $1$ or $-1$ depending on whether $\sigma(t,x)$ and $\sigma(t,y)$ share the same sign or not. Claim $(c)$ easily follows.\hfill$\Box$

\begin{Remark} One derives from claim $(b)$ the following  {sufficient} natural criterion for $\big({\mathcal Conv}\big)$ to hold when $q{=}d$:
for every $(t,x,y, \alpha)\in[0,T]\times\R^d\times\R^d\times[0,1],\;$
\begin{equation*}
   \label{eq:cxsig2}\sigma\sigma^*(t,\alpha x+(1-\alpha)y)\le \left(\alpha \sqrt{\sigma\sigma^*}(t,x)+(1-\alpha)\sqrt{\sigma\sigma^*}(t,y)\right)\!\!\left(\alpha \sqrt{\sigma\sigma^*}(t,x)+(1-\alpha)\sqrt{\sigma\sigma^*}(t,y)\right)^{\!*}.
 \end{equation*} 
  {Note that, when $q\le d$,  we can always assume that $q=d$ by replacing {\em mutatis mutandis} $\sigma$  by $\sigma_{\hookrightarrow}(t,x):= \sigma(t,x)\pi$ where $\pi\in\R^{q\times d}$ is the matrix with only non vanishing entries $\pi_{ii}=1$ for $i\in\{1,\cdots,q\}$ and adding to the $q$-dimensional standard  Brownian motion $(W_t)_{t\in[0,T]}$ $d-q$ coordinates corresponding to an independent $d-q$-dimensional standard Brownian motion.}\end{Remark}
 
The strategy of proof of both Theorems~\ref{thm:main}~and~\ref{thm:ci} will rely on two steps like in~\cite{Pag2016}:
 
\smallskip
 -- establish the result for the Euler schemes by combining propagation of convexity (resp. increasing convexity) in Section~\ref{sbpc} (resp. Section~\ref{sbpci}) and comparison in Section~\ref{sec:compar} (resp. Section~\ref{sec:compari})  for the convex ordering (resp. increasing convex ordering).
 
 \smallskip
 -- transfer the convex ordering   by making the schemes converge (see Section~\ref{sec:Cvgce}). 
 \subsection{Application to VIX options in  the quadratic rough Heston model}

 Let us consider the auxiliary variance process in the  quadratic rough Heston model (see~\cite{GaJuRo2020}):
 \[
 V_t =  a(Z_t-b)^2+c \quad \mbox{ with } \quad a>0,\,b,\,c\ge 0
 \]
 and, for $H \!\in (0,1/2)$, $\lambda\in\R$, $\sigma>0$ and $f:[0,T]\to\R$ H\"older continuous with exponent $\gamma\in (0,1]$,
 \begin{equation}\label{eq:Z}
 Z_t = Z_0 +\int_0^t (t-s)^{H-\frac 12} \lambda (f(s) -Z_s) ds +\sigma\int_0^t  (t-s)^{H-\frac 12}\sqrt{a(Z_s-b)^2+c}\,dW_s.
\end{equation}
Note that $z\mapsto \sqrt{a(z-b)^2+c}$ is convex and Lipschitz. Then the above Volterra equation has a unique strong solution, {temporarily} denoted by  $(Z^{(\sigma)}_ t)_{t\ge 0}$ to emphasize the dependence on the parameter $\sigma$. Let $V^{(\sigma)}_t=a(Z^{(\sigma)}_t-b)^2+c$ denote the resulting squared volatility.
 
According to Theorem \ref{thm:main}, one has $(Z^{(\sigma)}_ t)_{t\in [0,T]}\le_{cvx} (Z^{(\widetilde \sigma)}_ t)_{t\in [0,T]}$ for $\sigma \! \in (0,\widetilde\sigma]$.
 
 The price at time $0$ of a VIX contract with maturity  $T>0$ in such a  model is given by 
 \[
 \E \left ( \sqrt{\frac 1T \int_0^T V^{(\sigma)}_ t dt }\right)= \E \left ( F((Z^{(\sigma)}_t)_{t\in[0,T]})\right)\mbox{ where }F((x_t)_{t\in[0,T]})=\sqrt{\frac 1T \int_0^T\left(a(x_t-b)^2+c \right) dt}.
 \]
 
For $(x_t)_{t\in[0,T]},(y_t)_{t\in[0,T]}\in{\mathcal C}([0,T],\R)$ and $\alpha\in [0,1]$, by the convexity of $z\mapsto \sqrt{a(z-b)^2+c}$ and the monotonicity of the square function on $\R_+$ then the convexity of the $L^2(dt)$ norm, 
\begin{align*}
   F((\alpha x_t+(1-\alpha)y_t)_{t\in[0,T]})&\le \sqrt{\frac 1T \int_0^T\left(\alpha \sqrt{a(x_t-b)^2+c}+(1-\alpha)\sqrt{a(y_t-b)^2+c} \right)^2 dt}\\&\le \alpha F\big((x_t)_{t\in[0,T]}\big)+(1-\alpha)F\big((y_t)_{t\in[0,T]}\big).
\end{align*}
Moreover, the functional $F$ has at most at most affine growth with respect to the sup-norm. Hence, by Theorem~\ref{thm:main}, $\E \left ( F((Z^{(\sigma)}_t)_{t\in[0,T]})\right)\le\E \left ( F((Z^{(\widetilde \sigma)}_t)_{t\in[0,T]})\right)$ which also writes  
\[
\E \left ( \sqrt{\frac 1T \int_0^T V^{(\sigma)}_ t dt }\right) \leq  \E \left (\sqrt{ \frac 1T \int_0^T V^{(\widetilde \sigma)}_ t dt} \right) 
\]
i.e. the premium is a non-decreasing function of the parameter $\sigma$.

Moreover, owing to Theorem~\ref{thm:ci}, if we now denote by $Z^{z_0}$ the solution to \eqref{eq:Z} starting from $Z_0=z_0\!\in \R$ and $V^{z_0}_t=a\big(Z^{z_0}_t-b\big)^2+c $,  then $z_0\mapsto \E\, F(Z^{z_0})$ is convex for every convex functional $F: {\cal C}([0,T],\R)\to \R$ with at most polynomial growth with respect to the $\sup$-norm. In particular, $z_0\mapsto \E\, f(V^{z_0}_t)$ is convex for every non-decreasing convex function $f:\R\to\R$ with polynomial growth and every $t\!\in[0,T]$.
 \section{Convex ordering and convexity propagation for the Euler schemes}\label{sec:propag}
Let us recall in a synthetic way the two   Euler discretization schemes associated to the regular mesh $(t^n_k=\frac{kT}{n})_{0\le k\le n}$ with $n\ge 1$ steps of the Volterra stochastic differential equations~\eqref{eq:Volterra} and~\eqref{eq:Volterra2}:
 \begin{align*}
   \bar X_{t^n_k}&=X_0+ \sum_{\ell=1}^{k} b(t^n_{\ell-1},\bar X_{t^n_{\ell-1}})\int_{t^n_{\ell-1}}^{t^n_{\ell}}K_1(t^n_k,{\hat s})ds+\sum_{\ell=1}^{k} \sigma(t^n_{\ell-1},\bar X_{t^n_{\ell-1}})\int_{t^n_{\ell-1}}^{t^n_{\ell}} K_2(t^n_k,{\hat s})dW_s, \; k=0,\ldots,n,\\
   \bar Y_{t^n_k}&=Y_0+ \sum_{\ell=1}^{k}\widetilde  b(t^n_{\ell-1},\bar Y_{t^n_{\ell-1}})\int_{t^n_{\ell-1}}^{t^n_{\ell}}\widetilde K_1(t^n_k,\hat s)ds+\sum_{\ell=1}^{k}\widetilde  \sigma(t^n_{\ell-1},\bar Y_{t^n_{\ell-1}})\int_{t^n_{\ell-1}}^{t^n_{\ell}}\widetilde K_2(t^n_k,\hat s)dW_s, \; k=0,\ldots,n,
 \end{align*}
 with $\hat s=s$ for the $K$-integrated Euler scheme and $\hat s=\underline s$ for the $K$-discrete Euler scheme.
 
To recover Markovianity, we also introduce for $k\in\{1,\ldots,n\}$, the companion processes $(X^k_{t^n_\ell})_{0\le\ell\le k}$ and $(Y^k_{t^n_\ell})_{0\le\ell\le k}$ respectively starting from $X^k_0=X_0$ and $Y^k_0=Y_0$ and evolving inductively as follows
\begin{align}
   \hskip-0.5 cm X^k_{t^n_{\ell+1}}\!&=X^k_{t^n_\ell}+ b(t^n_\ell,X^\ell_{t^n_\ell})\! \int_{t^n_\ell}^{t^n_{\ell+1}} K_1(t^n_k,\hat s)ds+\sigma(t^n_\ell,X^\ell_{t^n_\ell})\!\int_{t^n_\ell}^{t^n_{\ell+1}}\!\!K_2(t^n_k,\hat s)dW_s, \, \ell= 0,\ldots, k-1, \label{eq:companion}\\
Y^k_{t^n_{\ell+1}}\!&=Y^k_{t^n_\ell}+ \widetilde  b(t^n_\ell,Y^\ell_{t^n_\ell}) \int_{t^n_\ell}^{t^n_{\ell+1}} \widetilde K_1(t^n_k,\hat s)ds+\widetilde \sigma(t^n_\ell,Y^\ell_{t^n_\ell})\int_{t^n_\ell}^{t^n_{\ell+1}} \widetilde K_2(t^n_k,\hat s)dW_s, \; \ell= 0,\ldots, k-1,\label{eq:companiony}\end{align}
 so that $(X^k_{t^n_k},\cdots,X^n_{t^n_k})$ (resp. $(Y^k_{t^n_k},\cdots,Y^n_{t^n_k})$) is a function of $(X^{k-1}_{t^n_{k-1}},\cdots,X^n_{t^n_{k-1}},(W_s-W_{t^n_{k-1}})_{s\in[t^n_{k-1},t^n_k]})$ (resp. $(Y^{k-1}_{t^n_{k-1}},\cdots,Y^n_{t^n_{k-1}},(W_s-W_{t^n_{k-1}})_{s\in[t^n_{k-1},t^n_k]})$).
 We then have 
 $$
 X^k_{t^n_k}=\bar X_{t^n_k}\quad\mbox{and}\quad Y^k_{t^n_k}=\bar Y_{t^n_k} \quad\mbox{for}\quad k\in\{1,\ldots,n\}.
 $$
In fact, the convex ordering of the Euler schemes will be obtained as a consequence of the convex ordering of their companion processes for which we introduce the following notation: $X^\bullet_{0:0}=X_0$, $Y^\bullet_{0:0}=Y_0$  and, for $k\in\{1,\ldots,n\}$,
\begin{align*}
  X^{\bullet}_{0:k}&=\big(X_0,X^1_{t^n_1},\cdots,X^n_{t^n_1},X^2_{t^n_2},\cdots,X^n_{t^n_2},\cdots,X^{k}_{t^n_{k}},\hdots,X^n_{t^n_{k}}\big),\\
  Y^{\bullet}_{0:k}&=\big(Y_0,Y^1_{t^n_1},\cdots,Y^n_{t^n_1},Y^2_{t^n_2},\cdots,Y^n_{t^n_2},\cdots,Y^{k}_{t^n_{k}},\hdots,Y^n_{t^n_{k}}\big).
\end{align*}

Let us suppose 
\begin{align}
                 &\forall k\in\{0,\cdots,n-1\},\;\sup_{x\in\R^d}\frac{|b(t^n_k,x)|+|\widetilde b(t^n_k,x)|+\|\sigma(t^n_{k},x)\|+\|\widetilde\sigma(t^n_{k},x)\|}{1+|x|}< +\infty,\label{sigcaff}\\                  &\forall k\in\{1,\cdots,n\},\;\int_0^{t^n_k} \left(K_1(t^n_k,\hat s)+K_2(t^n_k,\hat s)^2+\widetilde K_1(t^n_k,\hat s)+\widetilde K_2(t^n_k,\hat s)^2\right)ds< +\infty.\label{supteul}
\end{align}
Note that for the $K$-discrete Euler scheme, $\hat s=t^n_\ell$ for $s\in[t^n_\ell,t^n_{\ell+1})$ so that~\eqref{supteul} is always satisfied.             
Since the Gaussian random vector $ \left(\int_{t^n_\ell}^{t^n_{\ell+1}}K_2(t^n_{\ell+1},\hat s)dW_s,\cdots,\int_{t^n_\ell}^{t^n_{\ell+1}}K_2(t^n_{n},\hat s)dW_s\right)$ is independent from $X^\ell_{t^n_\ell}$, one easily checks by induction on $\ell$ that when $\|X_0\|_p< +\infty$ for some $p\ge 1$, then $\|X^{\bullet}_{0:n}\|_p< +\infty$. In the same way,  when $\|Y_0\|_p< +\infty$, then $\|Y^{\bullet}_{0:n}\|_p< +\infty$.
\begin{Proposition}\label{prop:convexeul} Assume $X_0\le_{cvx} Y_0$ with $X_0,Y_0\in L^1(\P)$,~\eqref{sigcaff}, that $\sigma$ or $\widetilde \sigma$ satisfies the convexity assumption~$\big({\mathcal Conv}\big)$ (see~\eqref{eq:cxsig}),
\begin{align}
   \forall k\in\{0,\cdots,n-1\},\;b(t^n_k,x)=\widetilde b(t^n_k,x)=\mu(t^n_k)+\nu(t^n_k)x\mbox{ with }(\mu(t^n_k),\nu(t^n_k))\in\R^d\times\R^{d\times d}\label{baff}
\end{align}
and either both  Euler schemes are $K$-discrete and $(K_2,\sigma)$ and $(\widetilde K_2, \widetilde \sigma)$ satisfy  the comparison inequality~$\big(\mathcal{CK}^{disc}_2\sigma\big)$ (see~\eqref{eq:sigmacKompar}) or both  Euler schemes are $K$-integrated, ~\eqref{supteul} holds, and $(K_2,\sigma)$ and $(\widetilde K_2, \widetilde \sigma)$ satisfy  the comparison inequality~$\big(\mathcal{CK}^{int}_2\sigma\big)$ (see~\eqref{eq:sigmacKompar2}). Then, for every $n\ge 1$,  
$$
X^{\bullet}_{0:n}\le_{cvx}Y^{\bullet}_{0:n}\mbox{ and, in particular, }(\bar X_{t^n_k})_{k=0,\ldots,n}\le_{cvx}(\bar{Y}_{t^n_k})_{k=0,\ldots,n}.
$$
\end{Proposition}
\begin{Remark}
  $\bullet$ Assumption~\eqref{baff} clearly implies that the part concerning $(b,\widetilde b)$ in assumption~\eqref{sigcaff} is satisfied.

\smallskip
\noindent $\bullet$ The conclusion still holds for the $K$-discrete (resp. $K$-integrated) Euler scheme when $\sigma$ or $\widetilde \sigma$ satisfies the convexity assumption~$\big({\mathcal Conv}\big)$ for $t$ restricted to $\{t^n_k:0\le k\le n-1\}$ and $(K_2,\sigma)$ and $(\widetilde K_2, \widetilde \sigma)$ satisfy  the comparison inequality~$\big(\mathcal{CK}^{int}_2\sigma\big)$ (resp.~$\big(\mathcal{CK}^{int}_2\sigma\big)$) with $(s_j,\cdots,s_1,s_0)$ of the form $(t^n_n,t^n_{n-1},\cdots,t^n_{n+1-j},t^n_{n-j})$ for $j\in\{0,\hdots,n\}$ . Note that the comparison inequalities are only used to derive~\eqref{usesigmacompar} below.
 \end{Remark}
 
\subsection{Backward propagation of convexity by the Euler schemes 
}\label{sbpc}
Let for $k\in\{1,\ldots,n\}$, $x^\bullet_{0:k}=(x_0,x^1_{1},\cdots,x^n_{1},x^2_{2},\cdots,x^n_{2},\cdots,x^{k}_{k},\hdots,x^n_{k})$ denote a generic element of ${\cal E}_k:=\R^d\times (\R^{ d})^n\times (\R^{ d})^{n-1}\times \cdots\times (\R^d)^{n+1-k}$ and let $x^\bullet_{0:0}=x_0$ be a generic element of ${\cal E}_0=\R^d$.   Set, for every $k\!\in \{0, \ldots,n-1\}$, $\ell\in\{k+1,\cdots,n\}$, $(x,y)\in\R^d\times\R^d$, 
\begin{equation}\label{eq:B(t,x)}
B_k(t^n_\ell,x,y) = x + \big(\mu(t^n_k)+\nu(t^n_k)y\big)\int_{t^n_k}^{t^n_{k+1}} K_1(t^n_\ell,\hat s)ds.
\end{equation}

 Let $\Phi_n:{\cal E}_n\to \R$ be a function with $p$-polynomial growth for some $p\ge 1$. We define  the functions $\Phi_k: {\cal E}_k\to \R$, $k\in\{0, \cdots,n-1\}$ and the auxiliary functions $\Psi_k$  by a backward induction as follows:
for every $(x^\bullet_{0:k},u)\!\in {\cal E}_k\times \R^{d\times q}$,
\begin{align}\label{eq:Psik}
   \Psi_k(x^\bullet_{0:k}, u) = \E\,\Phi_{k+1} \bigg(x^\bullet_{0:k}, B_k(t^n_{k+1}, x_k^{k+1},x_k^k)+u\int_{t^n_k}^{t^n_{k+1}} \hskip -0,5cm &K_2(t^n_{k+1},\hat s)dW_s,\cdots,\notag\\&B_k(t^n_n,x^n_k,x_k^k)+u\int_{t^n_k}^{t^n_{k+1}} \hskip -0,5cm K_2(t^n_{n},\hat s)dW_s\bigg)
\end{align}
\begin{equation}\label{eq:Phik}
\mbox{and } \hskip 3.85cm \forall x^\bullet_{0:k}\!\in {\cal E}_k,\;\Phi_k(x^\bullet_{0:k}) = \Psi_{k}\big(x^\bullet_{0:k},\sigma(t^n_k,x_k^k)\big),\hskip 3.5cm
\end{equation}
where for $k=0$, we use the convention $x^0_0=x^1_0=\cdots=x^n_0=x_0$ (consistent with the definition of the $X^\ell$ when $\ell\ge 1$).
Note that under the part in~\eqref{sigcaff} concerning $\sigma$ and the part in~\eqref{supteul} concerning $(K_1,K_2)$ (which appears in the assumptions of the proposition just below), the functions $x\mapsto \sigma(t_k^n,x)$ have at most affine growth and 
$$
\left(\int_{t^n_k}^{t^n_{k+1}} K_2(t^n_{k+1},\hat s)dW_s,\cdots,\int_{t^n_k}^{t^n_{k+1}}K_2(t^n_{n},\hat s)dW_s\right)
$$
is a Gaussian random vector and thus belongs to $L^p(\P)$. With the affine property of the functions $(x,y)\mapsto B_k(t^n_\ell,x,y)$ for $\ell\ge k+1$, one easily deduce by backward induction that the random variable in the expectation in the right-hand side of~\eqref{eq:Psik}
is integrable and the subsequently well-defined functions $\Psi_k$ and $\Phi_k$ have $p$-polynomial growth for $k\in\{0,\cdots,n-1\}$.
\begin{Proposition}[Backward propagation of convexity]\label{prop:convexprop} Assume~\eqref{baff},\begin{align*}
  \forall k\in\{0,\cdots,n-1\},\;&\int_0^{t^n_{k+1}} \left(K_1(t^n_{k+1},\hat s)+K_2(t^n_{k+1},\hat s)^2\right)ds+\sup_{x\in\R^d}\frac{\|\sigma(t^n_{k},x)\|}{1+|x|}< +\infty,\qquad\qquad\\
\hskip 2cm\quad\forall x\in\R^d,&\quad b(t^n_k,x)=\mu(t^n_k)+\nu(t^n_k)x,
 \end{align*}
  and that the function $\sigma$ satisfies~$\big({\mathcal Conv}\big)$ (see~\eqref{eq:cxsig}) for $t\in\{t^n_k:0\le k\le n-1\}$. Let $\Phi_n:{\cal E}_n\to \R$ be convex with $p$-polynomial growth for some $p\ge 1$. Then, the  functions ${\cal E}_k\ni x^\bullet_{0:k}\mapsto \Phi_k(x^\bullet_{0:k})$, $k=0,\ldots,n$,  defined by~\eqref{eq:Phik},~\eqref{eq:Psik} 
  are convex with $p$-polynomial growth. In particular, if $\Phi_n(x^{\bullet}_{0:n}) = F(x_0,x_1^1,\cdots,x_n^n)$ with $F: (\R^d)^{n+1}\to \R$ convex with polynomial growth, then, convexity is transferred to the function $\R^d\ni x_{0}\mapsto \Phi_0(x_{0})= \E\,F\big((\bar X^{x_0}_{t^n_k})_{k=0,\ldots,n}\big)$ where $(\bar X^{x_0}_{t^n_k})_{k=0,\ldots,n}$ denotes the (either $K$-discrete or $K$-integrated) Euler scheme starting from $X_0=x_0$.

\end{Proposition}

\noindent {\bf Proof.} The $p$-polynomial growth of the functions $\Phi_k$ has been established just before the proposition. Let us check by backward induction that they are convex.

For $u\in\R^{d\times q}$,  the $\R^{(n-k)d}$-valued random vector $\left(u\int_{t^n_k}^{t^n_{k+1}} \hskip -0,15cm K_2(t^n_{k+1},\hat s)dW_s,\cdots,u\int_{t^n_k}^{t^n_{k+1}} \hskip -0,15cm K_2(t^n_{n},\hat s)dW_s\right)$ is distributed according to the centered Gaussian distribution ${\cal N}_{(n-k)d}(0,\Gamma^k\otimes uu^*)$ where $$\left(\Gamma^k_{ij}=\int_{t^n_k}^{t^n_{k+1}}
  K_2(t^n_{k+i},\hat s) 
  K_2(t^n_{k+j},\hat s)ds\right)_{1\le i,j\le n-k}$$ is a symmetric semidefinite positive  matrix. As a consequence,
\begin{equation}
   \forall(x^\bullet_{0:k},u,v)\!\in {\cal E}_k\times \R^{d\times q}\times \R^{d\times q}\mbox{ with }uu^*=vv^*,\;\Psi_k(x^\bullet_{0:k},u)=\Psi_k(x^\bullet_{0:k},v)\label{invrotbis}.
\end{equation}

For $v\in\R^{d\times q}$ such that $u^*u - vv^*\!\in {\cal S}^+(d)$ then  $\Gamma^k\otimes vv^*-\Gamma^k\otimes uu^*=\Gamma^k\otimes(vv^*-uu^*)   \!\in {\cal S}^+((n-k)d)$   owing to Lemma~\ref{lem:Hadamard}.
                           
Since we know that two centered Gaussian distributions satisfy ${\cal N}(0, \Sigma_1)\le_{cvx}{\cal N}(0, \Sigma_2)$ if and only if $\Sigma_2\Sigma_2^*-\Sigma_1\Sigma_1^* \!\in {\cal S}^+(d)$
(see for instance~\cite[Lemma 3.2 and Remark 3.1]{JP}), we deduce that 
\begin{align*}
                           \Big(u\int_{t^n_k}^{t^n_{k+1}} K_2(t^n_{k+1},\hat s)dW_s,\cdots,u\int_{t^n_k}^{t^n_{k+1}}  & K_2(t^n_{n},\hat s)dW_s \Big)\\
                           &\le_{cvx}\Big(v\int_{t^n_k}^{t^n_{k+1}} K_2(t^n_{k+1},\hat s)dW_s,\cdots,v\int_{t^n_k}^{t^n_{k+1}} K_2(t^n_{n},\hat s)dW_s\Big).
\end{align*}                          
                           If, for $k\in\{0,\cdots,n-1\}$, $\Phi_{k+1}$ is convex, then so is  $${\cal E}_k\times(\R^d)^{n-k}\ni (x^\bullet_{0:k},w_{k+1},\hdots,w_n)\mapsto \Phi_{k+1}(x^\bullet_{0:k},B_k(t^n_{k+1},x_k^{k+1},x_k^k)+w_{k+1},\cdots,B_k(t^n_n,x_k^{n},x_k^k)+w_n)$$ by the affine property of $\R^d\times\R^d\ni (x,y)\mapsto B_k(t^n_\ell,x,y)$ for $\ell\in\{k+1,\cdots,n\}$ (see~\eqref{eq:B(t,x)}). We deduce that 
\begin{equation}
   \forall x^\bullet_{0:k}\in {\mathcal E}_k,\;\forall u,v\in \R^{d\times q}\mbox{ s.t. }uu^*\le vv^*,\;\Psi_k(x^\bullet_{0:k}, u)\le \Psi_k(x^\bullet_{0:k},v)\label{eqcroisG}
 \end{equation}
 and $\Psi_k$ is convex.
For $x^\bullet_{0:k},y^\bullet_{0:k}\in{\mathcal E}_k$ and $\alpha\in[0,1]$, with the existence of $O\in{\mathcal O}(q)$ (possibly depending on $(t^n_k,x^k_k,y^k_k,\alpha)$) such that
$$
\sigma\sigma^*(t^n_k,\alpha x^k_k+(1-\alpha)y^k_k)\le \left(\alpha \sigma(t^n_k,x^k_k)+(1-\alpha)\sigma(t^n_k,y^k_k)O\right)\left(\alpha \sigma(t^n_k,x^k_k)+(1-\alpha)\sigma(t^n_k,y^k_k)O\right)^*
$$
(see~\eqref{eq:cxsig}), the convexity of $\Psi_k$ and~\eqref{invrotbis}, we conclude that
\begin{align*}
\Phi_k\big(\alpha x^\bullet_{0:k}+(1-\alpha)y^\bullet_{0:k}\big) &=  \Psi_k\big(\alpha x^\bullet_{0:k}+(1-\alpha)y^\bullet_{0:k},\sigma(t^n_k,\alpha x^k_{k}+(1-\alpha)y^k_{k})\big)\\
                                                 &\le  \Psi_k\big(\alpha x^\bullet_{0:k}+(1-\alpha)y^\bullet_{0:k},\alpha \sigma(t^n_k, x^k_{k})+(1-\alpha)\sigma(t^n_k,y^k_{k})O\big)\\
                                                 &\le \alpha \,\Psi_k\big(x^\bullet_{0:k},\sigma(t^n_k,x^k_{k})\big)+(1-\alpha) \Psi_k\big(y^\bullet_{0:k},\sigma(t^n_k,y^k_{k})O\big).\\
						&= \alpha \,\Psi_k\big(x^\bullet_{0:k},\sigma(t^n_k, x^k_{k})\big)+(1-\alpha) \Psi_k\big(y^\bullet_{0:k},\sigma(t^n_k,y^k_{k})\big)\\  
						&=\alpha\,\Phi_k\big(x^\bullet_{0:k}\big)+(1-\alpha)\Phi_k\big(y^\bullet_{0:k}\big).
\end{align*}
In particular, $x_0\mapsto \Phi_0(x^\bullet_{0:0})=\phi_0(x_0)$ is convex.
~$\hfill\cqfd$

\subsection{Proof of Proposition~\ref{prop:convexeul}}\label{sec:compar}
Let ${\cal F}_t=\sigma(X_0,W_s,s\in[0,t], {\cal N}_{\P})$, $\Phi_n:{\cal E}_n\to\R$ be convex with at most affine growth. We define the sequence $(\Phi_k)_{k=0,\ldots,n}$ by the backward induction equalities~\eqref{eq:Phik},~\eqref{eq:Psik}.
It is clear by backward induction that, when $\|X_0\|_1< +\infty$, then
\[
  \Phi_k(X^\bullet_{0:k}) = \E\big( \Phi_n(X^\bullet_{0:n})\,|\, {\cal F}_{t^n_k})\;\mbox{ for }\;k\in\{0,\cdots,n\},
\]
and, in particular,
\begin{equation}
   \E\,  \Phi_0(X_0) =\cdots =  \E\, \Phi_{k}(X^\bullet_{0:k})   =\cdots= \E\, \Phi_n(X^\bullet_{0:n}).\label{eq:egalesp}
 \end{equation}
We also set $\widetilde\Phi_n=\Phi_n$ and define $(\widetilde\Phi_k)_{1\le k\le n}$  by backward induction: for $k\in\{0, \ldots,n-1\}$ and $x^\bullet_{0:k}\in {\cal E}_k$ (with the  convention $x^0_0=x^1_0=\cdots=x^n_0=x_0$ when $k=0$),
\begin{align*}
\widetilde \Phi_k(x^\bullet_{0:k})  = \E\,\widetilde \Phi_{k+1} \left(x^\bullet_{0:k},B_k(t^n_{k+1}, x_k^{k+1},x_k^k)\right.&\left.+\widetilde \sigma(t^n_k,x_k^k)\int_{t^n_k}^{t^n_{k+1}}
                                     \widetilde K_2(t^n_{k+1},\hat s)dW_s,\cdots,\right.\\
&\left. \qquad B_k(t^n_n, x^n_k,x_k^k)+\widetilde \sigma(t^n_k,x_k^k)\int_{t^n_k}^{t^n_{k+1}} 
                                                                                       \widetilde  K_2(t^n_{n},\hat s)dW_s\right),
\end{align*}
where $B_k$ is still given by~\eqref{eq:B(t,x)}. 
Note that these equalities still hold with $\widetilde \Phi$, $\widetilde\sigma$ and $\widetilde K_2$ replaced by $\Phi$, $\sigma$ and $ K_2$ and that, similarly to~\eqref{eq:egalesp}, 
$\E[\widetilde \Phi_0(Y_0)]=\E[\widetilde \Phi_n(Y^\bullet_{0:n})]=\E[\Phi_n(Y^\bullet_{0:n})]$.

Let us check by backward induction that for each $k\in\{0,\cdots,n\}$, $\Phi_k\le\widetilde\Phi_k$. The induction hypothesis holds with equality for $k=n$. Let us assume that it holds at rank $k+1$. When both the Euler schemes are $K$-integrated, then by~\eqref{eq:sigmacKompar2}, 
\begin{equation}
   \int_{t^n_k}^{t^n_{k+1}} {\mathbf K}_2  {\mathbf K}_2^*(t^n_n,\cdots,t^n_{k+1},\hat s)ds\otimes \sigma\sigma^*(t^n_k,x_k^k)\le \int_{t^n_k}^{t^n_{k+1}}\widetilde {\mathbf K}_2\widetilde  {\mathbf K}_2^*(t^n_n,\cdots,t^n_{k+1},\hat s)ds\otimes \widetilde\sigma\widetilde\sigma^*(t^n_k,x_k^k)\label{usesigmacompar}.
\end{equation}
This inequality still holds as a consequence of~\eqref{eq:sigmacKompar} when both the Euler schemes are $K$-discrete since then $\hat s=t^n_k$ for $s\in[t^n_k,t^n_{k+1})$. The left-hand and right-hand sides of the above inequality
are the respective covariance matrices of the centered Gaussian vectors $$\left(\sigma(t^n_k,x_k^k)\int_{t^n_k}^{t^n_{k+1}} K_2(t^n_{k+1},\hat s)dW_s,\cdots,\sigma(t^n_k,x_k^k)\int_{t^n_k}^{t^n_{k+1}} K_2(t^n_{n},\hat s)dW_s\right)$$ and $\left(\widetilde \sigma(t^n_k,x_k^k)\int_{t^n_k}^{t^n_{k+1}}\widetilde  K_2(t^n_{k+1},\hat s)dW_s,\cdots,\widetilde \sigma(t^n_k,x_k^k)\int_{t^n_k}^{t^n_{k+1}}\widetilde  K_2(t^n_{n},\hat s)dW_s\right)$ so that the former is smaller than the latter in the convex order. By Proposition~\ref{prop:convexprop}, when  $\sigma$ (resp. $\widetilde \sigma$) satisfies the convexity assumption~\eqref{eq:cxsig}, then the function $\Phi_{k+1}$ (resp. $\widetilde \Phi_{k+1}$) is convex so that
\begin{align*}
&\E\,\Phi_{k+1} \left(x^\bullet_{0:k},B_k(t^n_{k+1}, x_k^{k+1},x_k^k)+\sigma(t^n_k,x_k^k)\int_{t^n_k}^{t^n_{k+1}}
                                                                                                                                                                                                                K_2(t^n_{k+1},\hat s)dW_s,\cdots,\right.\\
&\hskip 7,5 cm  \left. B_k(t^n_n, x^n_k,x_k^k)+\sigma(t^n_k,x_k^k)\int_{t^n_k}^{t^n_{k+1}} 
                                                                                                                                                                 K_2(t^n_{n},\hat s)dW_s\right)\\
&\le \E\,\Phi_{k+1} \left(x^\bullet_{0:k},B_k(t^n_{k+1}, x_k^{k+1},x_k^k)+\widetilde \sigma(t^n_k,x_k^k)\int_{t^n_k}^{t^n_{k+1}}
                 \widetilde K_2(t^n_{k+1},\hat s)dW_s,\cdots,\right.\\
&\hskip 7,5 cm  \left. B_k(t^n_n, x^n_k,x_k^k)+\widetilde \sigma(t^n_k,x_k^k)\int_{t^n_k}^{t^n_{k+1}}
                                                                                                                                                                        \widetilde  K_2(t^n_{n},\hat s)dW_s\right)
\end{align*}
with the left-hand side equal to $\Phi_{k}(x^\bullet_{0:k})$ and the right-hand side not greater than $\widetilde \Phi_{k}(x^\bullet_{0:k})$ by the induction hypothesis (resp. the same inequality holds with $\Phi_{k+1}$ replaced by $\widetilde \Phi_{k+1}$, the right-hand side being equal to $\widetilde \Phi_{k}(x^\bullet_{0:k})$ and the left-hand side not smaller than $\Phi_{k}(x^\bullet_{0:k})$ by the induction hypothesis). Hence the induction hypothesis holds at rank $k$.

For $k=0$, with the convexity of either $\Phi_0$ (when $\sigma$ satisfies~\eqref{eq:cxsig}) or $\widetilde \Phi_0$ (when $\widetilde \sigma$ satisfies~\eqref{eq:cxsig}) and the inequality $X_0\le_{cvx} Y_0$, it yields $\E[\Phi_0(X_0)]\le \E[\widetilde \Phi_0(Y_0)]$ which also writes $\E[\Phi_n(X^\bullet_{0:n})]\le \E[\Phi_n(Y^\bullet_{0:n})]$. We conclude with Lemma~\ref{onlylineargrowth}.\hfill $\Box$

\section{Increasing convex ordering for the Euler schemes ($d=q=1$)}
The comparison for the Euler schemes is again obtained as a consequence of the comparison for the companion processes $X^{\bullet}_{0:n}$ and $X^{\bullet}_{0:n}$ where the notation has been introduced at the beginning of Section~\ref{sec:propag}. We recall that, under~\eqref{sigcaff} and~\eqref{supteul}, when $\|X_0\|_p+\|Y_0\|_p< +\infty$ for some $p\ge 1$, 
 then $\|X^{\bullet}_{0:n}\|_p+\|Y^{\bullet}_{0:n}\|_p< +\infty$.
 \begin{Proposition}[Increasing convex ordering]
   \label{prop:convexeul1d} Assume~\eqref{sigcaff}, that either both functions $\R\ni x\mapsto b(t^n_k,x)$ and $\R\ni x\mapsto |\sigma(t^n_k,x)|$ are convex and non-decreasing for every $k\in\{0,\cdots,n-1\}$ or both functions $\R\ni x\mapsto \widetilde b(t^n_k,x)$ and $\R\ni x\mapsto |\widetilde \sigma(t^n_k,x)|$ are   convex and non-decreasing for every $k\in\{0,\cdots,n-1\}$.
Also assume that either both the Euler schemes are $K$-discrete and the comparison assumptions~\eqref{comparkb} between $(K_1,b)$ and $(\widetilde K_1,\widetilde b)$ and $\big(\mathcal{CK}^{disc}_2\sigma_{1d}\big)$ (see~\eqref{comparksig})
between $(K_2,\sigma)$ and $(\widetilde K_2,\widetilde \sigma)$ hold or both the Euler schemes are $K$-integrated and~\eqref{supteul},~\eqref{compkbs} and~$\big(\mathcal{CK}^{int}_2\sigma\big)$ hold. Last assume that $X_0\le_{icv} Y_0$ with $X_0,Y_0\in L^1(\P)$. Then, for every $n\ge 1$,  
$$
X^{\bullet}_{0:n}\le_{icv}Y^{\bullet}_{0:n}\mbox{ and, in particular, }(\bar X_{t^n_k})_{k=0,\ldots,n}\le_{icv}(\bar{Y}_{t^n_k})_{k=0,\ldots,n}.
$$ \end{Proposition}
\begin{Remark}
   For the conclusion to hold, it is enough that \begin{itemize}
     \item  when both schemes are $K$-integrated, \eqref{compkbs} holds for $(s_0,s_1,s_2)\in\{(t^n_k,t^n_{k+1},t^n_\ell):0\le k<\ell\le n\}$ and~$\big(\mathcal{CK}^{int}_2\sigma\big)$ holds for $(s_j,\cdots,s_1,s_0)$ of the form $(t^n_n,t^n_{n-1},\cdots,t^n_{n+1-j},t^n_{n-j})$ for $j\in\{0,\hdots,n\}$,
      \item  when both schemes are $K$-discrete, the comparison assumptions~\eqref{comparkb} between $(K_1,b)$ and $(\widetilde K_1,\widetilde b)$ and $\big(\mathcal{CK}^{int}_2\sigma_{1d}\big)$ between $(K_2, \sigma)$ and $(\widetilde K_2, \widetilde \sigma)$ hold for $(s,t)\in\{(t^n_k,t^n_\ell):0\le k<\ell\le n\}$.
   \end{itemize}
\end{Remark}
\subsection{Backward propagation of increasing convexity by the Euler scheme}\label{sbpci}

 \smallskip
Let $\Phi_n: {\cal E}_n\to \R$ be a function with $p$-polynomial growth for some $p\ge 1$. We define  the functions $\Phi_k: {\cal E}_k\to \R$, $k\in\{0, \cdots,n-1\}$ using some auxiliary functions $\widehat \Psi_k$ by the following backward induction:  for 
every $(x^\bullet_{0:k},u,v)\!\in {\cal E}_k\times \R\times\R$,
\begin{align}
 \widehat \Psi_k(x^\bullet_{0:k}, u,v) &= \E\,\Phi_{k+1} \left(x^\bullet_{0:k}, x_k^{k+1}+u\int_{t^n_k}^{t^n_{k+1}}K_1(t^n_{k+1},\hat s)ds +v \int_{t^n_k}^{t^n_{k+1}} K_2(t^n_{k+1},\hat s)dW_s,\right.\notag\\
& \left. \hskip 2cm  \cdots,x^n_k+u\int_{t^n_k}^{t^n_{k+1}}K_1(t^n_{n},\hat s)ds+v\int_{t^n_k}^{t^n_{k+1}} K_2(t^n_{n},\hat s)dW_s\right)\label{defpsiphi}
\end{align}
\begin{equation}
\mbox{ and }\forall x^\bullet_{0:k}\!\in {\cal E}_k,\;\Phi_k(x^\bullet_{0:k}) = \widehat \Psi_{k}\big(x^\bullet_{0:k},b(t^n_k,x_k^k), \sigma(t^n_k,x_k^k)\big),\label{defphipsi}
 \end{equation}
 with the convention $x_0^0=x^1_0=\cdots=x^n_0=x_0$ when $k=0$.
Note that under what concerns $(K_1,K_2,b, \sigma)$ in~\eqref{sigcaff} and~\eqref{supteul} (which appears in the assumptions of the proposition just below), the functions $x\mapsto b(t_k^n,x)$ and $x\mapsto \sigma(t_k^n,x)$  have at most affine growth and 
$$\left(\int_{t^n_k}^{t^n_{k+1}} K_2(t^n_{k+1},\hat s)dW_s,\cdots,\int_{t^n_k}^{t^n_{k+1}} K_2(t^n_{n},\hat s)dW_s\right)$$
is a Gaussian random vector and thus belongs to $L^p(\P)$. One easily deduces still  by backward induction that the random variable in the expectation in the right-hand side of~\eqref{defpsiphi} is integrable and both functions $\Psi_k$ and $\Phi_k$ have $p$-polynomial growth for every $k\in\{0,\cdots,n-1\}$.
 \begin{Proposition}[Backward propagation of increasing convexity]\label{prop:backinconv} Assume 
\begin{align*}
  \forall k\in\{0,\cdots,n-1\},\;&\int_0^{t^n_{k+1}} \left(K_1(t^n_{k+1},\hat s)+K_2(t^n_{k+1},\hat s)^2\right)ds+\sup_{x\in\R^d}\frac{|b(t^n_k,x)|+\|\sigma(t^n_{k},x)\|}{1+|x|}< +\infty,\end{align*}
and that both functions $\R\ni x\mapsto b(t^n_k,x)$ and $\R\ni x\mapsto |\sigma(t^n_k,x)|$ are   convex and non-decreasing for every $k\in\{0,\cdots,n-1\}$. When $\Phi_n:{\cal E}_n\to\R$ is convex with $p$-polynomial growth for some $p\ge 1$ and non-decreasing in each of its variables, then so are the  functions ${\cal E}_k\ni x^\bullet_{0:k}\mapsto \Phi_k(x^\bullet_{0:k})$, $k=0,\ldots,n$,  defined by~\eqref{defphipsi} and~\eqref{defpsiphi}.

In particular, if $\Phi_n(x^{\bullet}_{0:n}) = F(x_0,x_1^1,\cdots,x_n^n)$ with $F: (\R^d)^{n+1}\to \R$ convex with polynomial growth and non-decreasing in each of its variables, then convexity and monotonicity are transferred to the function $\R^d\ni x_{0}\mapsto \Phi_0(x_{0})= \E\big(F((\bar X^{x_0}_{t^n_k})_{k=0,\ldots,n})\big)$ where $(\bar X^{x_0}_{t^n_k})_{k=0,\ldots,n}$ denotes any of the two Euler schemes (either $K$-discrete or $K$-integrated) starting from $X_0=x_0$.
\end{Proposition}
\noindent {\bf Proof.} 
If $\Phi_{k+1}$ is convex w.r.t.  $x^\bullet_{0:k+1}$, then $\widehat\Psi_k$ is convex w.r.t. $(x^\bullet_{0:k}, u,v )$ and, by a reasoning similar to the above derivation of~\eqref{invrotbis} and~\eqref{eqcroisG}, satisfies
   \begin{align}
     \forall x^\bullet_{0:k}\in {\mathcal E}_k,\;\forall u,v\in\R,&\;\widehat\Psi_k(x^\bullet_{0:k}, u,v)=\widehat\Psi_k(x^\bullet_{0:k},u,|v|),\label{inrotbibis}\\
     &\mbox{and }\forall w\in\R\mbox{ s.t. }|v|\le |w|,\;\widehat\Psi_k(x^\bullet_{0:k}, u,v)\le \widehat\Psi_k(x^\bullet_{0:k},u,w).\label{eqcroisGbis}
\end{align}
If  $\Phi_{k+1}$ is non-decreasing in each of its variables, then $\widehat\Psi_k$ is non-decreasing  in $u\in\R$ since $K_1\ge 0$.
 Hence when $\Phi_{k+1}$ is convex and non-decreasing in each of its variables, then using that $x\mapsto b(t^n_k,x)$ and $x\mapsto |\sigma(t^n_k,x)|$ are non-decreasing, we deduce that $\Phi_k$ is non-decreasing in each of its variables. Moreover, using~\eqref{defphipsi} and~\eqref{inrotbibis} for the equalities, the monotonicity of $\widehat\Psi_k$ in its two last variables and the convexity of $x\mapsto|\sigma(t^n_k,x)|$ and $x\mapsto b(t^n_k,x)$ for the first inequality, then the convexity of $\widehat\Psi_k$ for the second inequality, we obtain that for $x^\bullet_{0:k},y^\bullet_{0:k}\in{\mathcal E}_k$ and $\alpha\in[0,1]$,
\begin{align*}
\Phi_k\big(\alpha x^\bullet_{0:k}&+(1-\alpha)y^\bullet_{0:k}\big) =  \widehat\Psi_k\big(\alpha x^\bullet_{0:k}+(1-\alpha)y^\bullet_{0:k},b(t^n_k,\alpha x^k_{k}+(1-\alpha)y^k_{k}),|\sigma(t^n_k,\alpha x^k_{k}+(1-\alpha)y^k_{k})|\big)\\&\le  \widehat\Psi_k\big(\alpha x^\bullet_{0:k}+(1-\alpha)y^\bullet_{0:k},\alpha b(t^n_k, x^k_{k})+(1-\alpha)b(t^n_k,y^k_{k}),\alpha|\sigma(t^n_k,x^k_{k})|+(1-\alpha)|\sigma(t^n_k,y^k_{k})|\big)\\
                                                 &\le  \alpha\widehat\Psi_k\big(x^\bullet_{0:k},b(t^n_k, x^k_{k}),|\sigma(t^n_k,x^k_{k})|\big)+(1-\alpha)\widehat\Psi_k\big(y^\bullet_{0:k},b(t^n_k,y^k_{k}),|\sigma(t^n_k,y^k_{k})|\big)\\
						&=\alpha\,\Phi_k\big(x^\bullet_{0:k}\big)+(1-\alpha)\Phi_k\big(y^\bullet_{0:k}\big).
\end{align*}

By backward induction, we conclude that when $\Phi_n$ is convex and non-decreasing in each of its variables, then so is $\Phi_k$ for each $k\in\{0,\hdots,n\}$.\hfill $\Box$
\begin{Remark}
   A simpler argument for the derivation of~\eqref{eqcroisGbis} in this more elementary framework is that, when $f$ is convex and $Z$ is a centered random variable, then 
 \[
 (x,v)\longmapsto \E\, f(x +vZ)\quad\mbox{is convex and non-decreasing  in $v$ on $\R_+$}
 \] 
 as it is bounded from below by its value  at  $0$ by Jensen's inequality.
\end{Remark}

 \subsection{Proof of Proposition~\ref{prop:convexeul1d}}\label{sec:compari}

 Let $\Phi_n:{\cal E}_n\to \R$ be a convex function with at most affine growth, non-decreasing in each of its variables. We are going to check that $\E[\Phi_n(X^\bullet_{0:n})]\le \E[\Phi_n(Y^\bullet_{0:n})]$ so that the conclusion follows from Lemma~\ref{onlylineargrowth}. We define the sequence $(\widehat\Psi_k,\Phi_k)_{k=0,\ldots,n-1}$ by the backward induction equalities~\eqref{defpsiphi} and~\eqref{defphipsi}.
 It is clear by backward induction that
\[
  \Phi_k(X^\bullet_{0:k}) = \E\big( \Phi_n(X^\bullet_{0:n})\,|\, {\cal F}_{t^n_k})\mbox{ for }k\in\{0,\cdots,n\}.
\]
In particular,
\begin{equation}
   \E\,  \Phi_0(X_0) =\cdots =  \E\, \Phi_{k}(X^\bullet_{0:k})   =\cdots= \E\, \Phi_n(X^\bullet_{0:n}).\label{egalesp2}
 \end{equation}

 We also set $\widetilde\Phi_n=\Phi_n$ and define $(\widetilde\Phi_k)_{k=0,\ldots,n-1}$  by backward induction: for $k\in\{0, \ldots,n-1\}$ and $x^\bullet_{0:k}\in {\cal E}_k$ (with $x^0_0=x^1_0=\cdots=x^n_0=x_0$),
\begin{align*}
\widetilde \Phi_k(x^\bullet_{0:k}) & = \E\,\widetilde \Phi_{k+1} \left(x^\bullet_{0:k},x_k^{k+1}+\widetilde b(t^n_k,x_k^k)\int_{t^n_k}^{t^n_{k+1}}\widetilde K_1(t^n_{k+1},\hat s)ds+\widetilde \sigma(t^n_k,x_k^k)\int_{t^n_k}^{t^n_{k+1}}
                                     \widetilde K_2(t^n_{k+1},\hat s)dW_s,\cdots,\right.\\
&\hskip 4.2 cm  \left. x_k^{n}+\widetilde b(t^n_k,x_k^k)\int_{t^n_k}^{t^n_{k+1}}\widetilde K_1(t^n_n,\hat s)ds+\widetilde \sigma(t^n_k,x_k^k)\int_{t^n_k}^{t^n_{k+1}} 
                                                                                                                                                                                                                                              \widetilde  K_2(t^n_{n},\hat s)dW_s\right).
\end{align*}
Note that these equalities still hold with $\widetilde \Phi$, $\widetilde b$, $\widetilde K_1$, $\widetilde\sigma$ and $\widetilde K_2$ replaced by $\Phi$, $b$, $K_1$, $\sigma$ and $K_2$ and that, similarly to~\eqref{egalesp2}, 
$\E[\widetilde \Phi_0(Y_0)]=\E[\widetilde \Phi_n(Y^\bullet_{0:n})]=\E[\Phi_n(Y^\bullet_{0:n})]$.

Let us check by backward induction that for each $k\in\{0,\cdots,n\}$, $\Phi_k\le\widetilde\Phi_k$. The induction hypothesis holds with equality for $k=n$. Let us assume that it holds at rank $k+1$. When both Euler schemes are $K$-integrated, then by~\eqref{eq:sigmacKompar2}, 
\begin{equation*}
\sigma^2(t^n_k,x_k^k)\int_{t^n_k}^{t^n_{k+1}} {\mathbf K}_2  {\mathbf K}_2^*(t^n_n,\cdots,t^n_{k+1},\hat s)ds\le \widetilde\sigma^2(t^n_k,x_k^k)\int_{t^n_k}^{t^n_{k+1}}\widetilde {\mathbf K}_2\widetilde  {\mathbf K}_2^*(t^n_n,\cdots,t^n_{k+1},\hat s)ds.
\end{equation*}
This inequality still holds as a consequence of~\eqref{comparksig} when both Euler schemes are $K$-discrete since then $\hat s=t^n_k$ for $s\in[t^n_k,t^n_{k+1})$. The left-hand side and the right-hand side are the respective covariance matrices of the centered Gaussian random vectors $$\left(\sigma(t^n_k,x_k^k)\int_{t^n_k}^{t^n_{k+1}} K_2(t^n_{k+1},\hat s)dW_s,\cdots,\sigma(t^n_k,x_k^k)\int_{t^n_k}^{t^n_{k+1}} K_2(t^n_{n},\hat s)dW_s\right)$$ and $\left(\widetilde \sigma(t^n_k,x_k^k)\int_{t^n_k}^{t^n_{k+1}}\widetilde  K_2(t^n_{k+1},\hat s)dW_s,\cdots,\widetilde \sigma(t^n_k,x_k^k)\int_{t^n_k}^{t^n_{k+1}}\widetilde  K_2(t^n_{n},\hat s)dW_s\right)$ so that the former is smaller than the latter in the convex order. Moreover, when both Euler schemes are $K$-integrated, by~\eqref{compkbs}, $b(t^n_k,x_k^{k})\int_{t^n_k}^{t^n_{k+1}}K_1(t_\ell^n,\hat s)ds \le \widetilde b(t^n_k,x_k^{k})\int_{t^n_k}^{t^n_{k+1}}\widetilde K_1(t_\ell^n,\hat s)ds$ for $\ell\in\{k+1,\cdots,n\}$. This inequality still holds as a consequence of~\eqref{comparkb} when both the Euler schemes are $K$-discrete. Hence $$\forall \ell\in\{k+1,\cdots,n\},\;x^\ell_k+b(t^n_k,x_k^{k})\int_{t^n_k}^{t^n_{k+1}}K_1(t_\ell^n,\hat s)ds \le x^\ell_k+\widetilde b(t^n_k,x_k^{k})\int_{t^n_k}^{t^n_{k+1}}\widetilde K_1(t_\ell^n,\hat s)ds.$$ By Proposition~\ref{prop:backinconv}, when the functions $\R\ni x\mapsto b(t,x)$ and $\R\ni x\mapsto |\sigma(t,x)|$ (resp. $\R\ni x\mapsto \widetilde b(t,x)$ and $\R\ni x\mapsto  |\widetilde\sigma(t,x)|$) are  both convex and non-decreasing for every $t\!\in [0,T]$, then the function $\Phi_{k+1}$ (resp. $\widetilde \Phi_{k+1}$) is convex and non-decreasing in each of its variables so that
\begin{align*}
&\E\,\Phi_{k+1} \left(x^\bullet_{0:k},x_k^{k+1}+ b(t^n_k,x_k^k)\int_{t^n_k}^{t^n_{k+1}} K_1(t^n_{k+1},\hat s)ds+ \sigma(t^n_k,x_k^k)\int_{t^n_k}^{t^n_{k+1}}
                                      K_2(t^n_{k+1},\hat s)dW_s,\cdots,\right.\\
&\hskip 4.2 cm  \left. x_k^{n}+ b(t^n_k,x_k^k)\int_{t^n_k}^{t^n_{k+1}} K_1(t^n_n,\hat s)ds+ \sigma(t^n_k,x_k^k)\int_{t^n_k}^{t^n_{k+1}}  K_2(t^n_{n},\hat s)dW_s\right)\\
&\le \E\,\Phi_{k+1} \left(x^\bullet_{0:k},x_k^{k+1}+ \widetilde b(t^n_k,x_k^k)\int_{t^n_k}^{t^n_{k+1}} \widetilde K_1(t^n_{k+1},\hat s)ds+ \widetilde \sigma(t^n_k,x_k^k)\int_{t^n_k}^{t^n_{k+1}}
                                     \widetilde K_2(t^n_{k+1},\hat s)dW_s,\cdots,\right.\\
&\hskip 4.5 cm  \left. x_k^{n}+\widetilde b(t^n_k,x_k^k)\int_{t^n_k}^{t^n_{k+1}}\widetilde K_1(t^n_n,\hat s)ds+\widetilde \sigma(t^n_k,x_k^k)\int_{t^n_k}^{t^n_{k+1}}\widetilde  K_2(t^n_{n},\hat s)dW_s\right)
\end{align*}
with the left-hand side equal to $\Phi_{k}(x^\bullet_{0:k})$ and the right-hand side not smaller than $\widetilde \Phi_{k}(x^\bullet_{0:k})$ by the induction hypothesis (resp. the same inequality holds with $\Phi_{k+1}$ replaced by $\widetilde \Phi_{k+1}$, the right-hand side being equal to $\widetilde \Phi_{k}(x^\bullet_{0:k})$ and the left-hand side not smaller than $\Phi_{k}(x^\bullet_{0:k})$ by the induction hypothesis). Hence the induction hypothesis holds at rank $k$.

For $k=0$, with the convexity and monotonicity of either $\Phi_0$ (when $\R\ni x\mapsto b(t,x)$ and $\R\ni x\mapsto |\sigma(t,x)|$ are convex and non-decreasing) or $\widetilde \Phi_0$ (when $\R\ni x\mapsto \widetilde b(t,x)$ and $\R\ni x\mapsto |\widetilde \sigma(t,x)|$ are convex and non-decreasing) and the inequality $X_0\le_{icv} Y_0$, it yields $\E[\Phi_0(X_0)]\le \E[\widetilde \Phi_0(Y_0)]$ which also writes $\E[\Phi_n(X^\bullet_{0:n})]\le \E[\Phi_n(Y^\bullet_{0:n})]$. \hfill $\Box$
 
\section{From discrete  to  continuous time}\label{sec:Cvgce}
\subsection{$L^p$-convergence of the affine interpolation of the Euler scheme}
To transfer our discrete time results on convex ordering to continuous time, we combine Theorems~\ref{thm:Eulercvgce2} and~\ref{thm:Eulercvgce1} with the ${\cal C}([0,T], \R^d)$-valued interpolation operator $i_n$ associated to the mesh $(t^n_k = \frac{kT}{n})_{0\le k\le n}$, already introduced in~\cite{Pag2016}: 
\[
i_n: (\R^{d})^{n+1}\ni x_{0:n} \longmapsto \Big(t \mapsto  \sum_{k=1}^n \mbox{\bf 1}_{[t^n_{k-1},t^n_k)}(t) \left(\frac{t^n_{k}-t}{t^n_{k}-t^n_{k-1}}x_{k-1} + \frac{t-t^n_{k-1}}{t^n_{k}-t^n_{k-1}} x_k\right)+\mbox{\bf 1}_{\{t=t_n^n\}}x_n\Big)\!.
\]
\begin{Corollary}\label{coro:cvinterpol}Let $p\in[1,+\infty)$ be such that $\|X_0\|_p< +\infty$ and suppose either that $(\bar X_{t^n_k})_{k=0,\ldots,n}$ is the $K$-integrated Euler scheme and the assumptions of Theorem~\ref{thm:Eulercvgce2} $(b)$ are satisfied or  $(\bar X_{t^n_k})_{k=0,\ldots,n}$ is the $K$-discrete Euler scheme and the assumptions of Theorem~\ref{thm:Eulercvgce1} $(b)$ are satisfied. Then
  \[
 \Big\|\big\| i_n\big((\bar X_{t^n_k})_{k=0,\ldots,n}\big) -X\big \|_{\sup} \Big\|_p\to 0\quad \mbox{ as }\quad n\to +\infty.
\]   
\end{Corollary}
\noindent {\bf Proof.} For $x_{0:n}\in(\R^{d})^{n+1}$, the function $i_n(x_{0:n})$ being continuous and piecewise affine with affinity breaks at times $t^n_k$ where it is equal to $x_k$, \begin{equation*}\| i_n(x_{0:n})\|_{\sup}= \max_{k=0,\ldots,n} |x_k|.
\end{equation*}
For $x_{0:n},y_{0:n+1}\in(\R^{d})^{n+1}$, since $i_n(x_{0:n})-i_n(y_{0:n})=i_n(x_{0:n}-y_{0:n})$, we deduce that $i_n$ is Lipschitz with constant $1$ from $(\R^d)^{n+1}$ to $\big({\cal C}([0,T], \R^d), \|\cdot\|_{\sup}\big)$:
\[
\| i_n(x_{0:n}) -i_n(y_{0:n}) \|_{\sup}= \max_{k=0,\ldots,n} |x_k-y_k|.
\]
 As a consequence, one has
 \[
\big\| i_n\big((\bar X_{t^n_k})_{k=0,\ldots,n}\big) -i_n\big((X_{t^n_k})_{k=0,\ldots,n}\big)\big \|_{\sup} = \max_{k=0,\ldots,n} |\bar X_{t^n_k}-X_{t^n_k}| 
 \]  
Let $p>0$ be such that $\|X_0\|_p+\|Y_0\|_p< +\infty$. Calling upon either Theorem~\ref{thm:Eulercvgce1} or Theorem~\ref{thm:Eulercvgce2}, we deduce that 
 \[
 \Big\|\big\| i_n\big((\bar X_{t^n_k})_{k=0,\ldots,n}\big) -i_n\big((X_{t^n_k})_{k=0,\ldots,n}\big)\big \|_{\sup} \Big\|_p=\big \| \max_{k=0,\ldots,n} |\bar X_{t^n_k}-X_{t^n_k}| \big\|_p\to 0\quad \mbox{ as }\quad n\to +\infty.
 \] Now,  the uniform continuity modulus $w(\xi, \delta)=\sup_{0\le s\le t\le (s+\delta)\wedge T}|\xi(t)-\xi(s)|$ of a function $\xi\!\in {\cal C}([0,T], \R^d)$ converges to $0$ as $\delta\to 0$. One easily checks that 
 \[
 \big \| X-i_n\big( (X_{t^n_k})_{k=0,\ldots, n}\big)\big\|_{\sup} \le  w\big(X, \tfrac Tn\big)\le 2\big \| X\big\|_{\sup}.
 \]
As, according to the last statement in Theorem~\ref{prop:exunvolt}, $\big\| \|X\|_{\sup}\big\|_p <+\infty$, it follows by dominated convergence that  
\[
\Big\|  \big \| X-i_n\big( (X_{t^n_k})_{k=0,\ldots, n}\big)\big\|_{\sup}\Big\|_p\to 0\quad \mbox{ as }\quad n\to +\infty.
\] 
Combining these two convergences yields  
\[
 \hskip 4cm\Big\|\big\| i_n\big((\bar X_{t^n_k})_{k=0,\ldots,n}\big) -X\big \|_{\sup} \Big\|_p\to 0\quad \mbox{ as }\quad n\to +\infty.\hskip 4cm \Box
\]

\subsection{Proofs of Theorem~\ref{thm:main} and Theorem~\ref{thm:ci}}
\noindent {\bf Proof of Theorem~\ref{thm:main}.}
In this proof, $(\bar X_{t^n_k})_{k=0:n}$ denotes the $K$-discrete Euler scheme when $(K_1,K_2)$ and $(\widetilde K_1,\widetilde K_2)$ satisfy $(\underline{\cal K}^{cont})$ and the $K$-integrated Euler scheme when this condition is replaced by the weaker condition $({\cal K}^{cont})$ up to a reinforcement of the comparison assumptions between $(K_2,\sigma)$ and $(\widetilde K_2,\widetilde \sigma)$ (resp. $(K_1,b)$ and $(\widetilde K_1,\widetilde b)$).

\smallskip
\noindent $(a)$ By Corollary~\ref{coro:cvinterpol}, the  $1$-Wasserstein distance ${\cal W}_1\big( i_n\big((\bar X_{t^n_k})_{k=0:n}\big),X\big)$ between the law of $i_n\big((\bar X_{t^n_k})_{k=0:n}\big)$ and that of $X$ goes to $0$ as $n\to\infty$. For each ${\rm Lip}(F)$-Lipschitz functional $F: {\cal C}([0,T], \R^d)\to \R$, since $\left|\E\, F\circ  i_n\big((\bar X_{t^n_k})_{k=0:n}\big) \big)-\E\, F(X)\right|\le {\rm Lip}(F){\cal W}_1\big( i_n\big((\bar X_{t^n_k})_{k=0:n}\big),X\big)$, we have
$$\E\, F\circ  i_n\big((\bar X_{t^n_k})_{k=0:n}\big) \big)\to\E\, F(X)\mbox{ as }n\to\infty.$$ In the same way, $\E\, F\circ  i_n\big((\bar Y^n_{t^n_k})_{k=0:n}\big) \big)\to\E\, F(Y)$.

Let now $F: {\cal C}([0,T], \R^d)\to \R$ be  a Lipschitz convex functional.  The function $F_n: (\R^d)^{n+1} \to \R$ defined by $F_n(x_{0:n})= F\big (i_n(x_{0:n})\big)$ is ${\rm Lip}(F)$-Lipschitz since $i_n$ is $1$-Lipschitz, convex as the composition of a convex function with a linear application and such that $
\E\, F\circ  i_n\big((\bar X_{t^n_k})_{k=0:n}\big) \big)=\E\, F_n \big((\bar X_{t^n_k})_{k=0:n}\big)$ and $
\E\, F\circ  i_n\big((\bar Y^n_{t^n_k})_{k=0:n}\big) \big)=\E\, F_n \big((\bar Y^n_{t^n_k})_{k=0:n}\big)$. It follows from Proposition~\ref{prop:convexeul} that, under the assumptions made on the coefficients $(K_1,b,K_2,\sigma,\widetilde K_1,\widetilde b,\widetilde K_2,\widetilde \sigma)$ and the initial conditions $(X_0,Y_0)$ that, for every $n\ge 1$, 
\[
\E\, F_n\big((\bar X_{t^n_k})_{k=0:n} \big)\le  \E\, F_n\big((\bar{Y}^n_{t^n_k})_{k=0,\ldots,n}\big).
\]
Consequently, 
\[
\E\, F(X)  = \lim_n \E\, F_n\big((\bar X_{t^n_k})_{k=0:n} \big)\le \lim_n \E\, F_n\big((\bar Y^n_{t^n_k})_{k=0,\ldots,n}\big)   = \E\, F(Y).
\]
By Lemma~\ref{onlylipschitz}, we conclude that $X\preceq_{cvx} Y$.

\smallskip
\noindent $(b)$ One easily checks that if $d=1$ and the convex functional $F$ is non-decreasing for the pointwise partial order on ${\cal C}([0,T], \R)$, then $F_n $  is non-decreasing for the componentwise partial order on $\R^{n+1}$. One concludes by the same reasoning as in the proof of $(a)$ with Proposition~\ref{prop:convexeul} replaced by Proposition~\ref{prop:convexeul1d}.

\smallskip
\noindent $(c)$ Setting $(\widehat X_t)_{t\in[0,T]}=(-X_t)_{t\in[0,T]}$, $(\widehat Y_t)_{t\in[0,T]}=(-Y_t)_{t\in[0,T]}$ and $(\widehat W_t)_{t\in[0,T]}=(-W_t)_{t\in[0,T]}$, we have that $(\widehat W_t)_{t\in[0,T]}$ is a Brownian motion independent from $\widehat X_0$ and 
\begin{align*}
   \widehat X_t&=\widehat X_0+\int_0^t K_1(t,s) (-b(s,-\widehat X_s))ds+\int_0^t K_2(t,s)\sigma(s,-\widehat X_s)d\widehat W_s,\quad t\in[0,T],\\
  \widehat Y_t&=\widehat Y_0+\int_0^t \widetilde K_1(t,s) (-\widetilde b(s,-\widehat Y_s))ds+\int_0^t \widetilde K_2(t,s)\widetilde \sigma(s,-\widehat Y_s)d\widehat W_s,\quad t\in[0,T].\end{align*}
The inequality $X_0\preceq_{dcv} Y_0$ implies $\widehat X_0\preceq_{icv}\widehat Y_0$ (see the reasoning just below concerning the functionals $F$ and $G$). When $x\mapsto -b(t,x)$ and $x\mapsto |\sigma(t,x)|$ (resp. $x\mapsto -\widetilde b(t,x)$ and $x\mapsto |\widetilde \sigma(t,x)|$) are convex and non-increasing, then $x\mapsto -b(t,-x)$ and $x\mapsto |\sigma(t,-x)|$ (resp. $x\mapsto -\widetilde b(t,-x)$ and $x\mapsto |\widetilde \sigma(t,-x)|$) are convex and non-decreasing. Moreover~\eqref{comparkbb} implies that $$\forall\, 0\le s<t\le T,\;\forall x\in\R,\; K_1(t,s)(-b(s,-x))\le \widetilde K_1(t,s)(-\widetilde b(s,-x))$$  and~\eqref{compkbbs} implies that $$\forall\, 0\le s_0<s_1\le s_2\le T,\;\forall x\in\R,\;-b(s_0,-x)\int_{s_0}^{s_1}K_1(s_2,s)ds\le-\widetilde b(s_0,-x)\int_{s_0}^{s_1}\widetilde K_1(s_2,s)ds.$$
Let $F: {\cal C}([0,T], \R)\to \R$ be l.s.c., convex and non-increasing for the pointwise partial order on continuous functions. Then the functional $G$ defined by $G((x_t)_{t\in[0,T]})=F(-(x_t)_{t\in[0,T]})$ is also  l.s.c., convex but non-decreasing for the pointwise partial order. By $(b)$, $\E G(\widehat X)\le \E G(\widehat Y)$, which also writes $\E F(X)\le \E F(Y)$. This completes the proof.\hfill $\Box$

\bigskip
\noindent {\bf Proof of Theorem~\ref{thm:ci}.}
In this proof, $(\bar X^{x_0}_{t^n_k})_{k=0:n} $ denotes the $K$-integrated Euler scheme starting from $x_0\in\R^d$.

\smallskip
\noindent $(a)$ Let $p\in[1,+\infty)$ denote the polynomial growth order of the convex functional $F: {\cal C}([0,T], \R^d)\to \R$. We know from~\cite{Lucchetti2006} Lemma 2.1.1 that $F$ is then continuous for the sup-norm. By Corollary~\ref{coro:cvinterpol}, the  $p$-Wasserstein distance ${\cal W}_p\big( i_n\big((\bar X^{x_0}_{t^n_k})_{k=0:n}\big),X^{x_0}\big)$ between the law of $i_n\big((\bar X^{x_0}_{t^n_k})_{k=0:n}\big)$ and that of the solution $X^{x_0}$ to~\eqref{eq:Volterra} starting from $x_0\in\R^d$ goes to $0$ as $n\to\infty$. Hence 
$$
\E\, F\circ  i_n\big((\bar X^{x_0}_{t^n_k})_{k=0:n}\big) \big)\to\E\, F(X^{x_0})\mbox{ as }n\to\infty.
$$

 Since $(\R^d)^{n+1}\ni x_{0:n}\mapsto F_n(x_{0:n})=F(i_n(x_{0:n}))$ is convex with $p$-polynomial growth, by Proposition~\ref{prop:convexprop}, $\R^d\ni x_0\mapsto \E\, F\circ  i_n\big((\bar X^{x_0}_{t^n_k})_{k=0:n}\big) \big)$ is convex. We conclude that $\R^d\ni x_0\mapsto \E\, F(X^{x_0})$ is convex as the pointwise limit of convex functions.

\smallskip
\noindent $(b)$ If $d=1$ and the convex functional $F$ is non-decreasing for the pointwise partial order on ${\cal C}([0,T], \R)$, then $F_n $  is non-decreasing for the componentwise partial order on $\R^{n+1}$. One concludes by the same reasoning as in the proof of $(a)$ with Proposition~\ref{prop:convexprop} replaced by Proposition~\ref{prop:backinconv}.

\smallskip
\noindent $(c)$ The conclusion follows from $(b)$ combined with the change of sign argument used in the proof of Theorem~\ref{thm:main} $(c)$.\hfill $\Box$

\small
\bibliographystyle{plain} 
\bibliography{propconvex}
\appendix
\small 
\section{Matrices, covariances and random vectors}\label{sec:matrices}
\begin{Lemma} \label{lem:matrices} Let $A, B\!\in \R^{d\times q}$. Then 
\[
AA^* = BB^* \Longleftrightarrow \exists\, O\!\in {\cal O}(q) \; \mbox{ such that }\; A = BO.
\]
Moreover, $A=\sqrt{AA^*}Q$ where $Q\in\R^{d\times q}$ is such that $QQ^*$ is the orthogonal projection on ${\rm Im}\sqrt{AA^*}$ or equivalently $Q$ has $\mbox{\rm dim}({\rm Im}\sqrt{AA^*})$ columns which form an orthonormal basis of $\mbox{\rm Im}\sqrt{AA^*}$ and its other columns vanish. 
 \end{Lemma}      
 \begin{Corollary}  Let $A, B\!\in \R^{d\times q}$.  If $AA^*=BB^*$, then for any  radially symmetric $\R^q$-valued random vector $Z$, $AZ\stackrel{(d)}{=} BZ$ (such is the case for ${\cal N}(0,I_q)$-distributed random vectors).  
\end{Corollary}
   The proof of the above corollary is obvious since, by definition of radial symmetry, 
   $OZ\stackrel{(d)}{=}Z$ for every $O\!\in {\cal O}(q)$. Lemma~\ref{lem:matrices} is not needed when $Z\sim{\mathcal N}_q(0,I_q)$, since then $AZ\sim{\mathcal N}_d(0,AA^*)$.

The proof below is inspired by similar (though less general) results from~\cite{JP}.

\medskip
\noindent {\bf Proof of lemma~\ref{lem:matrices}.} $\Leftarrow$: If $A=BO$ with $O\in{\mathcal O}(q)$, then $AA^*=BOO^*B^*=BB^*$.
                               
\smallskip \noindent   $\Rightarrow$: In order to prove (simultaneously)   the necessary condition and the equality $A=\sqrt{AA^*}Q$ in the three cases $d=q$, $d<q$ and $d>q$, let us characterize the matrices  $Q\!\in\R^{d\times q}$ such that $QQ^*$ is the orthogonal projection on ${\rm Im }AA^*$.

                               Since clearly ${\rm Im }AA^*\subset{\rm Im }A$ and ${\rm Ker }AA^*={\rm Ker }A^*=({\rm Im }A)^\perp$ so that $$\mbox{dim Im}AA^*=d-\mbox{dim Ker}AA^*=d-\mbox{dim }({\rm Im}A)^\perp=\mbox{dim Im}A\le d\wedge q,$$ ${\rm Im}A={\rm Im}AA^*={\rm Im}\sqrt{AA^*}$ where the second equality is obtained by replacing $A$ by $\sqrt{AA^*}$.
                               \\
When $Q$ has $\mbox{dim Im}\sqrt{AA^*}$ columns which form an orthonormal basis of $\mbox{Im}\sqrt{AA^*}$ and its other columns vanish, it is clear that $QQ^*$ is the orthogonal projection on $\mbox{Im }\sqrt{AA^*}$. Let us conversely suppose that $Q\in\R^{d\times q}$ with columns $(Q_j)_{1\le j\le q}$ is such that $QQ^*$ is the orthogonal projection on ${\rm Im}\sqrt{AA^*}$. Since ${\rm Im} Q={\rm Im} QQ^*={\rm Im}\sqrt{AA^*}$, the vectors $Q_j$ belong to ${\rm Im}\sqrt{AA^*}$ and $Q_j=QQ^*Q_j$. Hence for $j\in\{1,\cdots,q\}$, $$|Q_j|^2=Q_j^*QQ^*Q_j=\sum_{i=1}^q (Q_i^*Q_j)^2=|Q_j|^2+\sum_{\stackrel{i=1}{i\neq j}}^q (Q_i^*Q_j)^2,$$ which ensures that the vectors $Q_j$ are orthogonal. Hence $|Q_j|^2Q_j=QQ^*Q_j=Q_j$ so that either $|Q_j|=0$ or $|Q_j|=1$. Therefore $Q$ has exactly $\mbox{dim Im}\sqrt{AA^*}$ orthonormal columns which belong to $\mbox{Im}\sqrt{AA^*}$ and its other columns vanish.

 \smallskip                        
\noindent{\em Case $d=q$}. By the singular value decomposition, $A=UDV$ and $B=\bar U \bar D\bar V$ for $U,V,\bar U,\bar V\in{\mathcal O}(q)$ and $D,\bar D\in\R^{q\times q}$ diagonal with non-negative diagonal elements. Then $AA^*=UD^2U^*$ and $\sqrt{AA^*}=UDU^*$ so that $A=\sqrt{AA^*}UV$. Denoting by $P\in\R^{q\times q}$ the matrix of the orthogonal projection on ${\rm Im} \sqrt{AA^*}$ in the canonical basis of $\R^q$, one has $\sqrt{AA^*}=\sqrt{AA^*}P$ since $({\rm Im }\sqrt{AA^*})^\perp={\rm Ker} \sqrt{AA^*}$. Hence $A=\sqrt{AA^*}PUV$ where $PUV(PUV)^*=P$. In the same way, $B=\sqrt{BB^*}\bar U\bar V$. When $AA^*=BB^*$, we conclude that $A=BO$ with $O=\bar V^*\bar U^*UV\in{\mathcal O}(q)$.
                               
 \smallskip 
 \noindent{\em Case $d<q$}. We define  $\widetilde A,\widetilde B\in \R^{q\times q}$ by 
$$
(\widetilde A_{ij},\widetilde B_{ij})=\begin{cases}(A_{ij},B_{ij})\mbox{ for }i=1:d,j=1:q\\
  (0,0)\mbox{ for }i=d+1:q,j=1:q\end{cases}.
  $$
  When $AA^*=BB^*$, $\widetilde A\widetilde A^*=\widetilde B\widetilde B^*$, so that, by the equality of dimensions step, there exists $O\in{\mathcal O}(q)$ such that $\widetilde A =\widetilde B O$ and, by restriction of this matrix equality to the $d$ first lines, $A=BO$. The equality of dimensions step also ensures that $\widetilde A=\sqrt{\widetilde A\widetilde A^*}\widetilde Q$ with $\widetilde Q\in\R^{q\times q}$ such that $\widetilde Q\widetilde Q^*$ is the orthogonal projection on ${\rm Im}\sqrt{\widetilde A\widetilde A^*}$.  Since $\sqrt{\widetilde A\widetilde A^*}=({\mathbf 1}_{\{i\vee j\le d\}}\sqrt{AA^*}_{ij})_{1\le i,j\le q}$, we deduce that $A=\sqrt{AA^*}Q$ where $Q=(\widetilde Q_{ij})_{1\le i\le d,1\le j\le q}$. By the above characterization, $\widetilde Q$ has $\mbox{dim Im}\widetilde A=\mbox{dim Im}A$ columns which form an orthonormal basis of $\mbox{Im}\widetilde A$ and, by definition of $\widetilde A$ have vanishing coordinates with indices larger than $d$ so that the corresponding columns of $Q$ form an orthonormal basis of ${\rm Im}\sqrt{AA^*}$. The other columns of $\widetilde Q$ and therefore of $Q$ vanish. Hence by the characterization, $Q$ is such that $QQ^*$ is the orthogonal projection on ${\rm Im}\sqrt{AA^*}$.

 \smallskip 
 \noindent{\em Case $d>q$}. The equality $AA^*=BB^*$ implies that $\mbox{Im} B=\mbox{Im} \sqrt{AA^*}$. Let $Q\in\R^{d\times q}$ be such that $QQ^*$ is the orthogonal projection on $\mbox{Im}\sqrt{AA^*}$. Then $\sqrt{AA^*}=QQ^*\sqrt{AA^*}$, $A=QQ^*A$, $B=QQ^*B$ and $Q^*AA^*Q=Q^*BB^*Q$ so that, by the equality of dimensions step, $Q^*A=Q^*BO$ for some $O\in{\mathcal O}(q)$. Left-multiplying by $Q$ the last equality leads to $A=BO$. Since $Q^*\sqrt{AA^*}QQ^*\sqrt{AA^*}Q=Q^*AA^*Q$ and $Q^*\sqrt{AA^*}Q$ is symmetric, $Q^*\sqrt{AA^*}Q=\sqrt{Q^*AA^*Q}$. By the equal dimensions step this implies that $Q^*A=Q^*\sqrt{AA^*}QV$ for some $V\in{\mathcal O}(q)$ and, by left-multiplication by $Q$,  $A=\sqrt{AA^*}QV$ where $QVV^*Q^*=QQ^*$ is the orthogonal projection on $\mbox{Im}\sqrt{AA^*}$.~$\hfill\cqfd$
\begin{Lemma}[Kronecker product]\label{lem:Hadamard} Let $S\!\in {\cal S}^+(d_1)$ and $T\!\in {\cal S}^+(d_2)$, $d_1, d_2\!\in \N$. Then (with obvious notations) the Kronecker product  $S\otimes T := [S_{ij}T]_{1\le i,j\le d_1} \!\in {\cal S}^+(d_1\times d_2)$. When $S$ and $T$ are both positive definite, then so is $S\otimes T$.
\end{Lemma}

\noindent {\bf Proof.} First note that $S\otimes T $ is clearly symmetric by construction. Let $u\in \R^{(d_1\times d_2)}$ be defined as column vector and for $i=1:d_1$, let $U_i=(u_{(i-1)d_1+k})_{k=1:d_2}\in\R^{d_2}$. Then
\[
u^*S\otimes T u = \sum_{i,j=1}^{d_1} S_{ij} U^*_{i} \,T\, U_{j}
\]
The matrix $T$ can be diagonalized as $T= P\Delta P^*$ where $P\!\in {\cal O}(d_2)$ and $\Delta = {\rm Diag}(\delta_{1:d_2})$, $\delta_k\ge 0$, $k=1,\ldots,d_2$ Consequently,
\[
u^*(S\otimes T) u = \sum_{i,j=1}^{d_1} S_{ij} (P^*U_{i})^*\Delta (P^*U_{j})=  \sum_{k=1}^{d_2}\delta_k \sum_{i,j =1}^{d_1}S_{ij} (P^*U_{i})_k^* (P^*U_{j})_k \ge 0 
\]
since the $\delta_k$ are non-negative and $S$ is positive semi-definite.\hfill $\Box$
\small
\section{First elements of proof of Theorem~\ref{prop:exunvolt}}\label{app:eu}
The aim of this appendix is to establish the  existence and uniqueness of pathwise continuous solutions to the Volterra equation~\eqref{eq:Volterra} when $X_0\!\in L^0(\P)$. The $L^p$ pathwise regularity when $X_0\in L^p(\P)$ for some $p>0$, will be established in Appendix \ref{app:A''}.

Taking advantage  of the structure of our stochastic Volterra equation with distinct {kernels, the assumptions ensuring existence and uniqueness of pathwise continuous (even $a$-H\"older continuous for some small enough $a>0$) solutions}   in~\cite[Theorem~3.3]{ZhangXi2010} {
when $X_0\!\in \cap_{r>0} L^r(\P)$} boil down to check that $K_1$ and $K_2^2$ satisfy assumptions $H'_1$, $H2$-$H3$-$H4$ with the notations of~\cite{ZhangXi2010}. Using again its notations, this corresponds to  $g(t) =X_0$ and $\kappa_1 =\kappa_2= \kappa:=K_1+K_2^2$ and
\[
\forall\, 0\le s<t<t',\quad\lambda(t',t,s) = \sum_{i=1,2} |K_i(t',s)-K_i(t,s)|^i.
\]
Thus, condition  $({\cal K}^{int}_{\beta})$ is a criterion for  $\kappa$ to lie in the class ${\cal H}_{>1}$ as requested in assumption $H2$ of~\cite[Theorem~3.3]{ZhangXi2010} {(this is mentioned in~\cite{ZhangXi2010}).}
Condition $({\cal K}^{cont}_{\theta})$ ensures that $H4$ is satisfied. 
Hence, when $X_0\!\in \cap_{r>0} L^r(\P)$, {all the assumptions of~\cite[Theorem~3.3]{ZhangXi2010}} are satisfied (the additional  H\"older condition on the random adapted  function $g$ is here empty since here $g(t)=X_0$). So~\cite[Theorem~3.3]{ZhangXi2010} directly applies.

\smallskip
 {Consequently, we may start from the  fact that, under  the assumptions of our theorem on the kernels and the coefficients of Equation~\eqref{eq:Volterra} and if $X_0\!\in L^{\infty}(\P)\subset \cap_{r>0}L^r(\P)$, this equation has a unique  pathwise  continuous solution.}

\smallskip 
$\rhd$ {\em Existence}. We now suppose that $X_0\!\in L^0(\Omega,{\cal F}_0,\P)$. For every $k\!\in \N^*$, let  $A_k= \{|X_0|<k\}$ and $X^{(k)}$ denote the   $({\cal F}_t)$-adapted continuous unique solution to~\eqref{eq:Volterra} starting from $X_0^{(k)} = X_0\mbox{\bf 1}_{A_k}\!\in L^{\infty}(\Omega, {\cal F}_0, \P)$. We set $X=\sum_{k\in\N^*}X^{(k)}\mbox{\bf 1}_{A_k\setminus A_{k-1}}$ where $A_0=\emptyset$. Let $t\in[0,T]$. Since
$$\int_0^tK_2(t,s)^2\|\sigma(s,X_s)\|^2ds=\sum_{k\in \N^*}\mbox{\bf 1}_{A_k\setminus A_{k-1}}\int_0^tK_2(t,s)^2\|\sigma(s,X^{(k)}_s)\|^2ds<+\infty,\;\P\mbox{-}a.s.,$$
the stochastic integral $\int_0^tK_2(t,s)\sigma(s,X_s)dW_s$ makes sense. Moreover, since stochastic integrals with respect to an $({\cal F}_s)$-Brownian motion commute with $\F_0$-measurable random variables, for each $k\in\N^*$, \begin{align*}
   {\bf 1}_{A_k\setminus A_{k-1}}&\int_0^tK_2(t,s)\sigma(s,X_s)dW_s=\int_0^tK_2(t,s){\bf 1}_{A_k\setminus A_{k-1}}\sigma(s,X_s)dW_s\\&=\int_0^tK_2(t,s){\bf 1}_{A_k\setminus A_{k-1}}\sigma(s,X^{(k)}_s)dW_s={\bf 1}_{A_k\setminus A_{k-1}}\int_0^tK_2(t,s)\sigma(s,X^{(k)}_s)dW_s,\;\P\mbox{-}a.s..
\end{align*}
Hence \begin{align*}
        \int_0^tK_2(t,s)\sigma(s,X_s)dW_s=\sum_{k\in \N^*}{\bf 1}_{A_k\setminus A_{k-1}}\int_0^tK_2(t,s)\sigma(s,X^{(k)}_s)dW_s,\;\P\mbox{-}a.s..
\end{align*}
Since, one clearly has $\int_0^tK_1(t,s)b(s,X_s)ds=\sum_{k\in \N^*}{\bf 1}_{A_k\setminus A_{k-1}}\int_0^tK_1(t,s)b(s,X^{(k)}_s)ds,\;\P\mbox{-}a.s.$, we deduce that
\begin{align*}
   X_t&=\sum_{k\in \N^*}{\bf 1}_{A_k\setminus A_{k-1}}\left(X^{(k)}_0+\int_0^tK_1(t,s)b(s,X^{(k)}_s)ds+\int_0^tK_2(t,s)\sigma(s,X^{(k)}_s)dW_s\right)\\
&=X_0+\int_0^tK_1(t,s)b(s,X_s)ds+\int_0^tK_2(t,s)\sigma(s,X_s)dW_s,\;\P\mbox{-}a.s..
\end{align*}
with the sum over $k\in\N^*$ providing a continuous modification of the right-hand side.

\smallskip
$\rhd$  {\em Uniqueness}.  Let $X$ and $\widetilde X$ denote two $({\cal F}_t)$-adapted continuous  solutions to~\eqref{eq:Volterra} starting from the same initial condition $X_0\!\in L^0(\Omega,{\cal F}_0,\P)$ and let $X^0$ denote the solution starting from $0$. For $k\in \N^*$, $\mbox{\bf 1}_{A_k} X +  \mbox{\bf 1}_{A_k^c} X^0$ and $\mbox{\bf 1}_{A_k} \widetilde X +  \mbox{\bf 1}_{A_k^c} X^0$ are both solutions to~\eqref{eq:Volterra} starting from $X_0\mbox{\bf 1}_{A_k}\!\in L^{\infty}(\Omega,{\cal F}_0, \P)$. By uniqueness, $\P(\{X=\widetilde X\}\cap A_k)=\P(A_k)$. Letting $k\to\infty$, we conclude that $
\P(X=\widetilde X)=1$. %
Hence uniqueness holds starting from any random variable lying in $L^0(\Omega, {\cal F}_0, \P)$.

\section{Properties  of the  $K$-integrated Euler scheme and $L^p$-convergence at fixed times}
\label{app:B} 
\medskip {$\rhd$ {\em Preliminaries}. As a first preliminary, we will establish   the following lemma  which provides a control in $L^p(\P)$ of the increments of  general Lebesgue or  stochastic integrals involving the kernels $K_i$. It  will be used several times  throughout the next appendices.
 \begin{Lemma} \label{lem:intHcont} Let $(H_t)_{t\in[0,T]}$ (resp. $(\widetilde H_t)_{t\in[0,T]}$) be an $({\cal F}_t)_{t\in[0,T]}$-progressively measurable process
 having values in $\R^d$ (resp. $\mathbb{M}_{d,q}(\R)$) such that $\displaystyle \sup_{t\in [0,T]}\|H_t\|_r+\sup_{t\in [0,T]}\|\widetilde H_t\|_r<+\infty$ for some $r\ge 2$. Assume that the kernels $K_i$ satisfy $({\cal K}^{int}_{\beta})$ for some $\beta>1$ and $({\cal K}^{cont}_{\theta})$ for some $\theta\in (0,1]$.
 Let $0\le T_0<T_1\le T$. Then
both processes 
 \[
[T_0,T] \ni t\longmapsto \int_{T_0}^{t\wedge T_1} K_1(t,s)H_sds\quad \mbox{ and }\quad [T_0,T] \ni t \longmapsto \int_{T_0}^{t\wedge T_1} K_2(t,s)\widetilde H_sdW_s
 \]
are $\theta\wedge\frac{\beta-1}{2\beta}$-H\" older from $[T_0,T]$ to $L^r(\P)$. To be more precise, there exists a real constant $C= C_{r,T,K_1,K_2}>0$ such that 
\begin{align*}
\Big\| \int_{T_0}^{t\wedge T_1} K_1(t,u)H_udu-\int_{T_0}^{s\wedge T_1} K_1(s,u)H_udu \Big\|_r + \Big \| \int_{T_0}^{t\wedge T_1} K_2(t,u)&\widetilde H_udW_u-\int_{T_0}^{s\wedge T_1} K_2(s,u)\widetilde H_udW_u \Big\|_r \\
& \le  C\sup_{t\in [0,T]}\|H_t\|_r\vee\|\widetilde H_t\|_r |t-s|^{\theta\wedge\frac{\beta-1}{2\beta}}.
\end{align*}
If   furthermore $(\widehat {\cal K}^{cont}_{\widehat \theta})$ holds true, one can replace  $\frac{\beta-1}{2\beta}$ by $\widehat \theta$ {\em mutatis mutandis} in the aboves claims.
\end{Lemma}
}

\noindent {{\bf Proof.} Let $s,\, t\!\in [T_0,T]$ with $s<t$. By using successively the Burkholder-Davis-Gundy, generalized Minkowski and H\"older (with exponents $(\beta ,\frac{\beta}{\beta-1})$) inequalities, one gets since $\frac r2 \ge 1$,
\begin{align*}
&\left\|\int_{T_0}^{t\wedge T_1} K_2(t,u)\widetilde H_udW_u- \int_{T_0}^{s\wedge T_1} K_2(s,u)\widetilde H_udW_u\right\|_r^2\\&\phantom{\|\int_{T_0}^{t\wedge T_1} K_2(t,u)}\le (C_r^{BDG})^2\left( \left\| \int_{s\wedge T_1}^{t\wedge T_1} K_2(t,u)^2\widetilde H_u^2du\right\|_{\frac r2}+ \left\| \int_{T_0}^{s\wedge T_1} \big(K_2(t,u)-K_2(s,u)\big)^2\widetilde H_u^2du\right\|_{\frac r2}\right)\\
& \phantom{\|\int_{T_0}^{t\wedge T_1} K_2(t,u)}\le (C_r^{BDG})^2\left( \int_s^t K_2(t,u)^2\|\widetilde H_u\|_r^2du+  \int_{T_0}^s \big(K_2(t,u)-K_2(s,u)\big)^2\|\widetilde H_u\|_r^2du\right)\\
&\phantom{\|\int_{T_0}^{t\wedge T_1} K_2(t,u)}\le  (C_r^{BDG})^2\sup_{u\in [T_0,T]}\|\widetilde H_u\|_r^2\Big( 
\sup_{t\in [0,T]} \Big(\int_0^t K_2(t,u)^{2\beta}  \Big)^{\frac{1}{\beta}}(t-s)^{\frac{\beta-1}{\beta}}\!+ \!\eta(t-s)^2 \Big).
\end{align*}
Hence there exists a real constant $\kappa= \kappa_{\eta, \beta,T}$ only depending on $\eta$, $\beta$ and $T$  such that
\[
\left\|\int_{T_0}^{t\wedge T_1} K_2(t,u)\widetilde H_udW_u-\int_{T_0}^{s\wedge T_1} K_2(s,u)\widetilde H_udW_u\right\|_r \le \kappa\,C_r^{BDG}\sup_{u\in [T_0,T]}\|\widetilde H_u\|_r (t-s)^{\theta\wedge \frac{\beta-1}{2\beta}}
\]
owing to~\eqref{eq:contK} and~\eqref{eq:contKtilde}. The variant under $(\widehat {\cal K}^{cont}_{\widehat \theta})$  is straightforward since $\int_{s}^t K_2(t,u)^2du\le \widehat \eta(t-s)^2$ for $0\le s\le t$. One proceeds likewise for  $t\mapsto \int_{T_0}^{t\wedge T_1} K_1(t,u)H_udu$ (only using the generalized Minkowski inequality). \hfill$\Box$}

Throughout the remaining of this section, we assume that the assumptions of Theorem~\ref{thm:Eulercvgce2} are in force and that $X_0\in L^p(\P)$ for the value of $p$ under consideration.

First, we define to alleviate notations $\displaystyle\varphi_{a}(t) = \Big(\int_0^t K_1(t,s)^{a}ds\Big)^{1/a}$ and   $\displaystyle\psi_{a}(t) =\Big(\int_0^t K_2(t,s)^{a}ds\Big)^{1/a}$, $\varphi^*_a = \sup_{t\in[0,T]} \varphi_a(t)$ and $\psi^*_a$ likewise for $a >0$. Note  that $ \varphi^*_1+\psi_2^*<+\infty$ owing to Assumption~$({\cal K}^{int}_{\beta})$ and H\"older's inequality since $\beta >1$.

\medskip 
\noindent {\sc Property~{\bf 1}.} \label{pp:1} {\em $L^p$-integrability and pathwise continuity of $\bar X_t$, $p>0$}. {First we prove that $\bar X_{t^n_k}\!\in L^p(\P)$ by induction on $k$.  If $\bar X_{t^n_{\ell}}\!\in L^p(\P)$ for $\ell=0, \ldots,k-1$, it follows easily from equation~\eqref{eq:Euler2} that $\bar X_{t^n_k}\!\in L^p(\P)$ since both $b$ and $\sigma$ have linear growth on their space variable uniformly w.r.t. their time variable as a consequence of~\eqref{eq:HolLipbsig}  and  $\bar X_{t^n_{\ell-1}}$ and the Gaussian Wiener integral $\displaystyle \int_{t^n_{\ell-1}}^{t^n_\ell}K_2(t^n_k,s)dW_s$ are independent. The induction relies on the Minkowski inequality when $p\ge 1$ and on its subadditive counterpart for $\|\cdot\|_p^p$ when $p\!\in (0,1)$. }

For the continuous time   $K$-integrated Euler scheme, one derives likewise  from~\eqref{eq:Eulergen2bis} and Lemma \ref{lem:intHcont} that, as $\bar X_{t_\ell}\!\in L^p(\P)$ for every $\ell=0,\ldots,k$, one has $\sup_{t\in[t_k,t_{k+1}]}\|\bar X_t\|_p<+\infty$. As a consequence, $\sup_{t\in[0,T]}\|\bar X_t\|_p<+\infty$.

  \smallskip
 \noindent 
 As for the pathwise continuity, first note that   the functions $\left(\int_{t^n_{\ell}}^{t^n_{\ell+1}}  K_1(t,s)ds\right)_{t\in[t^n_{\ell+1},T]}$, $\ell=0, \ldots,n-1$,  are continuous owing to~$({\cal K}^{cont}_{\theta})$.  So are       $\left(\int_{t^n_\ell}^t K_1(t,s)ds\right)_{t\in[t^n_\ell,t^n_{\ell+1}]}$, $\ell=0,\ldots,n-1$, owing to~$({\cal K}^{int}_{\beta})$  and $({\cal K}^{cont}_{\theta})$.

\medskip Applying the above Lemma~\ref{lem:intHcont} with  $\widetilde H\equiv 1$, $T_0=t^n_k$, $T_1=t^n_{k+1}$ with $k\in\{0,\cdots,n-1\}$ and $r>\frac{1}{\theta\wedge\frac{\beta-1}{2\beta}}$ yields that the processes $\left(\int_{t^n_k}^{t\wedge t^n_{k+1}} K_2(t,s)dW_s\right)_{t\in[t^n_k,T]}$ have a continuous modification (null at $t^n_k$) owing to  Kolmogorov's criterion (see~\cite[Theorem 2.1, p.26, $3^{rd}$ edition]{RevuzYor}). Then one concludes that the $K$-integrated Euler scheme has a continuous modification.

\medskip 
\noindent {\sc Property~{\bf ~2}.} \label{pp:2} {\em Moment control of $\bar X_t$, $p\ge 2$}.   $\rhd$
Let $\rho \!\in (0,1)$ and $A>0$ be two real numbers to be specified later on. Under Assumption~$({\cal LH}_{\gamma})$, $b$ and $\sigma $ have  linear growth on $x$ uniformly in $t\!\in [0, T]$ i.e. there exists a real constant
$C^{(0)}= C_{b,\sigma,T}>0$  such that, for every $t\!\in [0,T]$ and every $x\!\in \R$, 
\begin{equation}\label{eq:Lineargrowth}
|b(t,x)| +\|\sigma(t,x)\|\le C^{(0)} (1+|x|^2)^{1/2} \le C^{(0)}(1+|x|).
\end{equation}
Set $\displaystyle \bar f_p(t) = \sup_{0\le s \le t} \|\bar X_s\|_p$. We already know that $\bar f_p(T)<\infty$ and are now going to derive an estimation not depending on the number $n$ of time-steps. It follows from the above equation and the generalized Minkowski inequality that 
\begin{align*}
\Big\| \int_0^t K_1(t,s)b(s,\bar X_{\underline s}) ds \Big\|_p& \le \int_0^t K_1(t,s) \big\| b(s,\bar X_{\underline s})\|_pds\\
&\le  C^{(0)} \int_0^t K_1(t,s) \big( 1+\|\bar X_{\underline s}\|_p \big)ds\\
& \le C^{(0)}\Big(\varphi_1(t) +   \int_0^t K_1(t,s)  \|\bar X_{\underline s}\|_p  ds\Big).
\end{align*}
 As   $\underline s \le s$, $\|\bar X_{\underline s}\|_p\le \bar f_p(s)$ and the above  inequality implies 
\[
\bar f_p(t) \le \|X_0\|_p +C^{(0)}\Big(\varphi_1(t) + \int_0^t K_1(t,s) \bar f_p(s)ds\Big) + \sup_{s\le t}G_p(s),
\]
where $ \displaystyle G_p(t) = \Big\|\int_0^t K_2(t,s)  \sigma(s,\bar X_{\underline s})dW_s  \Big\|_p$.  {Set $\rho =  \frac 12 (1+ \frac {1}{\beta})=\frac{\beta+1}{2\beta}\!\in (0,1)$ since  $\beta >1$.} As $\bar f_p$ is non-decreasing, one derives
 \begin{align}
 \nonumber \bar f_p(t) & \le \bar f_p(0) + C^{(0)}\Big(\varphi_1(t) +   \bar f_p^{\rho}(t) \int_0^t K_1(t,s) \bar f_p(s)^{1-\rho} ds\Big)+\sup_{s\le t}G_p(s)\\
 \nonumber   & \le \bar f_p(0) +C^{(0)}  \varphi_1(t) + C^{(0)} \bar f_p^{\rho}(t)  {\varphi_{\frac{2\beta}{\beta+1}}(t) }\Big( \int_0^t \bar f_p(s) ds\Big)^{1-\rho}+\sup_{s\le t}G_p(s)\\
 \label{eq:L^pbound-b} & =  \bar f_p(0) + C^{(0)}\varphi_1(t) + C^{(0)}   {\varphi_{\frac{2\beta}{\beta+1}}(t) }\frac{ \bar f_p^{\rho}(t)}{A} \cdot A \Big(\int_0^t \bar f_p(s) ds\Big)^{1-\rho} + \sup_{s\le t}G_p(s),
 \end{align}
where we used H\"older's inequality in the second line with conjugate exponents {$(\frac{2\beta}{\beta +1}, \frac{2\beta}{\beta-1})$}.  Applying now Young's inequality with the same conjugate exponents yields 
\begin{equation}\label{eq:fbarp}
\bar f_p(t) \le \bar f_p(0) +C^{(0)}\varphi_1(t) +C^{(0)} { \varphi_{\frac{2\beta}{\beta+1}}(t)}\left(\frac{\rho}{A^{1/\rho}} \bar f_p(t)+(1-\rho) A^{1/(1-\rho)}\int_0^t \bar f_p(s)ds  \right)+\sup_{s\le t}G_p(s).
\end{equation}
 
 We proceed likewise with the Volterra pseudo-diffusion term, up to the use of the $L^p$-Burkholder-Davis-Gundy ({\em BDG}) inequality. First we get  for every $t\!\in [0,T]$,  
  \begin{align*}
  G_p(t)& \le C^{BDG}_p \left\|\int_0^t K_2(t,s)^2 \sigma(s,\bar X_{\underline s})^2 ds   \right\|^{1/2}_{p/2}.
  \end{align*}
  
  As $p\ge 2$, it follows from the linear growth property for $\sigma$ and both  regular and generalized Minkowski inequalities that
   \begin{align*}
   \nonumber G^2_p(t)& = (\widetilde C^{(0)})^2 \left(\int_0^t K_2(t,s)^2 ds + \Big\|\int_0^t K_2(t,s)^2|\bar X_{\underline s}|^2  ds \Big\|_{p/2} \right) \\
   \nonumber  &\le  (\widetilde C^{(0)})^2\left(\int_0^t K_2(t,s)^2 ds + \int_0^t K_2(t,s)^2 \big\|\bar X_{\underline s}\big\|^2_{p}  ds \right)
       \end{align*}
 where $\widetilde C^{(0)}= \widetilde C^{(0)}_{b,\sigma,p}:= C^{BDG}_pC^{(0)}$. Using that $\sqrt{a+b}\le \sqrt{a}+\sqrt{b}$, $a$, $b\ge 0$, yields
   \begin{align*}
    \nonumber  G_p(t) &\le \widetilde C^{(0)} \left[  \left(\int_0^t K_2(t,s)^2 ds \right)^{1/2} +   \left(\int_0^t K_2(t,s)^2 \big\|\bar X_{\underline s}\big\|^2_{p}ds\right)^{1/2}\right]\\
    \nonumber   &   \le  \widetilde C^{(0)}\psi_2(t)  +  \widetilde  C^{(0)} \left(\int_0^t K_2(t,s)^2 \bar f_p(s)^2d s\right)^{1/2}.
    \end{align*}
   Then, using again  that $\bar f_p$ is non-decreasing, it follows, still with $\rho= \frac 12(1+\frac {1}{\beta})$ that 
     \begin{align*}
  G_p(t)&\le  \widetilde C^{(0)} \psi_{2}(t) + \widetilde C^{(0)} \bar f_p^{\rho}(t)  \left(\int_0^t K_2(t,s)^2   \bar f_p(s)^{2(1-\rho)} ds\right)^{1/2}\\
  	    & \le \widetilde C^{(0)}  \psi_{2}(t) +\widetilde C^{(0)} \psi_{2\beta}(t)  \bar f_p^{\rho}(t)\Big(\int_0^t \bar f_p(s)^{\frac{2(1-\rho)\beta}{\beta-1}} ds\Big)^{\frac{\beta-1}{2\beta}},
      \end{align*}
where we applied H\" older's inequality with conjugate exponents $(\beta, \frac{\beta}{\beta-1})$.  Noting that $\frac{2(1-\rho)\beta}{\beta-1}= 1$, applying Young's inequality with conjugate exponents $\big(1/\rho,1/(1-\rho)\big)$ and introducing $A>0$ like for the drift term, we obtain
\begin{equation}\label{eq:L^pbound-sigma}
G_p(t) \le   \widetilde C^{(0)}\psi_{2}(t) + \widetilde C^{(0)} \psi_{2\beta}(t) \left( \rho \frac{\bar f_p(t)}{A^{1/\rho}} +(1-\rho) A^{1/(1-\rho)}\Big(\int_0^t \bar f_p(s) ds  \Big) \right).
\end{equation} 

Taking advantage of the fact that $\bar f_p$ is non-decreasing, we derive that
\[
\sup_{s\le t} G_p(s) \le   \widetilde C^{(0)}\!\sup_{0\le s\le t}\!\psi_{2}(s) + \widetilde C^{(0)} \!\sup_{0\le s\le t}\!\psi_{2\beta}(s) \left( \rho \frac{\bar f_p(t)}{A^{1/\rho}} +(1-\rho) A^{1/(1-\rho)}\Big(\int_0^t \bar f_p(s)  ds  \Big) \right).
\]
 {Plugging this inequality in~\eqref{eq:fbarp} yields}
\[
\bar f_p(t) \le \bar f_p(0) + C^{(1)}  +C^{(2)}   \left( \rho \frac{\bar f_p(t)}{A^{1/\rho}} +(1-\rho) A^{1/(1-\rho)} \int_0^t \bar f_p(s) ds \right)
\]
 where  
 \[
 C^{(1)}= C^{(1)}_{K_1,K_2,\beta,b,\sigma, T} = C^{(0)}\varphi^*_1 \,+ \,\widetilde   C^{(0)} \psi^*_{2}
\]
 and
 \[
 C^{(2)}=C^{(2)}_{K_1,K_2,p,\beta,b,\sigma, T} = C^{(0)}  {\varphi_{\frac{2\beta}{\beta+1}}^*} \,+\, \widetilde   C^{(0)}\psi_{2\beta}^*.
\] 
Both  positive real constants $C^{(1)}$  and $C^{(2)}$ are  finite owing to $({\cal K}^{int}_{\beta})$.
 
 \medskip We choose now $A$ large enough, namely $ A= (2 \rho\, C^{(2)}  )^{\rho}$, so that $1- \rho \,C^{(2)} A^{-1/\rho}= \frac 12< 1$. Consequently, for every $t\!\in [0,T]$, 
 \[
 \bar f_p(t)\le 2\bar f_p(0) +  2 C^{(1)} + 2 C^{(2)}  (1-\rho) A^{1/(1-\rho)} \int_0^t \bar f_p(s) ds.
 \] 
 Since $\bar f_p(0)= \|X_0\|_p$, one concludes by Gr\"onwall's lemma that 
 \[
 \bar f_p(t) \le 2\big(\|X_0\|_p+  C^{(1)}  \big) e^{C^{(3)}  t}
 \]
with $C^{(3)}= C^{(3)}_{K_1,K_2,p,\beta,b,\sigma,T} =  2 C^{(2)}   (1-\rho) A^{1/(1-\rho)}$. As a conclusion to this step, we get
\[
\sup_{t\in [0,T]}\|\bar X_t\|_p=\bar f_p(T) \le C^{(4)}(1+\|X_0\|_p)
\]
 with $C^{(4)} = 2\max \big( 1,C^{(1)}) e^{C^{(3)}T} \geq 2$. Since, according to \cite[Theorem 3.1]{ZhangXi2010}, $\sup_{t\in [0,T]}\|X_t\|_p<\infty$, by the same reasoning applied to the Volterra process, we also derive
 \[
\sup_{t\in [0,T]}\|X_t\|_p\le C^{(4)}(1+\|X_0\|_p).
\]
 \smallskip 
 \noindent {\sc Property~{\bf ~3}.}\label{pp:3}  {\em  Control of the increments $\|\bar X_t-\bar X_s\|_p$ of the $K$-integrated Euler scheme~I (for $p\ge 2$)}. 
{Our aim in this step  is to prove   that there exists a positive real constant  $C_{p,T}$ not depending on $n$ such that 
  \[
\forall\, s,\, t\!\in [0,T], \quad  \|\bar X_t-\bar X_s\|_p\le C_{p,T}(1+\|X_0\|_p) |t-s|^{\theta\wedge\widehat\theta}.
\]
We recall the definition~\eqref{eq:Eulergen2}
 \begin{equation*}
\bar X_t =X_0+\int_0^t K_1(t,s)b(\underline s,\bar X_{\underline s})ds +\int_0^t K_2(t,s)\sigma(\underline s,\bar X_{\underline s})dW_s,\quad t\in[0,T].
 \end{equation*}
 which strongly suggests to call upon Lemma~\ref{lem:intHcont} (with $r=p$) 
 with $H_s = b(\underline s, \bar X_{\underline s})$ and $\widetilde H_s =  \sigma(\underline s, \bar X_{\underline s})$.
First note that
\[
 \big\| b(\underline{s},\bar X_{\underline{s}})\big\|_p \le C^{(0)}(1+\|\bar X_{\underline s}\|_p) \le C^{(0)}\big(1+C^{(4)}(1+\|\bar X_0\|_p) \big)\le 2C^{(0)}C^{(4)}\big(1+  \|\bar X_0\|_p\big) 
\]  
owing to~\eqref{eq:Lineargrowth} and Property~{\bf 2} and the fact that $C^{(4)}\geq 1$. One shows likewise that,
for every $u\!\in [0,T]$, 
\[
\big\| \sigma(\underline{s},\bar X_{\underline{s}}) \big\|_{p}\le 2C^{(0)}C^{(4)}(1+\|\bar X_0\|_p).
\]
Consequently, it follows from Lemma~\ref{lem:intHcont}  under $\widehat{\cal K}^{cont}_{\widehat \theta}$ that
\begin{align*}
 \|\bar X_t-\bar X_s\|_p & \le C_{r,T,K_1,K_2} \Big(\sup_{t\in [0,T]}\big\| b(\underline{u},\bar X_{\underline{u}})\big\|_p  +\sup_{t\in [0,T]}\big\| \sigma(\underline{u},\bar X_{\underline{u}})\big\|_p \Big)(t-s)^{\theta\wedge \widehat \theta}\\
 &\le C^{(5)}(1+\|\bar X_0\|_p)(t-s)^{\theta\wedge \widehat \theta},
\end{align*}
where $C^{(5)} = 4C_{p,T,K_1,K_2}C^{(0)}C^{(4)}$ with $ C_{r,T,K_1,K_2}$ coming from Lemma~\ref{lem:intHcont}. Otherwise, if only  $({\cal K}^{int}_{\beta})$  is in force, replace $\widehat \theta$ by $\frac {\beta-1}{2\beta}$.
}

\medskip
 \noindent{\sc Property~{\bf ~4}.} \label{pp:4}{\em 
   Rate of convergence of $\|X_t-\bar X_t\|_p$  at fixed time $t$ for $p>p_{eu}$}.   Have in mind that  $p_{eu}= \frac{1}{\theta}\vee \frac{2\beta}{\beta-1}$ and that for such values of $p$ we know from a careful reading of the proof of~\cite[Theorem~3.3]{ZhangXi2010} that  Equation~\eqref{eq:Volterra} has  a unique  adapted  pathwise continuous solution $(X_t)_{t\in [0,T]}$. Then, for every  $t\!\in [0,T]$, 
 \[
 X_t-\bar X _t = \int_0^t K_1(t,s) \big(b(s,X_s)-b(\underline s, \bar X_{\underline s})\big)ds + \int_0^t K_2(t,s) \big(\sigma(s,X_s)-\sigma(\underline s, \bar X_{\underline s})\big)dW_s\quad\P\mbox{-}a.s.
 \]
We denote by $\bar A^n_1(t)$ and $\bar A^n_2(t)$ the two terms of the sum  in the right-hand side  of the above equation. Set 
$$
g_p(t) = \displaystyle \sup_{s\in [0,t]} \|X_s-\bar X_s\|_p, \quad t\!\in [0,T].
$$
 This non-decreasing function is finite owing to Theorem~\ref{prop:exunvolt} and Property~{\bf 1} since $X_0\!\in L^p(\P)$. We straightforwardly deduce that 
\begin{equation}
   g_p(t)\le \sup_{s\le t}\big(\| \bar A^n_1(s)\|_p + \| \bar A^n_2(s)\|_p\big).\label{eq:majogpt}
\end{equation}
Let $[b]_{H,L}$ be the mixed H\"older-Lipschitz coefficient such that, for every $s$, $t\!\in [0,T]$ and every $x$, $y\!\in \R$,
\[
|b(s,x)-b(t,y)| \le [b]_{H,L}\big(|t-s|^{\gamma}(1+|x|+|y|) +|x-y|\big).
\] 
The coefficient $[\sigma]_{H,L}$ is defined likewise.  

\smallskip
As for $\bar A^n_1(t)$, the   generalized Minkowski inequality implies (with in mind the notation $h= \frac Tn$ for the time step) 
  \begin{align*}
\big \| \bar A^n_1(t)\big\|_p &\le  \int_0^t  \Big\|  K_1(t,s) \big|b(s,X_s) - b(\underline s, \bar X_{\underline s}) \big|\Big\|_p ds\\
 & \le  [b]_{H,L}\int_0^t  K_1(t,s)\Big( | s-\underline s|^{\gamma}\big(1+\|X_s\|_p +\|\bar X_s\|_p) + \|X_s -\bar X_{\underline s}\|_p  \Big)ds\\
 & \le  [b]_{H,L} \left( \big(1+ 2 \sup_{t\in [0,T]} (\|X_s\|_p \vee \|\bar X_s\|_p)\big)\varphi^*_1  \, h^{\gamma} + \int_0^t K_1(t,s) \|X_s -\bar X_{\underline s}\|_p   ds  \right)\\
 & \le 3\,C^{(4)} \varphi^*_1 [b]_{H,L}   (1+\|X_0\|_p) h^{\gamma} + [b]_{H,L} \left(\int_0^t K_1(t,s) g_p(s)   ds +\int_0^t K_1(t,s) \|\bar X_s -\bar X_{\underline s}\|_p   ds \right)  
 \end{align*}
 where we used that  $ \|X_s -\bar X_{\underline s}\|_p\le  g_p(s)+ \|\bar X_{ s} -\bar X_{\underline s}\|_p $.
 
\smallskip
Following the strategy developed in~\eqref{eq:L^pbound-b}, we get  for every $A>0$ and  {$\rho= \frac 12(1+\frac{1}{\beta})\!\in (0, 1)$},
\begin{align}
\nonumber \big\| \bar A^n_1(t)\big\|_p &\le    [b]_{H,L} \varphi_{1/\rho}(t) \left(\frac{\rho}{A^{1/\rho}}g_p(t) +(1-\rho) A^{1/(1-\rho)}\int_0^t g_p(s)ds  \right)\\
\label{eq:barAn1t} & \qquad+ [b]_{H,L}\int_0^t K_1(t,s) \|\bar X_s -\bar X_{\underline s}\|_p   ds +  3\,C^{(4)} \varphi^*_1 [b]_{H,L}(1+\|X_0\|_p)h^{\gamma}.
\end{align}
 
\smallskip
As for $\bar A^n_2(t)$, combining   the {\em BDG} inequality and both  regular and generalized Minkowski inequalities implies
  \begin{align*}
 \| \bar A^n_2(t)\|^2_p & \le (C^{BDG}_p)^2 \left\|  \int_0^t  K_2(t,s)^2 \big\|\sigma(s,X_s) - \sigma (\underline s, \bar X_{\underline s})\big\|^2   ds\right\|_{p/2}\\
 & \le  (C^{BDG}_p)^2[\sigma]^2_{H,L}\int_0^t  K_2(t,s)^2\Big( | s-\underline s|^{\gamma}\big(1+\|X_s\|_p +\|\bar X_s\|_p) + \|X_s -\bar X_{\underline s}\|_p \Big)^2ds\\
 & \le  2  (C^{BDG}_p)^2[\sigma]^2_{H,L} \left( \big(3C^{(4)}\big)^2 \psi_2(t)^2(1+\|X_0\|_p)^2   h^{2\gamma} + \int_0^t K_2(t,s)^2 \|X_s -\bar X_{\underline s}\|^2_p   ds  \right).
 \end{align*}
 
  \medskip
  Following this time the strategy developed to  establish~\eqref{eq:L^pbound-sigma} and  setting again $\rho =  \frac 12 (1+ \frac {1}{\beta})\!\in (0,1)$, we derive
 \begin{align}\nonumber
  \big \| \bar A^n_2(t)\big\|_p & \le \sqrt{2}C^{BDG}_p[\sigma]_{H,L}\left[ \psi_{2\beta}(t) \left( \rho \frac{g_p(t)}{A^{1/\rho}} +(1-\rho) A^{1/(1-\rho)} \int_0^t g_p(s) ds \right)\right.\\
 \label{eq:barAn2t} & \hskip3cm\left.+ \Big(\int_0^t K_2(t,s)^2 \|\bar X_s -\bar X_{\underline s}\|^2_p   ds\Big)^{1/2} + 3C^{(4)}\psi_2(t)(1+\|X_0\|_p)h^{\gamma} \right].
  \end{align}

  \smallskip
  Combining~\eqref{eq:barAn1t} and~\eqref{eq:barAn2t} yields 
 \begin{equation}\label{eq:g_pGronwall1}
  \big \| \bar A^n_1(t)\big\|_p +  \big \| \bar A^n_2(t)\big\|_p   \le C^{(6)} \left( \rho \frac{g_p(t)}{A^{1/\rho}} +(1-\rho) A^{1/(1-\rho)} \int_0^t g_p(s) ds \right) +  R(t)
  \end{equation}
  with $\displaystyle C^{(6)} = C^{(6)}_{K_1,K_2, \beta, b,\sigma, p } = [b]_{H,L} {\varphi^*_{\frac{2\beta}{\beta+1}}}+\sqrt{2} C^{BDG}_p[\sigma]_{H,L} \psi^*_{2\beta}$ and 
 \begin{equation}\label{eq:R(t)}
R(t) =  C^{(7)}\left((1+\|X_0\|_p) h^{\gamma}  +\int_0^t K_1(t,u) \|\bar X_u -\bar X_{\underline u}\|_p   du  
+ \Big(\int_0^t K_2(t,u)^2 \|\bar X_u -\bar X_{\underline u}\|^2_p   du\Big)^{1/2} \right)
  \end{equation}
with
$$
 C^{(7)}=C^{(7)}_{K_1,K_2, \beta, b,\sigma, p } = \displaystyle \max\big(3C^{(4)}[b]_{H,L} \varphi_1^*+3\sqrt{2} C^{(4)}C^{BDG}_p [\sigma]_{H,L}\psi^*_2, [b]_{H,L}, \sqrt{2}C^{BDG}_p  [\sigma]_{H,L}\big).
$$ 

Then setting $A= \big(2 \rho \,C^{(6)}\big)^{\rho}$ yields for every $t\!\in [0,T]$ (since $g_p(0)=0$) 
 \begin{align}\label{eq:g_pGronwall2}
\big \| \bar A^n_1(t)\big\|_p +  \big \| \bar A^n_2(t)\big\|_p \le \frac {g_p(t)}2+ C^{(6)}(1-\rho)A^{1/(1-\rho)}\int_0^t g_p(s)ds +R(t).
\end{align}  
Consequently, it follows from~\eqref{eq:majogpt} and  Gr\"onwall's lemma that  
 \begin{equation*}
g_p(t)  \le 2 C^{(6)}(1-\rho)A^{1/(1-\rho)}\int_0^t g_p(s)ds +2 \,\sup_{s\le t} R(s) \le 2 \,e^{C^{(8)} t} \sup_{s\le t} R(s)
  \end{equation*}
  with $C^{(8)} = C^{(8)} _{K_1,K_2, \beta, b,\sigma, p } =2 C^{(6)}(1-\rho)A^{1/(1-\rho)}$. In particular 
  \begin{equation}\label{eq:g_pGronwall3}
 \sup_{t\in [0,T]} \|X_t-\bar X_t\|_p= g_p(T) \le 2\,e^{C^{(8)} T} \sup_{t\in [0,T]} R(t).
  \end{equation}
   Now, coming back to $R(t)$, we have that
 \begin{align*}
  \Big(\int_0^t K_2(t,s)^2 \|\bar X_s -\bar X_{\underline s}\|^2_p   ds\Big)^{1/2}  & \le \psi_{2}^*\sup_{s\in [0,t]} \big\| \bar X_s-\bar X_{\underline s}\big\|_p 
 \end{align*}
 whereas, obviously,
 \[
 \int_0^t K_1(t,s) \|\bar X_s -\bar X_{\underline s}\|_p ds \le \varphi_{1}^*\sup_{s\in [0,t]} \big\| \bar X_s-\bar X_{\underline s}\big\|_p. 
 \]
Hence, for every $t\in [0,T]$, 
  \[
 \sup_{s\le t} R(s) \le C^{(7)}(1+\|X_0\|_p)h^{\gamma}+  C^{(7)}\big(\varphi^*_1+\psi^*_{2}\big)  \sup_{s\in [0,t]} \big\| \bar X_s-\bar X_{\underline s}\big\|_p.
  \]
  
By combining these bounds with Property~{\bf 3},   one concludes  that, under $(\widehat {\cal K}^{cont}_{\widehat \theta})$
  \[
 \forall\, t\!\in [0,T], \quad  \|X_t-\bar X_t\|_p \le C^{(9)}(1+\|X_0\|_p) \left(h^{\gamma}+h^{\theta\wedge\widehat\theta}\right)
  \]
with $C^{(9)}= C^{(9)} _{K_1,K_2, \beta, b,\sigma, p,T } = 2\,e^{C^{(8)} T}C^{(7)}\max\left(1,(\varphi^*_1+\psi^*_{2}\big)C^{(5)}\right)$.
  
Note that if only $({\cal K}^{int}_{\beta})$ is in force, then the final bound holds {\em mutatis mutandis} with $\frac{\beta-1}{2\beta}$ instead of $\widehat \theta$.

\section{Splitting lemma and representation of the   solution of  the Volterra equation {and its Euler schemes}} \label{app:B'} 
Theorem~\ref{thm:Eulercvgce1} is a revisited version of~\cite[Theorem~2.2]{RiTaYa2020}. Surprisingly, in its original formulation the error bounds only hold for starting values $X_0\!\in L^p(\P)$ when  $p$ is large enough depending on some integrability characteristics of the kernels. This unusual restriction can be circumvented thanks to  the lemma below applied to the functional $\Phi:C([0,T], \R^d)^2 \to \R$  defined by 
\[
\Phi(x,y)= \sup_{t\in[0,T]}|x(t)-y(t)|.
\]
\begin{Lemma}[Splitting lemma]\label{lem:gap filled}   Assume that the functions $b$ and $\sigma$ are Lipschitz in space uniformly in time, that $\displaystyle\sup_{t\in[0,T]}(|b(t,0)|+\|\sigma(t,0)\|)<+\infty$ and that the kernels $K_i$, $i=1,2$, satisfy $({\cal K}^{int}_{\beta})$ and $({\cal K}^{cont}_{\theta})$ for some $\beta>1$ and $\theta\!\in (0, 1]$ respectively {and, if dealing with the $K$-discrete Euler scheme,} {$(\underline{\cal K}^{cont}_{\underline\theta})$ and  $(\underline{\widehat {\cal K}}^{cont}_{\underline{\widehat{\theta}}})$ for some $\underline\theta,\underline{\widehat{\theta}}\in(0,1]$}. Let $\Phi:C([0,T], \R^d)^2 \to \R$ be a Borel functional and let $n \!\in \N$. Assume  there exists  $\bar p>0$ and a real constant $C>0$ possibly depending on $n$, such that, for every $x_0\!\in \R^d$, 
\[
\| \Phi(X^{x_0}, \bar X^{x_0}) \|_{\bar p} \le C (1+|x_0|),
\]
where $X^{x_0}$ and $\bar X^{x_0}$  respectively denote the solutions of the Volterra equation and any of its two (genuine)  Euler schemes starting from $x_0$. 
Then, for every $p\!\in (0,\bar p]$ and every random vector $X_0$ independent of $W$, the solution $X=(X_t)_{t\in [0,T]}$ and the Euler scheme under consideration starting from $X_0$ satisfy
\[
\| \Phi(X, \bar X) \|_p \le 2^{(1/p-1)^+} C (1+\|X_0\|_p).
\]
\end{Lemma}

\noindent {\bf Proof}. According to Theorem~\ref{lem:flotVolterra} below (whose proof is an adaptation  to Volterra equations of the proof of Blagove$\check{\rm   s}\check{\rm  c}$enkii-Freidlin's theorem~\cite[Theorem 13.1, Section V.12-13, p.136]{RoWI2000}),
there exists  a bi-measurable functional $F: \R^d\times {\cal C}_0([0,T], \R^q)\to {\cal C}([0,T], \R^d)$ such that  the processes $X^{x_0}$ reads
\[
X^{x_0} = F\big(x_0,(W_t)_{t\in [0,T]}\big)
\]
and such that, for any starting random value $X_0\!\in L^p(\P)$, $p>0$,  $F\big(X_0,(W_t)_{t\in [0,T]}\big)$ solves the Volterra equation~\eqref{eq:Volterra} starting from $X_0$.  

The same holds true for (any of) the two Euler schemes (both denoted   $\bar X$ here):  there exists  a measurable functional $\bar F_n: \R^d\times {\cal C}_0([0,T], \R^q)\to {\cal C}([0,T], \R^d)$ such that 
\[
\bar X = \bar F_n\big(X_0,(W_t)_{t\in [0,T]}\big),
\]
{(see Proposition~\ref{lem:K-intF} further on for the $K$-integrated scheme whereas,  for the $K$-discrete Euler scheme, this can be proved by induction in an elementary way, the continuity of this scheme under $(\underline{\cal K}^{cont}_{\underline\theta})$ and  $(\underline{\widehat {\cal K}}^{cont}_{\underline{\widehat{\theta}}})$ for some $\underline\theta,\underline{\widehat{\theta}}>0$, being established in Appendix \ref{app:contconvKdisc} below).}

\smallskip
This entails that the distribution $\P_{(X,\bar X)}$ on ${\cal C}([0,T], \R^d)^2$   of
$(X,\bar X)= \big(F(X_0,W), \bar F_n(X_0,W)\big)$ satisfies
\[
\P_{(X,\bar X)}(dx, d\bar x) = \int_{\R^d} \P_{X_0}(dx_0) \P_{(X^{x_0},\bar X^{x_0})}(dx, d \bar x).
\]
Consequently, using   the monotonicity of probabilistic $L^r(\P)$-norms and pseudo-norms (in the third line) and the elementary inequality $(a+b)^{\rho}\le a^{\rho}+b^{\rho}$, for $a$, $b \ge 0$ and $\rho =\frac{p}{\bar p} \!\in [0,1]$, we derive 
\begin{align*}
\|\Phi(X,\bar X)\|_p^p = \E \,|\Phi(X,\bar X^n)|^p &= \int_{\R^d}\P_{X_0}(dx_0)\E\, |\Phi(X^{x_0},\bar X^{x_0})|^p\\
		              &\le  \int_{\R^d}\P_{X_0}(dx_0)\Big(\E\, |\Phi(X^{x_0},\bar X^{x_0})|^{\bar p}\Big)^{p/\bar p} \\
			        &\le  \int_{\R^d}\P_{X_0}(dx_0)\Big(C^{\bar p}\big(1+|x_0|)^{\bar p}\Big)^{p/\bar p}\\
			        &\le  C^p \int_{\R^d} \P_{X_0}(dx_0)(1+|x_0|^p)   = C^p (1+ \|X_0\|_p^p)\\
			        &\le C^p 2^{(1-p)^+}(1+ \|X_0\|_p)^p
\end{align*} 
so that, finally,
\[
\hskip 5,1cm \|\Phi(X,\bar X)\|_p\le  2^{(1/p-1)^+}C(1+ \|X_0\|_p).\hskip 5,1cm\Box
\]

\smallskip
Our next task is to establish representations for the Volterra process and its two Euler schemes in order to be able to apply the splitting lemma.

\bigskip
The next theorem deals with the flow  $x\mapsto (X^x_t)_{t\in [0,T]}$ of the Volterra equation. It proves the existence of a bi-measurable functional $F: \R^d\times {\cal C}_0([0,T], \R^q)\to {\cal C}([0,T], \R^d) $    such that the solution $X=(X_t)_{t\in [0,T]}$ of equation~\eqref{eq:Volterra} reads $X= F(X_0,W)$ for any (finite) starting value  $X_0$.
\begin{Theorem}[Blagove$\check{\rm \bf s}\check{\rm \bf c}$enkii-Freidlin like theorem: representation  of  Volterra's flow]\label{lem:flotVolterra} Assume that the functions $b$ and $\sigma$ are lipschitz in space uniformly in time, that $\displaystyle\sup_{t\in[0,T]}(|b(t,0)|+\|\sigma(t,0)\|)<+\infty$ and that the kernels $K_i$, $i=1,2$, satisfy $({\cal K}^{int}_{\beta})$ and $({\cal K}^{cont}_{\theta})$ for some $\beta>1$ and $\theta\!\in (0, 1]$ respectively. 

\smallskip
\noindent $(a)$  Let $X^x$ denotes the solution to the Volterra equation~\eqref{eq:Volterra} starting from $x\!\in \R^d$ and let $\lambda \!\in (\frac 12, 1)$. There exists {$p^*=p^*_{\beta, \theta, \lambda,d}$ (made explicit in the proof)}  
such that for every $p> p^*$,   
  \[
\forall \,x,y\in\R^d,\;\Big\| \sup_{t\in [0,T]} |X^x_t-X^y_t| \Big\|_p \le C |x-y|^{\lambda},
\]
for some positive real constant $C= C_{p, b,\sigma, K_1,K_2,\beta, \theta}$.

\smallskip
\noindent $(b)$  There exists a functional $F: \R^d\times {\cal C}_0([0,T], \R^q)\ni(x,w) \mapsto F(x,w) \in {\cal C}([0,T], \R^d)$ bi-measurable (where the spaces of continuous functions are equipped with the Borel sigma-field induced by the uniform convergence topology) and continuous in $x$ such that, for any stochastic basis $(\Omega, {\cal A}, \P, (\F_t)_{t\in [0,T]})$, any $q$-dimensional $(\F_t)_t$-Brownian motion $W$ and any $\F_0$-measurable $\R^d$-valued random vector $X_0\!\in L^0(\P)$, the solution to the Volterra equation~\eqref{eq:Volterra} is $X=F(X_0,W)$. 
\label{thm:BF}\end{Theorem}

\noindent {\bf Proof}. The proofs of claims~$(a)$ and $(b)$ are intertwined.\\{\sc Step~1} ({\em Continuity of the flow}).  It follows from~\cite[Theorem 3.3]{ZhangXi2010} that,  {for every $p> p_{eu}=\frac{1}{\theta}\vee \frac{2\beta}{\beta-1}$}, the Volterra equation~\eqref{eq:Volterra} has a unique strong solution starting from any random vector $X_0\!\in L^p(\P)$ which  can be proved to be pathwise continuous (and even $a$-H\"older pathwise continuous for some small enough $a>0$) and satisfying furthermore  $\E \sup_{t\in [0,T]} |X_t|^p <+\infty$.

{By mimicking the proof   of Property~{\bf 4} in 
Appendix~\ref{app:B}  to control  $\|X_t-\bar X_t\|_p$ (see also the proof of~Theorem~3.3 in~\cite{ZhangXi2010}), one derives  that there exists a positive real constant $C_{p,T}= 2\,e^{C^{(8)}T}$ such that,
\begin{equation}\label{eq:LipMarginal}
\forall\,p>p_{eu}, \quad \sup_{t\in [0,T]}\| X^x_t-X^y_t\big\|_p \le C_{p,T} |x-y|.
\end{equation}
}
{ We also know from  Property~{\bf 3} in Appendix~\ref{app:B}  (see also~\cite{ZhangXi2010}) that for every $p\ge 2$,  there exists a positive real constant $C'_{p,T}$ (not depending on $n$) such that  the $K$-integrated Euler scheme with step $\frac Tn$ satisfies 
\[
\forall\, s,\, t\!\in [0,T], \, \forall\, x\!\in \R^d, \quad  \| \bar X^x_t-\bar X^x_s\big\|_p \le C'_{p,T}(1+|x|)|t-s|^{\tilde \theta} 
\]
with $\tilde \theta = \theta \wedge \frac{\beta-1}{2\beta}$. On the other hand, it follows from Property~{\bf 4}  in  Appendix~\ref{app:B}
 that,  for every $t\!\in [0,T]$,  $\bar X_t\to X_t$ as $n\to +\infty$  in $L^p(\P)$ if $p>p_{eu}$. Hence, for every $p>p_{eu}$,
\[
\forall\, s,\, t\!\in [0,T], \, \forall\, x\!\in \R^d, \quad  \| X^x_t-X^x_s\big\|_p \le C'_{p,T}(1+|x|)|t-s|^{ \tilde \theta}.
\]
}
Now, let  $\lambda \!\in \big(0,1\big)$ be fixed. One has 
\begin{align*}
\big\| X^x_t-X^y_t-(X^x_s-X^y_s)\big\|_p &\le \big( \|X^x_t-X^y_t \|_p+\|X^x_s-X^y_s\|_p\big)^{\lambda} \big( \|X^x_t-X^x_s\|_p + \|X^y_t-X^y_s\|_p \big)^{1-\lambda}\\
& \le 2  C_{p,T}^{\lambda}(C'_{p,T})^{1-\lambda}|x-y|^{\lambda}(1+|x|+|y|)^{1-\lambda}|t-s|^{ \tilde \theta(1-\lambda)}
\end{align*}
or, equivalently, 
\[
\forall\, s,t\!\in [0,T], \, \forall\, x\!\in \R^d, \quad \E | X^x_t-X^y_t-(X^x_s-X^y_s)|^p \le  \tilde C_{p,T}(1+|x|+|y|)^{p(1-\lambda)}|x-y|^{p\lambda}|t-s|^{p \tilde \theta(1-\lambda)}.
\]

Assume from now on that $p>p^*= p_{eu}\vee \frac{1}{ \tilde \theta(1-\lambda)}\vee \frac{d}{\lambda}$, where the second bound from below (by $\frac{1}{ \tilde \theta(1-\lambda)}$) will permit to apply Kolmogorov's criterion in time and the third to apply this criterion in space. It follows from the proof of Kolmogorov's criterion~\cite[Theorem 2.1]{RevuzYor} where tracking the constants ensures that the final constant has linear growth in the constant which appears in the hypothesis, that, for any $a\!\in\big(0,(1-\lambda) \tilde \theta-\frac 1p\big)$, there exists a positive real constant  $C_{a,p,T}$ such that 
\[
 \forall\, x,y\!\in \R^d, \quad \E \sup_{s,t\in [0,T]} \left(\frac{ | X^x_t-X^y_t-(X^x_s-X^y_s)|}{|t-s|^{a}}\right)^p \le   C_{a,p,T}(1+|x|+|y|)^{p(1-\lambda)} |x-y|^{p\lambda}.
\]
Setting  $s=0$ yields
\begin{align*}
 \forall\, x,\, y\!\in \R^d, \quad \E\sup_{t\in [0,T]} \big |X^x_t-X^y_t\big|^p& \le  2^{p-1}\left(|x-y|^p+ T^a C_{a,p,T}(1+|x|+|y|)^{p(1-\lambda)} |x-y|^{p\lambda}\right)\\
 & \le C'_{a,p,T} (1+|x|+|y|)^{p(1-\lambda)} |x-y|^{p\lambda}
\end{align*}
where we used the rough inequality $|x-y|^{p(1-\lambda)} \le (1+|x|+|y|)^{p(1-\lambda)}$.

\smallskip  Let $R\!\in \N$ be fixed. For every $x\!\in B(0,R)$, we consider a pathwise continuous modification $(X^x_t)_{t\in [0,T]}$ so that the mapping $B(0,R) \ni x\mapsto (X^x_t)_{t\in [0,T]}$ can be viewed as having values in the Polish space ${\cal C}([0,T], \R^d)$, equipped with the distance $\rho_{_T}$ derived from the  $\sup$-norm on $[0,T]$. 
Then for every $x, y\!\in B(0,R)$
\begin{align}\label{eq:flowsup}  
 \forall\, x,\, y\!\in B(0,R), \quad \E\sup_{t\in [0,T]} \big |X^x_t-X^y_t\big|^p
& \le C'_{a,p,T} (1+2R)^{p(1-\lambda)}|x-y|^{p\lambda}.
\end{align}
As $p\lambda >d$, a new application of Kolmogorov's criterion~\cite[Theorem 2.1]{RevuzYor} ensures the existence of  a $\P$-modification $(\widetilde X^{(R),x})_{x\in B(0,R)}$ of the ${\cal C}([0,T], \R^d)$-valued process $(X^x)_{x\in B(0,R)}$ such that  $\P$-$a.s.$ $x\mapsto \widetilde X^{(R),x}$  is $\rho_{T}$-continuous and more precisely  $a'$-H\"older-continuous for $a'\!\in \big( 0, \lambda-\frac dp
\big)$.

It is clear that the restriction of $(\widetilde X^{(R+1),x})_{x\in B(0,R+1)}$ to $B(0,R)$ is also a $\P$-modification with  continuous ``paths'' in $x$. Hence they are $\P$-indistinguishable. Consequently, there exists a $\P$-modification $(\widetilde X^x)_{x\in \R^d}$ of $(X^x)_{x\in \R^d}$ such that $\P(d\omega)$-a.s. $x\mapsto (\widetilde X^x_t)_{t\in [0,T]}$ is continuous from $\R^d$ to ${\cal C}([0,T], \R^d)$ and $(\widetilde X^x_t)_{t\in [0,T]}$ is solution to~\eqref{eq:Volterra} for $X_0=x$. 

\smallskip
\noindent  {\sc Step~2} ({\em The functional $F$)}. Now we adapt to Volterra equations the classical proof of Blagove$\check{\rm s}\check{\rm c}$enkii-Freidlin's theorem (see e.g.~\cite[Theorem 10.4]{RogersWilliamsII}) originally written for standard Brownian diffusion processes. We  only take advantage of the strong existence and uniqueness of pathwise continuous  strong solutions for every   $\mu$-distributed starting random vector $X_0$  having a finite $L^p$-moment to exhibit  a functional $\widetilde F_{\mu}: \R^d\times   {\cal C}_0([0,T], \R^q)\to  {\cal C}([0,T], \R^d)$ (adapted w.r.t. to the canonical filtrations of both spaces) such that, for any strong   solution to~\eqref{eq:Volterra}  driven by a   $q$-dimensional $(\F_t)$-Brownian motion $B=(B_t)_{t\in [0,T]}$  on a filtered probability space $(\Omega, {\cal A}, (\F_t), \P)$ and any $\F_0$-measurable starting  random vector $X_0$ with distribution $\mu$:
\[
\P\mbox{-} a.s. \quad X= \widetilde F_{\mu}(X_0,B).
\]
Let us be more precise. We consider the canonical space $ \widetilde \Omega_0 = \R^d\times{\cal C}_0([0,T], \R^q)$ equipped with the product $\sigma$-field 
{${\cal B}or(\R^d)\otimes {\cal B}or({\cal C}_{q,0})$} and the product probability measure $\P^0= \mu\otimes \Q_{_W}$ where $\mu$ has a finite $p$-th moment for some $p>p_{eu}$ and $\Q_{_W}$ denotes the $q$-dimensional Wiener measure.  Let $(X_0,B)$ denote the canonical process on $ \widetilde \Omega_0$ defined by $(X_0, B)(x,w)= (x,w)$. As $p> p_{eu}$ the Volterra 
SDE~\eqref{eq:Volterra} has a unique pathwise strong solution  $(\xi_t)_{t\in [0,T]}$ defined on  $ \widetilde \Omega_0$ (and adapted to the $\P^0$-completed natural filtration  $\F^0_t = \sigma({\cal N}_{\P_0}, X_0, B_s,\,0\le s\le t), t\!\in [0,T]$, etc). We set 
$$
\widetilde F_\mu(x,w) = \xi(x,w).
$$
 Then one remarks that on any filtered probability space $(\Omega, {\cal A}, ({\cal F}_t) , \P)$ on which  lives a $q$-dimensional $(\F_t)$-Brownian motion $W$ and $X_0$ an $\F_0$-measurable $\R^d$-valued random vector, then $\widetilde F_\mu(X_0, W)$ is a strong solution to the Volterra equation related to $W$ and by strong uniqueness it is $\P$-indistinguishable of the ``natural" strong solution $X$ starting from $X_0$.

Thus, the  functional $\widetilde F_{\mu}$ being defined on the canonical space $\widetilde \Omega_0$, does not depend on the selected stochastic basis (nor the Brownian motion under consideration). In particular, for every $x\!\in \R^d$, setting$$
F(x,w):= \widetilde F_{\delta_x}(x,w),\;(x,w)\!\in \R^d\times{\cal C}_0([0,T], \R^q),
$$ one  has $\widetilde X^x = F (x,W)$ $\P$-a.s. (by uniqueness of strong solutions of the  Volterra equation since the starting value $x$ is deterministic).

If  $\Q_{_W}$ denotes the Wiener measure on ${\cal C}([0,T], \R^q)$
\[
\int \,\rho_{_T}\big(F(x,w), F(x', w))^p \Q_{_W}(dw) = \E\, \rho_{_T}(\widetilde X^x,\widetilde X^y)^p
\]
so that, we show  like for $\widetilde X^x$, that $x\mapsto F(x,w)$ admits an $a'$-H\"older  $\Q_{_W}$-modification since $p>p^*$ owing to Step~1. Hence $x\mapsto F(x,W)$ and $(\widetilde X^x)_{x\in \R^d}$ are $\P$-indistinguishable. Note that since $F$ is continuous in $x$ and measurable in $w$, this functional is bi-measurable according to Lemma 4.51~\cite{AliBord}.

\smallskip
\noindent  {\sc Step~3} ({\em Representation for $X_0\!\in L^p(\P)$, $p$ large}).  First let us consider the case where $X$ is the pathwise continuous solution to the Volterra equation starting from a random vector $X_0\in L^p(\P)$, $p>p^*$, independent of $W$ given by~\cite[Theorem 3.3]{ZhangXi2010}. We will prove that $X=F(X_0,W)$.
Mi\-micking the proof to establish~\eqref{eq:LipMarginal}, one shows that, if $X'$ denotes the pathwise continuous solution to the Volterra equation starting from $X'_0\in L^p(\P)$ independent of $W$, 
\begin{equation}\label{eq:flotgeneral}
\sup_{t\in [0,T]} \|X_t-X'_t\|_p \le C_{p,T} \|X_0-X'_0\|_p.
\end{equation}
Let $\varphi_k:\R^d\to \mathbb{Q}^d$ such that $|\varphi_k(x)|\le |x|$  for every $x\!\in \R^d$ and $\sup_{x \in \R^d}|\varphi_k(x)-x|\to 0$ as $k\to+\infty$ and set $X^{(k)}_0=\varphi_k(X_0)$. The set $N=  \cup_{x\in \mathbb{Q}^d}\{\omega\!\in \Omega: \widetilde X^x(\omega)\neq F (x,W(\omega)) \}\cup N_0$ is $\P$-negligible by construction. Since 
$X^{(k)}_0$ is a $\mathbb{Q}^d$-valued, $\F_0$-measurable random vector, independent of $W$, $F(X^{(k)}_0,W)$ is a pathwise continuous solution to the Volterra equation so that, by~\eqref{eq:flotgeneral}, 
$$ 
\sup_{t\in [0,T]} \|X_t-F(X_0^{(k)},W)_t\|_p\le C_{p,T} \|X_0-X^{(k)}_0\|_p\to 0\mbox{ as }k\to +\infty.
$$
With the triangle inequality, one deduces that
\begin{align*}
\limsup_{k\to\infty}\sup_{t\in [0,T]} \| X_t-F(X_0,W)_t\|_p &\le \limsup_{k\to\infty}\big\|\rho_{_T}\big(F(X_0,W),F(X_0^{(k)},W)\big)\big\|_p.
\end{align*}
Using that $X_0$ and $W$ are independent   and the bound~\eqref{eq:flowsup}, we get 
\begin{align*}
\E\, \rho_{_T}\big(F(X_0,W),F(X_0^{(k)},W)\big)^p& = \int_{\R^d} \P_{X_0}(dx_0) \E \,\rho_{_T}\big(F(x_0,W),F(\varphi_k(x_0),W)\big)^p\\
& = \int_{\R^d}\P_{X_0}(dx_0) \E \,\rho_{_T}\big(X^{x_0},X^{\varphi_k(x_0)}\big)^p\\
&\le C'_{a,p,T}\sum_{R\in\N^*}(1+2R)^{p(1-\lambda)}\E \Big[|X_0-\varphi_k(X_0)|^{\lambda p} \mbox{\bf 1}_{\{R-1\le |X_0|\le R\}}\Big]\\
&\le C'_{a,p,T}\Big( 3^{p(1-\lambda)}\E \Big[|X_0-\varphi_k(X_0)|^{\lambda p} \mbox{\bf 1}_{\{ |X_0|\le 1\}}\Big]\\
&\qquad+ \sum_{R\ge 2}(1+2R)^{p(1-\lambda)}\E \Big[|X_0-\varphi_k(X_0)|^{\lambda p}\mbox{\bf 1}_{\{R-1\le |X_0|\le R\}}\Big]\Big)\\  
&\le  C_{a,p,T,X_0}\Big(1+\sum_{R\ge 2}\frac{(1+2R)^{p(1-\lambda)}}{(R-1)^p}\Big)\E \Big[|X_0-\varphi_k(X_0)|^{\lambda p} |X_0|^p \Big]
\end{align*}
where we used Markov inequality in the last line.
The series is converging since  $\frac{(1+2R)^{p(1-\lambda)}}{(R-1)^p}\sim \frac{1}{R^{p\lambda}}$ and $p\lambda>d\ge 1$. Hence
\[
\E\, \rho_{_T}\big(F(X_0,W),F(X_0^{(k)},W)\big)^p\le C'_{a,p,T,X_0}  \E \Big[|X_0-\varphi_k(X_0)|^{\lambda p} |X_0|^p \Big]\le C'_{a,p,T,X_0} \sup_{x \in \R^d}|\varphi_k(x)-x|\E\, |X_0|^p\to 0 
\]
as $k\to +\infty$. Consequently, for every $t\!\in [0,T]$, $X_t = F(X_0,W)_t$ $\P$-$a.s.$. Both being pathwise continuous $X= F(X_0,W)$ $\P$-a.s.

\smallskip
\noindent  {\sc Step~4} ({\em Representation for $X_0\!\in L^0(\P)$}). Let $X_0\!\in L^0(\P)$.   We consider like in Appendix~\ref{app:eu} the sequence of  truncated starting values $X^{(k)}_0= X_0\mbox{\bf 1}_{A_k}$ with $A_k=\{|X_0|< k\}$, $k\ge 1$ and the resulting pathwise continuous solution to~\eqref{eq:Volterra} still denoted $X^{(k)}$. As $X^{(k)}_0\!\in L^p(\P)$ for $p>p^*$, then $\P$-$a.s.$ 
\[
X^{(k)}= F(X_0^{(k)},W),\quad k\ge 1.
\]
Hence, under the convention $A_0=\emptyset$,
\begin{align*}
   \hskip3.25cm X=\sum_{k\ge 1}X^{(k)}\mbox{\bf 1}_{A_k\setminus A_{k-1}}=\sum_{k\ge 1}F(X_0^{(k)},W)\mbox{\bf 1}_{A_k\setminus A_{k-1}}=F(X_0,W).\hskip3.25cm\Box
\end{align*}

{\begin{Proposition}[Representation formula for  the $K$-integrated Euler scheme] \label{lem:K-intF} Assume $({\cal K}^{int}_{\beta})$ and $({\cal K}^{cont}_{\theta})$ are in force. Let $n\ge 1$. There exists a bi-measurable functional $\bar F_n : \R^d\times {\cal C}_0([0,T], \R^q)\to  {\cal C}([0,T], \R^d)$ such that, for any stochastic basis $(\Omega, {\cal A}, \P, (\F_t)_{t\in [0,T]})$, any $q$-dimensional $(\F_t)_t$-Brownian motion $W$ and any $\F_0$-measurable $\R^d$-valued random vector $X_0\!\in L^0(\P)$,  the $K$-integrated Euler scheme $\bar X$ with step $\frac Tn$ starting from $X_0$ defined by~\eqref{eq:Eulergen2} writes
\[
\bar X = \bar F_n (X_0,W).
\]
\end{Proposition}
}
\noindent {\bf Proof.} As for the $K$-integrated scheme, it amounts to prove that, for every $k\in\{0,\cdots,n-1\}$, the process
\[
 [t^n_k,T]     \ni t\mapsto \int_{t^n_k}^{t\wedge t^n_{k+1}} K_2(t,u)dW_u
\]
can be written as a functional of $W$. First we can assume temporarily  that $W$ is defined on the canonical space $\Omega_0= {\cal C}_0([0,T], \R^q)$~--~ i.e. $W_t(w)= w(t)$, for every $w\!\in \Omega_0$~--~equipped with its Borel $\sigma$-field {${\cal B}or({\cal C}_{q,0})$} induced by the $\sup$-norm topology and the Wiener measure $\Q_{_W}$. For every $t\!\in [t^n_k,T]$,
 $\displaystyle  \int_{t^n_k}^{t\wedge t^n_{k+1}} K_2(t,u)dW_u$ is $\sigma\big({\cal N}_{\P},W_{u}-W_{t^n_k},\, t^n_k\le u\le t\wedge t^n_{k+1}\big)$-measurable hence $\sigma\big({\cal N}_{\P},W_u,\, 0\le u\le T\big)$-measurable  so that 
\begin{equation}\label{eq:Wienerrepresentation}  \int_{t^n_k}^{t\wedge t^n_{k+1}} K_2(t,u)dW_u = F^{k}(t,W).
\end{equation}
(A more precise result involving adaptation to the canonical filtration can be obtained likewise but it is not useful for our purpose here.) 
The pathwise continuity of $\left(\int_{t^n_k}^{t\wedge t^n_{k+1}} K_2(t,u)dW_u\right)_{t\in[t^n_k,T]}$ follows from Lemma \ref{lem:intHcont} applied with $\widetilde H\equiv 1$ and implies that,  for  $\Q_{_W}(dw)$-almost every $w\!\in  \Omega_0$,  $[t^n_k,T]\ni t\mapsto F^{k}(t,w)\!\in {\cal C}_0([0,T],\R^d)$.
Then for any standard Brownian motion $W$ on a stochastic basis $(  \Omega, {\cal A}, (\F_t)_{t\in [0,T]}, \P)$ the representation~\eqref{eq:Wienerrepresentation} holds true (with the same functional) since the Gaussian law of $\left(W,\int_{t^n_k}^{t\wedge t^n_{k+1}} K_2(t,u)dW_u\right)$ does not depend on the choice of $W$.    

 This being done, one checks by induction on the discretization times  using~\eqref{eq:Eulergen2bis} that  the $K$-integrated Euler scheme starting from any  (finite) $\R^d$-valued starting random vector $X_0$ admits a representation 
 \[
 \bar X= \bar F_n (X_0,W),
 \]
where {$\bar F_n : \big(\R^d\times \Omega_0, {\cal B}or(\R^d)\otimes {\cal B}or({\cal C}_{q,0})\big) \to \big({{\cal C}([0,T],\R^d), {\cal B}or({\cal C}_{d})}\big)$ is  measurable.}

\section{End of the proof of Theorem~\ref{prop:exunvolt}}\label{app:A''}
\smallskip
\noindent {\sc Step~1} ({\em $a$-H\" older version/modification when $X_0\!\in L^p(\P)$, $p$ large enough}). Our bound~\eqref{eq:Holderpaths} is more general and more precise  than its counterpart in~\cite[Theorem~3.3]{ZhangXi2010}  since it takes the form 
\begin{equation}\label{eq:improved}
\left \| \sup_{s\neq t\in [0,T]}\frac{|X_t-X_s|}{|t-s|^{a}}\right\|_p \le C_{a,p,T} (1+\|X_0\|_p),
\end{equation}
 where $C_{a,p,T}$ does not depend on $X_0$, $a$ is any element of $(0, \theta\wedge \frac{\beta-1}{2\beta})$ and only $X_0\!\in L^p(\P)$ is requested.  This turns out to be the key for the extension to lower exponents $p$ in what follows.

At this stage we no longer rely on~\cite{ZhangXi2010} whose proof in Theorem~3.1 is sub-optimal. We could provide a direct proof to evaluate $\|X_t-X_s\|_p$, first  for $p\ge 2$ {by just mimicking that of Property~{\bf 3} in Appendix~\ref{app:B}}. In fact it follows from Properties {\bf 3} and {\bf 4} from  Appendix~\ref{app:B} the following two results: first, for  $p\ge 2$ there exists a positive real constant  $C$ (not depending on $n$) such that for every $n\ge 1$ the genuine $K$-integrated Euler scheme $\bar X$ with step $\frac Tn$ satisfies  the H\"older regularity property
\[
\forall\, s,\, t\in [0,T], \quad \|\bar X_t-\bar X_s\|_p\le C(1+\|X_0\|_p)|t-s|^{\theta \wedge \widehat \theta}
\]
and, for $p>p_{eu}=\frac 1\theta\vee\frac{2\beta}{\beta-1}>2$, the ``marginal'' convergence property of the $K$-integrated Euler scheme toward $X$ holds i.e. $\|X_t-\bar X_t\|_p\to 0$ as $n\to +\infty$ for every $t\!\in [0,T]$.
As a consequence, letting $n\to 0$, yields that, for $p>p_{eu}$,  
\begin{equation}\label{eq:improvedbis}
\forall\, s,\, t\in [0,T], \quad \| X_t- X_s\|_p\le C(1+\|X_0\|_p)|t-s|^{\theta \wedge \widehat \theta}.
\end{equation}

At this stage, we can proceed as follows: let   $a\!\in (0,\theta  \wedge \widehat \theta)$  and let    $p>\frac{1}{\theta \wedge \widehat \theta-a} \vee p_{eu} $ so that $a<\theta\wedge \widehat\theta -\frac 1p$. Tracking once again the constants   in the proof of Kolmogorov's criterion, e.g. in~\cite[Chapter 2, p.26]{RevuzYor}, to check that the final constant  depends linearly in the initial one,
one derives that~\eqref{eq:improved} does hold true for such $p$.

\smallskip 

\smallskip
\noindent {\sc Step~2} ({\em $L^p$-pathwise regularity  $p\!\in (0, +\infty)$}). To extend~\eqref{eq:improved} to every $p>0$ such that $\|X_0\|_p<+\infty$, one   can take advantage of the form of the control 
 to apply  the splitting Lemma~\ref{lem:gap filled}  with  $p<\bar p$, $\bar p>\frac{1}{\theta   \wedge \widehat \theta-a}\vee p_{eu}$  and the functional $\Phi: {\cal C}([0,T], \R^d)^2\to \R$ defined by   
\[
  \Phi(x, y) = \sup_{s\neq t \in [0,T]} \frac{|x(t)-x(s)|}{|t-s|^{a} }. 
\]
The marginal bound~\eqref{eq:improvedbis} can be  extended likewise by considering    $\Phi(x,y) = |x(t)-x(s)|$. \hfill $\Box$

\section{Proofs of Theorems~\ref{thm:Eulercvgce2}  and~\ref{thm:Eulercvgce1}} \label{app:E}

\subsection{Proof of Theorem~\ref{thm:Eulercvgce2} (convergence of the $K$-integrated Euler scheme)}
In this section, we complete the proof of Theorem~\ref{thm:Eulercvgce2}, based on the splitting Lemma 
 established in Appendix~\ref{app:B'}.

 \medskip
 \noindent {\sc Step~1} ({\em Moment control, $p\!\in (0,2)$}). Let us denote 
 $\bar X^{x_0}= (\bar X^{x_0}_t)_{t\in [0,T]}$ the 
 Euler scheme~\eqref{eq:Eulergen2}  starting from $x_0$ at $t=0$.
 Property~{\bf 2} in Appendix \ref{app:B} applied with $p=2$ 
 implies that there exists a real constant $C=C_{K_1, K_2, \beta, b,\sigma, T}>0$ (not depending on the time step $\frac Tn$ i.e. on $n$) such that,  if $X_0= x_0$, 
 \[
\sup_{t\in [0,T]}\|\bar X^{x_0}_t\|_2  
\le C (1 +|x_0|).
 \]
Then it follows from Lemma~\ref{lem:gap filled} applied with $\bar p=2$ and $\Phi(x,y) =\sup_{t\in[0,T]}|y(t)|$ that
\[
\forall\, p\!\in (0,2], \quad\sup_{t\in [0,T]}\|\bar X_t\|_p 
 \le 2^{(1/p-1)^+}C (1+\|X_0\|_p).
 \]

\medskip 
 \noindent {\sc Step~2} ({\em  $L^p$-control of the increments of the $K$-integrated Euler scheme for $p\!\in (0,2)$}). 
  If {$p\!\in (0,2)$}, we proceed likewise replacing Property~{\bf 2} in Appendix \ref{app:B} by Property{\bf 3}, with the functional $\Phi(x,y) = |y(t)-y(s)|$, $s,\, t\!\in [0,T]$ in order  to prove that
  \[
\forall\, s, t\!\in [0,T], \quad  \|\bar X_t-\bar X_s\|_p\le  2^{(1/p-1)^+}C_{2,T}(1+\|X_0\|_p) |t-s|^{\theta\wedge\widehat\theta}.
\]
 
 \noindent {\sc Step~3} ({\em $L^p$-rate of convergence for the marginals at fixed time $t$ {for $p>0$}}).  
We fix $\bar p>p_{eu}$ and combine Property~{\bf 4} in Appendix \ref{app:B} and Lemma~\ref{lem:gap filled} applied with  $\Phi(x,y) =  |x(t)-y(t)|$ to $p\in (0,\bar p]$.

\smallskip
 \noindent{\sc Step~4} ({\em $L^p$-rate of convergence for the sup-norm, $p$ large enough}). To switch from the $L^p$-convergence of marginals, namely $\|X_t-\bar  X_t\|_p$
  to this one, we rely on Corollary 4.4~\cite{RiTaYa2020} which is a form close to Kolmogorov's $C$-tightness criterion of the Garsia-Rodemich-Rumsey lemma (see~\cite{GRR1070}).
\begin{Theorem}[GRR lemma] Let $(Y^n)_{n \ge 1}$ be a sequence of continuous processes where the processes $Y^n= (Y^n_t)_{t\in[0,T]}$ are defined on a probability space $(\Omega,{\cal A}, \P)$. Let $p\ge 1$. Assume  there exists  $a>1$, a sequence $(\delta_n)_{n\ge 1}$ of positive real numbers  converging to $0$ and  a real constant $C>0$ such that 
\begin{equation}\label{eq:GRR}
\forall\, n\ge 1,\; \forall\, s, t\!\in [0,T], \quad \E\, |Y^n_t-Y^n_s|^p \le C|t-s|^{a} \delta^p_n.
\end{equation}
Then
there exists a real constant $C_{p,T}>0$ such that
\[
\forall\, n\ge 1,\quad \E\, \sup_{t\in [0,T]}|Y^n_t-Y^n_0|^{p} \le C_{p,T}\,\delta^p_n.
\]  
\end{Theorem}
 
 Let $\varepsilon\!\in (0,1)$. We aim at applying the above inequality to the sequence of processes $Y^n = X-\bar X$ to prove that the order of convergence of the Euler scheme derived at Property~{\bf 4} from Appendix~\ref{app:B} is preserved for the supremum over time up to multiplication by the factor $1-\varepsilon$.

\smallskip
Set $\theta^{*} =\theta\wedge\widehat\theta $. We first deal with large values of $p$, namely {$p>p_{\varepsilon} =\frac{2}{\theta^{*}\varepsilon}\vee p_{eu}$. Then let $\lambda := \frac{1}{p\theta^{*}}+\frac{\varepsilon}{2}\!\in (0,\varepsilon)$}. For every $s$, $t\!\in [0,T]$, one has
 \[
\|Y^n_t-Y^n_s\|_p \le \|Y^n_t-Y^n_s\|_p^{\lambda}\|Y^n_t-Y^n_s\|_p^{1-\lambda} \le  \big(\|X_t-X_s\|_p +\|\bar X_t-\bar X_s\|_p  \big)^{\lambda} \big(\|X_t-\bar X_t\|_p+\|X_s-\bar X_s \|_p \big)^{1-\lambda} . 
\]
By Property~{\bf 4} from Appendix~\ref{app:B}, there exists a positive real constant $C_1$ such that 
\[
\sup_{t\in [0,T]} \|X_t-\bar X_t\|_p\le C_1  (1+\|X_0\|_p) \big( \tfrac Tn \big)^{\theta^{*}\wedge \gamma}.
\]
With Property~{\bf 3}, still from Appendix~\ref{app:B}, we deduce the existence of a positive real constant $C_2$ such that 
\[
\forall\, s,\, t\!\in [0,T], \quad \| X_t-X_s\|_p +\|\bar X_t-\bar X_s \|_p \le C_2 (1+\|X_0\|_p)|t-s|^{\theta^{*}}.
\]
Hence, 
using that   $a:=\lambda \theta^{*} p=1+ \frac{\varepsilon}{2}\,\theta^{*}\,p >1$,  we derive
\[
\E\, |Y^n_t-Y^n_s|^p \le C^p_3  (1+\|X_0\|_p)^p|t-s|^{a} \big( \tfrac Tn \big)^{p(\theta^{*}\wedge \gamma)(1-\lambda)} 
\]
with  $C_3= C_1^{1-\lambda}C_2^{\lambda}$. Hence, as $Y^n_0=0$ and $a>1$, it follows from the above {\em GRR} lemma that there exists 
a positive real constant $\kappa_{\varepsilon,p} =\kappa_{\varepsilon,p,\beta, \gamma, \theta, b,\sigma,T}$ such that, for every $n\ge 1$, 
\[
\Big\| \sup_{t\in [0,T]} |X_t-\bar X_t|\Big\|_{p} = \Big\| \sup_{t\in [0,T]} |Y^n_t-Y^n_0| \Big\|_{p}  \le \kappa_{\varepsilon,p}  (1+\|X_0\|_p) \big( \tfrac Tn \big)^{(\theta^{*}\wedge \gamma)(1-\lambda)} .
\]
As $\lambda \le \varepsilon$, it is clear that, up to an updating of the constant $\kappa_{\varepsilon,p}$ to take into account the values of $n\!\in\{1,\ldots,\lfloor T\rfloor\}$, one has for every $n\ge 1$, 
\[
\Big\| \sup_{t\in [0,T]} |X_t-\bar X_t|\Big\|_{p}    \le \kappa_{\varepsilon,p}  (1+\|X_0\|_p) \big( \tfrac Tn \big)^{(\theta^{*}\wedge \gamma)(1-\varepsilon)} .
\]
%
Now we take advantage of Appendix~\ref{app:B'}  to extend this result to the low integrability setting of $X_0$, i.e. $X_0\!\in L^p(\P)$, $p\in (0,p_\varepsilon]$. 

 \medskip
 \noindent {{\sc Step~5} ({\em $L^p$-rate of convergence of the supremum over time, $p\!\in(0, p_{\varepsilon})$})}. It follows from Lemma~\ref{lem:gap filled} {applied with $\bar p= p_{\varepsilon}+1$ and 
 $\Phi(x,y)= \sup_{t\in [0,T]}|x(t)-y(t)|$  that 
  \begin{align*}
  \Big\| \sup_{t\in [0,T]} |X_t-\bar X_t|\Big\|_p 
  &\le 2^{(1/p-1)^+} (\kappa_{\varepsilon,p_{\varepsilon}+1}) \Big(\frac Tn\Big)^{(\gamma\wedge \theta^{*})(1-\varepsilon)}(1+\|X_0\|_p).
  \end{align*}}
  
This  completes the proof of this step and of the theorem.~\hfill $\Box$

\subsection{Proof of Theorem~\ref{thm:Eulercvgce1} (continuity and convergence of the $K$-discrete time Euler scheme)}\label{app:contconvKdisc}
Since $(\underline{{\cal K}}^{cont}_{\underline{\theta}})$ holds for some $\underline{\theta}>0$, then for each $\ell\in\{0,\cdots,n-1\}$, $t\mapsto (K_1(t,t^n_\ell),K_2(t,t^n_\ell))$ is locally H\"older continuous with exponent $\underline{\theta}$ on $(t^n_\ell,T]$. For $k\in\{0,\cdots,n-1\}$, since the $K$-discrete Euler scheme reads on every interval $(t^n_k,t^n_{k+1}]$
\begin{align}
   \forall\, t\!\in(t^n_k,t^n_{k+1}],\quad\bar X_t=X_0&+\sum_{\ell=0}^{k-1} \left(K_1(t,t^n_\ell)b(t^n_\ell,\bar X_{t^n_\ell})\tfrac{T}{n}+K_2(t,t^n_\ell)\sigma(t^n_\ell,\bar X_{t^n_\ell})(W_{t^n_{\ell+1}}-W_{t^n_\ell})\right)\notag\\&+K_1(t,t^n_k)b(t^n_k,\bar X_{t^n_k})(t-t^n_k)+K_2(t,t^n_k)\sigma(t^n_k,\bar X_{t^n_k})(W_t-W_{t^n_k}),\label{Kdischemek}
\end{align}
one easily deduces that it is continuous on the interval $(t^n_k,t^n_{k+1}]$ with the sum of the two first terms in the right-hand side even continuous on the closed interval $[t^n_k,t^n_{k+1}]$. Since $(\underline{\widehat{\cal K}}^{cont}_{\underline{\widehat\theta}})$ holds for some $\underline{\widehat\theta}>0$, the choice $\delta=t-t^n_k$ in this condition ensures that $\forall t\in(t^n_k,t^n_{k+1}]$, $K_i(t,t^n_k)\le \widehat\kappa (t-t^n_k)^{\underline{\widehat\theta}-1/i}$. Using the law of the iterated logarithm satisfied by the Brownian motion, we deduce that the sum of the two last terms in the right-hand side of~\eqref{Kdischemek} goes to $0$ as $t\to t^n_k+$ and the $K$-discrete Euler scheme has continuous paths. 
Under the assumptions, it is established in Theorem~\ref{prop:exunvolt} that $X$ is itself pathwise continuous. The authors in~\cite{RiTaYa2020} claim the conclusions of the theorem are  valid for $p$ large enough as  far as~\eqref{eq:incrementXbar2} and~\eqref{eq:Xt-barX_t2} are concerned and  whenever $X_0\!\in \cap_{r>0}L^r(\P)$ for~\eqref{eq:sup|Xt-barX_t|2}. It is straightforward to show they hold true for any $p\!\in (0, +\infty)$ by  applying Lemma~\ref{lem:gap filled}
with the functionals $\Phi(x,y)= |y(t)-y(s)|$, $\Phi(x,y)= |x(t)-y(t)|$ and $\Phi(x,y)= \sup_{t\in [0,T]}|x(t)-y(t)|$, $x,y\!\in {\cal C}([0,T], \R^d)$  respectively. The proof is quite similar to that of Theorem~\ref{thm:Eulercvgce2} at this stage. \hfill$\Box$ 

\end{document}